\magnification=1200
\baselineskip=14pt

\hsize=118mm
\hoffset=7mm
\vsize=185mm
\voffset=10mm

\def\Sc{{\cal S}}
\def\Xc{{\cal X}}

\def\fl{\forall}
\def\ify{\infty}
\def\lgl{\langle}
\def\ot{\otimes}
\def\ov{\overline}
\def\part{\partial}
\def\ra{\rightarrow}
\def\rgl{\rangle}
\def\sbs{\subset}

\def\ts{\times}
\def\wh{\widehat}
\def\wt{\widetilde}
\def\wdg{\wedge}

\def\Covol{\mathop{\rm Covol}\nolimits}
\def\Im{\mathop{\rm Im}\nolimits}
\def\Ker{\mathop{\rm Ker}\nolimits}
\def\mod{\mathop{\rm mod}\nolimits}
\def\Supp{\mathop{\rm Supp}\nolimits}
\def\Re{\mathop{\rm Re}\nolimits}

\def\hra{\hookrightarrow}
\def\longra{\longrightarrow}
\def\Ra{\Rightarrow}
\def\Lra{\Leftrightarrow}

\def\xx{\vrule height 0.3em depth 0.2em width 0.3em}

\def\build#1_#2^#3{\mathrel{
\mathop{\kern 0pt#1}\limits_{#2}^{#3}}}

\def\a{\alpha}
\def\b{\beta}
\def\d{\delta}
\def\D{\Delta}
\def\g{\gamma}
\def\G{\Gamma}
\def\lb{\lambda}
\def\L{\Lambda}
\def\Om{\Omega}
\def\s{\sigma}
\def\Si{\Sigma}
\def\t{\theta}
\def\ve{\varepsilon}
\def\vp{\varphi}
\def\z{\zeta}

\font\tenbb=msbm10
\font\sevenbb=msbm7
\font\fivebb=msbm5
\newfam\bbfam
\textfont\bbfam=\tenbb \scriptfont\bbfam=\sevenbb
\scriptscriptfont\bbfam=\fivebb
\def\bb{\fam\bbfam}

\def\Rb{{\bb R}}
\def\Qb{{\bb Q}}
\def\Cb{{\bb C}}
\def\Fb{{\bb F}}
\def\Zb{{\bb Z}}
\def\Tb{{\bb T}}

\def\Hc{{\cal H}}

\def\Sc{{\cal S}}
\def\Xc{{\cal X}}
\def\Uc{{\cal U}}

\def\build#1_#2^#3{\mathrel{
\mathop{\kern 0pt#1}\limits_{#2}^{#3}}}

\catcode`\@=11
\def\displaylinesno #1{\displ@y\halign{
\hbox to\displaywidth{$\@lign\hfil\displaystyle##\hfil$}&
\llap{$##$}\crcr#1\crcr}}

\def\ldisplaylinesno #1{\displ@y\halign{
\hbox to\displaywidth{$\@lign\hfil\displaystyle##\hfil$}&
\kern-\displaywidth\rlap{$##$}
\tabskip\displaywidth\crcr#1\crcr}}
\catcode`\@=12

\def\lgl{\langle}
\def\rgl{\rangle}
\def\ov{\overline}
\def\fl{\forall}
\def\ify{\infty}
\def\ot{\otimes}
\def\part{\partial}
\def\sbs{\subset}
\def\sm{\simeq}
\def\ts{\times}

\def\ra{\rightarrow}

\def\wt{\widetilde}
\def\wh{\widehat}

\def\semi{\mathop{>\!\!\!\triangleleft}}


\vglue 5mm

\centerline{\bf Trace formula in noncommutative Geometry and}

\smallskip

\centerline{\bf the zeros of the Riemann zeta function}

\bigskip

\centerline{Alain CONNES}

\vglue 1cm
\vglue 1cm

\noindent {\bf Abstract } 
{\it We give a spectral interpretation of the critical zeros of the 
Riemann zeta function as an absorption spectrum, while eventual 
noncritical zeros appear as resonances. We give a geometric 
interpretation of the explicit formulas of number theory as a trace 
formula on the noncommutative space of Adele classes. This reduces the 
Riemann hypothesis to
 the validity of the trace formula and eliminates the parameter $\d$ of 
our previous approach.}
\bigskip
\vglue 1cm

\centerline{\bf Table of contents}

\vglue 1cm
 
\noindent {\bf Introduction. }

\noindent {\bf I Quantum chaos and the hypothetical Riemann flow.}

\noindent {\bf II Algebraic Geometry and global fields of non zero 
characteristic.}

\noindent {\bf III Spectral interpretation of critical zeros.}

\noindent {\bf IV The distribution trace formula for flows on
manifolds.}

\noindent {\bf V The action $(\lb ,x) \ra \lb \, x$ of $K^*$ on a local 
field 
$K$.}

\noindent {\bf VI The global case, and the formal trace computation.}

\noindent {\bf VII Proof of the trace formula in the $S$-local case.}

\noindent {\bf VIII The trace formula in the global case, and 
elimination of $\d$.}

\noindent {\bf Appendix I, Proof of theorem 1.}

\noindent {\bf Appendix II, Explicit formulas.}

\noindent {\bf Appendix III, Distribution trace formulas. }
\vglue 1cm
\eject
\noindent {\bf Introduction }
\smallskip
We shall give in this paper a spectral interpretation of the zeros of 
the Riemann zeta function and a geometric framework to which one can 
transpose the ideas of algebraic geometry involving the action of the 
Frobenius and the Lefchetz formula. The spectral interpretation of the 
zeros of zeta will be as an absorption spectrum, i.e. as missing 
spectral lines. All zeros will play a role in the spectral side of the 
trace formula, but while the critical zeros will appear per-se, the non 
critical ones will appear as resonances and 
enter in the trace formula through their harmonic potential with respect 
to the critical line. Thus the spectral side is entirely canonical, and 
by proving positivity of the Weil distribution, we shall show that its 
equality with the geometric side, i.e. the global trace formula, is 
equivalent to the Riemann Hypothesis for all $L$-functions with 
Gr\"ossencharakter. We shall model our discussion on the Selberg trace 
formula, but it differs from the latter in several important respects. 
We shall first explain in particular why a crucial negative sign in the 
analysis of the statistical fluctuations of the zeros of zeta indicates 
that the spectral interpretation should be as an absorption spectrum, or 
equivalently should be of a cohomological nature. As it turns out, the 
geometric framework involves an innocent looking space, the space $X$ of 
Adele classes, where two adeles which belong to the same orbit of the 
action of $GL_1(k)$ ($k$ a global field), are considered equivalent. The 
group $C_k= \,GL_1(A)/GL_1(k)$ of Idele classes (which is the class 
field theory counterpart of the Galois group) acts by multiplication on 
$X$.
\smallskip
 Our first preliminary result (theorem 1 of section III) gives a 
spectral interpretation of the critical 
zeros of zeta and $L$ functions on a global field $k$ from the action of 
the Idele class group on a space of square integrable functions on
the space $X = A/k^*$ of Adele classes. Corollary 2 gives the 
corresponding computation of the spectral trace.
This result is only preliminary because it requires the use of an 
unnatural parameter $\d$ which plays the role of a Sobolev exponent and 
allows to see the absorption spectrum as a point spectrum. 
\smallskip
Our second preliminary result is a formal computation (section VI) of 
the character of the representation of the Idele class group on the 
above $L^2$ space. This formal computation gives the Weil distribution 
which is the essential ingredient of the Riemann-Weil explicit formula. 
At this point (which was the situation in [Co]), the main problems are 
to give a rigorous meaning to the formal trace computation and to 
eliminate the unwanted parameter $\d$.
\smallskip
These two problems will be solved in the present paper. We first prove a 
trace formula (theorem 3 of section V) for the action of the 
multiplicative group $K^*$ of a local field $K$
on the Hilbert space $L^2 (K)$, and (theorem 4 of section VII) a trace 
formula for the action of the multiplicative group $C_S$ of Idele 
classes associated to a finite set $S$ of places of a global field $k$, 
on the Hilbert space of square integrable functions $L^2 (X_S)$, where 
$X_S$ is the quotient of $\prod_{v \in S} k_v$ by the action of the 
group $O^*_S$ of S-units of $k$. In both cases we obtain exactly the 
terms of the Weil explicit formulas which belong to the finite set of 
places. This result is quite important since the space $X_S$ is highly 
non trivial as soon as the cardinality of $S$ is larger or equal to $3$. 
Indeed this quotient space is non type I in the sense on Noncommutative 
Geometry and it is reassuring that the trace formula continues to hold 
there.
\smallskip
We check in detail (theorem 6 of Appendix II) that the rewriting of the 
Weil explicit formulas which is predicted by the global trace formula is 
correct.
\smallskip
 Finally, we eliminate in section VIII, using ideas which are common 
both to the Selberg trace formula and to the standard explanation of the 
absorption lines in physics, the unpleasant parameter $\d$ which 
appeared as a label of the function spaces of section III. We write the 
global trace formula as an analogue of the Selberg trace formula.
The validity of the trace formula for any finite set of places follows 
from theorem 4 of section VII, but in the global case is left open and 
shown ( Theorem 5 of section VIII) to be equivalent to the validity of 
the Riemann Hypothesis for all $L$ functions with Gr\"ossencharakter. 
This equivalence, together with the plausibility of a direct proof of 
the trace formula along the lines of theorem 4 (section VII) constitute 
the main result of this paper.
The elimination of the parameter $\d$ is the main improvement of the 
present paper with respect to [Co].\smallskip
It is an old idea, due to Polya and Hilbert that in order to
understand the location of the zeros of the Riemann zeta
function, one should find a Hilbert space $\Hc$ and an
operator $D$ in $\Hc$ whose spectrum is given by the non trivial
zeros of the zeta function. The hope then is that suitable
selfadjointness properties of $D$ (of $i \left( D-{1\over 2}
\right)$ more precisely) or positivity properties of $\D =
D(1-D)$ will be easier to handle than the original
conjecture. The main reasons why this idea should be taken
seriously are first the work of A. Selberg ([Se]) in which a
suitable Laplacian $\D$ is related in the above way to an
analogue of the zeta function, and secondly the theoretical
([M][B][KS]) and experimental evidence ([O][BG]) on the fluctuations of
the spacing between consecutive zeros of zeta. The number of
zeros of zeta whose imaginary part is less than $E > 0$,
$$
N(E) = \# \ \hbox{of zeros} \ \rho \ , \ 0 < {\rm Im} \,
\rho < E \leqno (1)
$$
has an asymptotic expression ([R]) given by
$$
N(E) = {E \over 2\pi} \ \left( \log \left( {E \over 2\pi}
\right) -1 \right) + {7 \over 8} + o(1) + N_{\rm osc} (E)
\leqno (2)
$$
where the oscillatory part of this step function is
$$
N_{\rm osc} (E) = {1\over \pi} \ {\rm Im} \, \log \, \z \,
\left( {1\over 2} + iE \right) \leqno (3)
$$
assuming that $E$ is not the imaginary part of a zero and
taking for the logarithm the branch which is $0$ at $+\ify$.

\medskip

One shows (cf. [Pat]) that $N_{\rm osc} (E)$ is $O(\log E)$. In
the decomposition (2) the two terms $\lgl N(E) \rgl = N(E) -
N_{\rm osc} (E)$ and $N_{\rm osc} (E)$ play an independent
role. The first one $\lgl N(E) \rgl$ which gives the average
density of zeros just comes from Stirling's formula and is
perfectly controlled. The second $N_{\rm osc} (E)$ is a
manifestation of the randomness of the actual location of the
zeros, and to eliminate the role of the density one returns to
the situation of uniform density by the transformation
$$
x_j = \lgl N(E_j) \rgl \quad (E_j \ \hbox{the $j^{\rm th}$
imaginary part of zero of zeta}) \, .\leqno (4)
$$
Thus the spacing between two consecutive $x_j$ is now 1 in
average and the only information that remains is in the
statistical fluctuation. As it turns out ([M][O]) these fluctuations
are the same as the fluctuations of the eigenvalues of a
random hermitian matrix of very large size.

\smallskip

\noindent H. Montgomery [M] proved (assuming RH) a weakening of the 
following conjecture (with $\a,\b >0$),
$$
\ldisplaylinesno{
{\rm Card} \, \{ (i,j) \, ; \, i,j \in 1,\ldots ,M \, ; \,
x_i - x_j \in [\a ,\b ]\} \cr
\sim M \int_{\a}^{\b} \left( 1-
\left( {\sin (\pi u) \over \pi u} \right)^2 \right) \, du
&(5) \cr
}
$$
 This law (5) is precisely the same as the
correlation between eigenvalues of hermitian matrices of the
gaussian unitary ensemble ([M]). Moreover,
numerical tests due to A. Odlyzko ([O][BG]) have confirmed with great 
precision the
behaviour (5) as well as the analoguous behaviour for more
than two zeros. In [KS], N. Katz and P. Sarnak proved an analogue of the 
Montgomery-Odlyzko law for zeta and L-functions of function fields over 
curves.

\medskip

It is thus an excellent motivation to try and find a natural
pair $(\Hc ,D)$ where naturality should mean for instance
that one should not even have to define the zeta function, let alone 
its analytic continuation, in
order to obtain the pair (in order for instance to avoid the
joke of defining $\Hc$ as the $\ell^2$ space built on the
zeros of zeta).

\vglue 1cm

\noindent {\bf I Quantum chaos and the hypothetical Riemann flow}

\medskip
\noindent Let us first describe following [B] the direct atempt to 
construct the Polya-Hilbert space from quantization of a classical 
dynamical system. The original motivation for the theory of random
matrices comes from quantum mechanics. In this theory the
quantization of the classical dynamical system given by the
phase space $X$ and hamiltonian $h$ gives rise to a Hilbert
space $\Hc$ and a selfadjoint operator $H$ whose spectrum is
the essential physical observable of the system. For
complicated systems the only useful information about this
spectrum is that, while the average part of the counting
function, $$
N(E) = \# \ \hbox{eigenvalues of $H$ in} \ [0,E] \leqno (1)
$$
is computed by a semiclassical approximation mainly as a
volume in phase space, the oscillatory part,
$$
N_{\rm osc} (E) = N(E) - \lgl N(E) \rgl \leqno (2)
$$
is the same as for a random matrix, governed by the statistic
dictated by the symmetries of the system.

\smallskip

\noindent In the absence of a magnetic field, i.e. for a
classical hamiltonian of the form,
$$
h = {1\over 2m} \, p^2 + V(q) \leqno (3)
$$
where $V$ is a real-valued potential on configuration space,
there is a natural symmetry of classical phase space, called time 
reversal symmetry,
$$
T(p,q) = (-p,q) \leqno (4)
$$
which preserves $h$, and entails that the correct ensemble on
the random matrices is not the above GUE but rather the
gaussian orthogonal ensemble: GOE. Thus the oscillatory part
$N_{\rm osc} (E)$ behaves in the same way as for a random {\it
real symmetric} matrix.

\smallskip

\noindent Of course $H$ is just a specific operator in $\Hc$
and, in order that it behaves {\it generically} it is
necessary (cf. [B]) that the classical hamiltonian system
$(X,h)$ be {\it chaotic} with isolated {\it periodic orbits}
whose instability exponents (i.e. the logarithm of the
eigenvalues of the Poincar\'e return map acting on the
transverse space to the orbits) are different from 0.

\smallskip

\noindent One can then ([B]) write down an asymptotic
semiclassical approximation to the oscillatory function
$N_{\rm osc} (E)$
$$
N_{\rm osc} (E) = {1\over \pi} \ {\rm Im} \int_0^{\ify} {\rm
Trace} (H-(E+i\eta))^{-1} \, id\eta \leqno (5)
$$
using the stationary phase approximation of the
corresponding functional integral. For a system whose
configuration space is 2-dimensional, this gives ([B] (15)),
$$
N_{\rm osc} (E) \sm {1\over \pi} \sum_{\g_p} \sum_{m=1}^{\ify}
{1\over m} \, {1\over 2{\rm sh} \left( {m\lb_p \over
2}\right)} \, \sin (S_{\rm pm} (E)) \leqno (6)
$$
where the $\g_p$ are the primitive periodic orbits, the label
$m$ corresponds to the number of traversals of this orbit,
while the corresponding instability exponents are $\pm
\lb_p$. The phase $S_{\rm pm} (E)$ is up to a constant equal
to $m \, E \, T_{\g}^{\#}$ where $T_{\g}^{\#}$ is the period
of the primitive orbit $\g_p$.

\smallskip

\noindent The formula (6) gives very precious information ([B]) on
the hypothetical ``Riemann flow'' whose quantization should
produce the Polya-Hilbert space. The point is that the Euler
product formula for the zeta function yields (cf. [B]) a
similar asymptotic formula for $N_{\rm osc} (E)$ (3),
$$
N_{\rm osc} (E) \sm {-1 \over \pi} \sum_p \sum_{m=1}^{\ify}
{1\over m} \, {1 \over p^{m/2}} \, \sin \, (m \, E \, \log \,
p) \, . \leqno (7)
$$
Comparing (6) and (7) gives the following information,

\medskip

\item{(A)} The periodic primitive orbits should be labelled
by the prime numbers $p=2,3,5,7,\ldots$, their periods should
be the $\log p$ and their instability exponents $\lb_p = \pm
\log p$.

\medskip

\noindent Moreover, since each orbit is only counted once, the
Riemann flow should not possess the symmetry $T$ of (4) whose
effect would be to duplicate the count of orbits. This last
point excludes in particular the geodesic flows since they have
the time reversal symmetry $T$. Thus we get 
\medskip

\item{(B)} The Riemann flow cannot satisfy time reversal symmetry.

\medskip
\noindent However there are two important mismatches  (cf.
[B]) between the two formulas (6) and (7). The first one is
the overall {\it minus sign} in front of formula (7), the
second one is that though $2 {\rm sh} \, \left( {m \lb_p
\over 2} \right) \sim p^{m/2}$ when $m \ra \ify$, we do not
have an equality for finite values of $m$.

\smallskip

\noindent These are two fundamental difficulties and in order
to overcome them we shall use the well known strategy of
extending the problem to the case of arbitrary {\it global fields}. By
specialising to the function field case we shall then obtain
additional precious information.
 
\vglue 1cm

\noindent {\bf II Algebraic Geometry and global fields of non zero 
characteristic}

\medskip

The basic properties of the Riemann zeta function extend to zeta
 functions associated to an arbitrary global field, and it is unliquely 
that one can settle the problem of the spectral interpretation of the 
zeros,
let alone find the Riemann flow, for the particular case of the global 
field
 $\Qb$ of rational numbers without at the same time settling these 
problems 
for all global fields.
The conceptual definition of such fields $k$, is the following:
\smallskip

\noindent {\it A field $k$ is a global field iff it is } discrete {\it 
and} cocompact {\it in a (non
discrete) locally compact semisimple abelian ring $A$.}
\smallskip
 As it turns out $A$ then depends functorially on $k$ and is called
the Adele ring of $k$, often denoted by $k_A$. Thus though the field $k$ 
itself has no 
interesting topology, there is a canonical and highly non trivial 
topological ring which 
is canonically associated to $k$.
 When the
characteristic $p$ of a global field $k$ is $>0$, the field
$k$ is the function field of a non singular algebraic curve
$\Si$ defined over a finite field $\Fb_q$ included in $k$ as
its maximal finite subfield, called the field of constants.
One can then apply the ideas of algebraic geometry, first
developed over $\Cb$, to the geometry of the curve $\Si$ and
obtain a geometric interpretation of the basic properties of
the zeta function of $k$; the dictionary contains in
particular the following lines
\medskip

$$
\ldisplaylinesno{
\matrix{
\hbox{Spectral interpretation of} \quad &\quad \hbox{Eigenvalues of 
action} \cr
\hbox{the zeros} \quad &\quad \hbox{of Frobenius on $\ell$-adic 
cohomology}
\cr \cr
\hbox{Functional equation} \quad &\quad \hbox{Riemann Roch
theorem} \cr
&\quad \hbox{(Poincar\'e duality)} \cr \cr
\hbox{Explicit formulas of} \quad &\quad \hbox{Lefchetz
formula} \cr
\hbox{number theory} \quad &\quad \hbox{for the Frobenius}
\cr \cr
\hbox{Riemann hypothesis} \quad &\quad \hbox{Castelnuovo
positivity} \cr } &(1)
}
$$
\medskip
Since $\Fb_q$ is not algebraically closed, the points of
$\Si$ defined over $\Fb_q$ do not suffice and one needs to
consider $\bar{\Si}$, the points of $\Si$ on the algebraic
closure $\bar{\Fb}_q$ of $\Fb_q$, which is obtained by
adjoining to $\Fb_q$ the roots of unity of order prime to
$q$. This set of points is a countable union of periodic
orbits under the action of the Frobenius automorphism, these
orbits are parametrized by the set of places of $k$ and their
periods are indeed given by the analogues of the $\log p$ of
(A). Being a countable set it does not qualify for analogue
of the Riemann flow and it only aquires an interesting
structure from algebraic geometry. The minus sign which was
problematic in the above discussion admits here a beautiful
resolution since the analogue of the Polya-Hilbert space is
given, if one replaces $\Cb$ by $\Qb_{\ell}$ the field of
$\ell$-adic numbers $\ell \not= p$, by the cohomology group
$$
H_{\rm et}^1 (\bar{\Si} ,\Qb_{\ell}) \leqno (2)
$$
which appears with an overall minus sign in the Lefchetz
formula 
$$
\sum (-1)^j \, \hbox{Trace} \, \varphi^* / H^j = \sum_{\varphi (x) = x} 
\, 1. \leqno (3)
$$
\smallskip

\noindent For the general case this suggests

\medskip\item{(C)} The Polya-Hilbert space $\Hc$ should appear from
its negative  $\ominus \Hc$.

\medskip
\noindent In other words, the spectral interpretation of the zeros 
of the Riemann zeta function should be as an absorption spectrum rather 
than as an emission spectrum, to borrow the language of spectroscopy.
\smallskip
\noindent The next thing that one learns from this excursion
in characteristic $p > 0$ is that in that case one is not
dealing with a flow but rather with a single transformation.
In fact taking advantage of abelian covers of $\Si$ and of
the fundamental isomorphism of class field theory one finds
that the natural group that should replace $\Rb$ for the
general Riemann flow is the Idele class group:
$$
C_k = {\rm GL}_1 (A) / k^* \, . \leqno (D)
$$
We can thus collect the information (A) (B) (C) (D) that we
have obtained so far and look for the Riemann flow as an
action of $C_k$ on an hypothetical space $X$.

\vglue 1cm

\noindent {\bf III Spectral interpretation of critical zeros}

\medskip

\noindent There is a third approach to the problem of the
zeros of the Riemann zeta function, due to G. P\'olya [P] and M. Kac [K] 
and pursued
further in [J] [BC]. It is based on statistical mechanics and
the construction of a quantum statistical system whose {\it
partition function} is the Riemann zeta function. Such a
system was naturally constructed in [BC] and it does indicate
using the first line of the dictionary of Noncommutative Geometry 
(namely the correspondence between quotient spaces and noncommutative 
algebras) what the space $X$
should be in general:
$$
X = A/k^* \leqno (1)
$$
namely the quotient of the space $A$ of adeles, $A = k_A$ by
the action of the multiplicative group $k^*$,
$$
a \in A \ , \ q \in k^* \ra aq \in A \, . \leqno (2)
$$
This space $X$ already appears in a very implicit manner in the work of 
Tate and Iwasawa
on the functional equation. It is a noncommutative space in that, even 
at the level of measure theory,
it is a tricky quotient space. For instance
 at the measure theory level, the corresponding von
Neumann algebra,
$$
R_{01} = L^{\ify} (A) \semi k^* \leqno (3)
$$
where $A$ is endowed with its Haar measure as an additive
group, is the hyperfinite factor of type ${\rm II}_{\ify}$.

\smallskip

\noindent The idele class group $C_k$ acts on $X$ by
$$
(j,a) \ra ja \qquad \fl \, j \in C_k \ , \ a \in X \leqno (4)
$$
and it was exactly necessary to divide $A$ by $k^*$ so that
(4) makes good sense.

\smallskip

\noindent We shall come back later to the analogy between the
action of $C_k$ on $R_{01}$ and the action of the Galois
group of the maximal abelian extension of $k$.

\smallskip

\noindent What we shall do now is to construct the Hilbert
space $L_{\d}^2$ of functions on $X$ with growth
indexed by $\d > 1$. Since $X$ is a quotient space we shall
first learn in the usual manifold case how to obtain the
Hilbert space $L^2 (M)$ of square integrable functions on a manifold $M$
by working only on the universal cover $\wt M$ with the
action of $\G = \pi_1 (M)$. Every function $f \in C_c^{\ify}
(\wt M )$ gives rise to a function $\wt f$ on $M$ by
$$
\wt f (x) = \sum_{\pi (\wt x ) =x} f(\wt x ) \leqno (5)
$$
and all $g \in C^{\ify} (M)$ appear in this way. Moreover,
one can write the Hilbert space inner product $\int_M
\wt{f}_1 (x) \, \wt{f}_2 (x) \, dx$, in terms of $f_1$ and
$f_2$ alone. Thus $\Vert \wt f \Vert^2 = \int \left\vert
\sum_{\g \in \G} \, f(\g x) \right\vert^2 \, dx$ where the
integral is performed on a fundamental domain for $\G$ acting
on $\wt M$. This formula defines a prehilbert space norm on
$C_c^{\ify} (\wt M )$ and $L^2 (M)$ is just the completion of
$C_c^{\ify} (\wt M )$ for that norm. Note that any function
of the form $f-f_{\g}$ has vanishing norm and hence
disappears in the process of completion. In our case of $X =
A/k^*$ we thus need to define the analoguous norm on the Bruhat-Schwartz 
space $\Sc (A)$ of functions on $A$ (cf Appendix I for the general 
definition 
of the Bruhat-Schwartz space). Since 0 is
fixed by the action of $k^*$ the expression $\sum_{\g \in
k^*} \, f (\g x)$ does not make sense for $x=0$ unless we
require that $f(0) =0$. Moreover, when $\vert x \vert \ra 0$,
the above sums approximate, as Riemann sums, the product of
$\vert x \vert^{-1}$ by $\int f \, dx$ for the additive Haar
measure, thus we also require $\int f \,
dx = 0$. We can now define the Hilbert space $L_{\d}^2 (X)_0$
as the completion of the codimension 2 subspace
$$
\Sc (A)_0 = \{ f \in \Sc (A) \ ; \ f(0) = 0 \ , \ \int f \,
dx =0 \} \leqno (6)
$$
for the norm $\Vert \ \Vert_{\d}$ given by
$$
\Vert f \Vert_{\d}^2 = \int \Bigl\vert \sum_{q\in k^*} f (qx)
\Bigl\vert^2 \,  (1+\log^2 \vert x \vert)^{\d /2} \,\vert x \vert \,
d^* x \leqno (7)
$$
where the integral is performed on $A^* /k^*$ and $d^* x$ is
the multiplicative Haar measure on $A^* / k^*$. The ugly term
$(1+\log^2 \vert x \vert)^{\d /2}$ is there to control the
growth of the functions on the non compact quotient. We shall see how to 
remove it later in section VII. Note that $\vert qx\vert = \vert x\vert$
 for any $q\in k^*$. 
\smallskip

\noindent   The key point is that we use the
measure $\vert x \vert \, d^* x$ instead of the additive Haar
measure $dx$. Of course for a local field $K$ one has $dx =
\vert x \vert \, d^* x$ but this fails in the above global
situation. Instead one has,
$$
dx = \, \build \lim_{\ve \ra 0}^{} \, \ve \, \vert x
\vert^{1+\ve} \, d^* x \, . \leqno (8)
$$
One has a natural representation of $C_k$ on $L_{\d}^2 (X)_0$
given by
$$
(U(j) \, f) \, (x) = f(j^{-1} \, x) \qquad \fl \, x \in A \ ,
\ j \in C_k \leqno (9)
$$
and the result is independent of the choice of a lift of $j$
in $J_k = {\rm GL}_1 (A)$ because the functions $f-f_q$ are
in the kernel of the norm. The conditions (6) which define
$\Sc (A)_0$ are invariant under the action of $C_k$ and give
the following action of $C_k$ on the 2-dimensional supplement
of $\Sc (A)_0 \subset \Sc (A)$; this supplement is $\Cb \oplus
\Cb (1)$ where $\Cb$ is the trivial $C_k$ module
(corresponding to $f(0)$) while the Tate twist $\Cb (1)$ is
the module
$$
(j,\lb) \ra \vert j \vert \, \lb \leqno (10)
$$
coming from the equality
$$
\int f(j^{-1} \, x) \, d \, x = \vert j \vert \int f(x) \, dx
\, . \leqno (11)
$$
In order to analyse the representation (9) of $C_k$ on
$L_{\d}^2 (X)_0$ we shall relate it to the left regular
representation of the group $C_k$ on the Hilbert space
$L_{\d}^2 (C_k)$ obtained from the following Hilbert space
square norm on functions,
$$
\Vert \xi \Vert_{\d}^2 = \int_{C_k} \vert \xi (g)\vert^2 \,
(1+\log^2 \vert g \vert)^{\d /2} \, d^* g \leqno (12)
$$
where we have normalized the Haar measure of the
multiplicative group $C_k$, with module,
$$
\vert \ \vert : C_k \ra \Rb_+^* \leqno (13)
$$
in such a way that (cf. [W3])
$$
\int_{\vert g \vert \in [1,\L]} \, d^* g \sim \log \L \quad
\hbox{when} \quad \L \ra +\ify \, . \leqno (14)
$$
The left regular representation $V$ of $C_k$ on $L_{\d}^2
(C_k)$ is
$$
(V(a) \, \xi) \, (g) = \xi (a^{-1} \, g) \qquad \fl \, g,a
\in C_k \, . \leqno (15)
$$
Note that because of the weight $(1+\log^2 \vert x \vert)^{\d
/2}$, this representation is {\it not} unitary but it
satisfies the growth estimate 
$$
\Vert V(g)\Vert = 0 \, (\log \vert g \vert)^{\d /2} \quad
\hbox{when} \quad \vert g \vert \ra \ify \leqno (16)
$$
which follows from the inequality (valid for $u,v \in \Rb$)
$$
\rho (u+v) \leq 2^{\d /2} \, \rho (u) \, \rho (v) \ , \ \rho
(u) = (1+u^2)^{\d /2} \, . \leqno (17)
$$
We let $E$ be the linear isometry from $L_{\d}^2 (X)_0$ into
$L_{\d}^2 (C_k)$ given by the equality,
$$
E(f) \, (g) = \vert g \vert^{1/2} \sum_{q\in k^*} f(qg)
\qquad \fl \, g \in C_k \, . \leqno (18)
$$
By comparing (7) with (12) we see that $E$ is an isometry
and the factor $\vert g \vert^{1/2}$ is dictated by comparing
the measures $\vert g \vert \, d^* g$ of (7) with $d^* g$ of
(12).

\smallskip

\noindent One has $E(U(a) \, f) \, (g) = \vert g \vert^{1/2}
\, \sum_{k^*} \, (U(a) \, f) \, (qg) = \vert g \vert^{1/2} \,
\sum_{k^*} \, f(a^{-1} \, qg) = \vert a \vert^{1/2} \,
\vert a^{-1} \, g \vert^{1/2} \, \sum_{k^*} \, f(q \, a^{-1}
\, g) = \vert a \vert^{1/2} \, (V(a) \, E(f))$ $(g)$.

\smallskip

\noindent Thus,
$$
E \, U (a) = \vert a \vert^{1/2} \, V(a) \, E \, . \leqno (19)
$$
This equivariance shows that the range of $E$ in $L_{\d}^2
(C_k)$ is a closed invariant subspace for the representation
$V$.

\smallskip

\noindent The following theorem and its corollary show that
the cokernel $\Hc = L_{\d}^2 (C_k) /$ ${\rm Im} (E)$ of the
isometry $E$ plays the role of the Polya-Hilbert space. Since
${\rm Im} \, E$ is invariant under the representation $V$ we
let $W$ be the corresponding representation of $C_k$ on $\Hc$.

\smallskip

\noindent The abelian locally compact group $C_k$ is (non
canonically) isomorphic to $K \ts N$ where
$$
K = \{ g \in C_k \ ; \ \vert g \vert =1 \} \ , \ N = {\rm
range} \ \vert \ \vert \sbs \Rb_+^* \, . \leqno (20)
$$
For number fields one has $N=\Rb_+^*$ while for fields of non
zero characteristic $N \sm \Zb$ is the subgroup $q^{\Zb} \sbs
\Rb_+^*$ .(Where $q=p^{\ell}$ is the cardinality of the field of 
constants).

\smallskip

\noindent We choose (non canonically) an isomorphism
$$
C_k \sm K \ts N \, . \leqno (21)
$$
By construction the representation $W$ satisfies (using (16)),
$$
\Vert W (g) \Vert = 0 (\log \vert g \vert)^{\d /2} \leqno (22)
$$
and its restriction to $K$ is unitary. Thus $\Hc$ splits as a
canonical direct sum of pairwise orthogonal subspaces,
$$
\Hc = \ \build \oplus_{\chi \in \wh K}^{} \, \Hc_{\chi} \ , \
\Hc_{\chi} = \{ \xi \ ; \ W(g) \, \xi = \chi \, (g) \, \xi \ , \
\fl \, g \in K \} \leqno (23)
$$
where $\chi$ runs through the Pontrjagin dual group of $K$,
which is the discrete abelian group $\wh K$ of characters of
$K$. Using the non canonical isomorphism (21), i.e. the
corresponding inclusion $N \sbs C_k$ one can now restrict the
representation $W$ to any of the sectors $\Hc_{\chi}$. When
${\rm char} (k) > 0$, then $N \sm \Zb$ and the condition (22)
shows that the action of $N$ on $\Hc_{\chi}$ is given by a
single operator with {\it unitary} spectrum. (One uses the
spectral radius formula $\vert {\rm Spec} \, w \vert =
\overline{\rm Lim} \, \Vert \,w^n \Vert^{1/n}$.) When ${\rm
Char} (k) = 0$, we are dealing with an action of $\Rb_+^* \sm
\Rb$ on $\Hc_{\chi}$ and the condition (22) shows that this
representation is generated by a closed unbounded operator $D_{\chi} $
with purely imaginary spectrum. The resolvent $R_{\lb} =
(D_{\chi} -\lb)^{-1}$ is given, for ${\rm Re} \, \lb > 0$, by the
equality
$$
R_{\lb} = \int_0^{\ify} W_{\chi} (e^s) \, e^{-\lb s} \, ds
\leqno (24)
$$
and for ${\rm Re} \, \lb < 0$ by,
$$
R_{\lb} = \int_0^{\ify} W_{\chi} (e^{-s}) \, e^{\lb s} \, ds
\leqno (25)
$$
while the operator $D_{\chi} $ is defined by
$$
D_{\chi}  \, \xi = \build {\rm lim}_{\ve \ra 0}^{} \, {1\over \ve} \,
(W_{\chi} (e^{\ve}) -1) \, \xi \, . \leqno (26)
$$

\bigskip

\noindent {\bf Theorem 1.} {\it Let $\chi \in \wh K$, $\d > 1$,
$\Hc_{\chi}$ and $D_{\chi} $ be as above. Then $D_{\chi} $ has discrete
spectrum, ${\rm Sp} \, D_{\chi}  \sbs i \, \Rb$ is the set of imaginary 
parts of zeros of
the $L$ function with Gr\"ossencharakter $\wt{\chi}$ which have
real part equal to ${1\over 2}$; $\rho \in {\rm Sp} \, D
\Leftrightarrow L \left(\wt{\chi} , {1 \over 2} +\rho \right)
=0$ and $\rho \in i \, \Rb$, where $\wt{\chi}$ is the unique
extension of $\chi$ to $C_k$ which is equal to $1$ on $N$.
Moreover the multiplicity of $\rho$ in ${\rm Sp} \, D$ is equal
to the largest integer  $ n < {1+\d \over 2}$ , $ n\leq $ multiplicity 
of
${1\over 2} + \rho$ as a zero of $L$.}

\bigskip
 Theorem 1 has a similar formulation when the characteristic of $k$ is 
non zero. The following corollary is valid for global fields $k$ of 
arbitrary characteristic.
\bigskip
\noindent {\bf Corollary 2.} {\it For any Schwartz function $h
\in \Sc (C_k)$ the operator $W(h) = \int W(g) \, h(g) \, d^*
\, g$ in $\Hc$ is of trace class, and its trace is given by
$$
\hbox{Trace} \, W(h) = \sum_{{L\left(\wt{\chi} , {1\over 2} + \rho
\right) =0 \atop \rho \in { i \, \Rb / N^{ \bot}}}} \wh h (\wt{\chi} 
,\rho)
$$
where the multiplicity is counted as in Theorem 1 and where
the Fourier transform $\wh h$ of $h$ is defined by,
$$
\wh h (\wt{\chi} ,\rho) = \int_{C_k} h(u) \, \wt{\chi} (u) \, \vert u
\vert^\rho \, d^* \, u \, . $$ }

\bigskip

Note that we did not have to define the $L$ functions,
 let alone their analytic continuation,
 before
stating the theorem, which shows that the pair
$$
(\Hc_{\chi} ,D_{\chi}) \leqno (27)
$$
certainly qualifies as a Polya-Hilbert space.

\smallskip

\noindent The case of the Riemann zeta function corresponds to
the trivial character $\chi =1$ for the global field $k=\Qb$
of rational numbers.

\smallskip

\noindent In general the zeros of the $L$ functions can have
multiplicity but one expects that for a fixed Gr\"ossencharakter
 $\chi$ this multiplicity is bounded, so that for a
large enough value of $\d$ the spectral multiplicity of $D$
will be the right one. When the characteristic of $k$ is $>0$
this is certainly true.

\smallskip

\noindent If we modify the choice of non canonical isomorphism
(21) this modifies the operator $D$ by
$$
D' = D - i \, s \leqno (28)
$$
where $s \in \Rb$ is determined by the equality
$$
\wt{\chi}' (g) = \wt{\chi} (g) \, \vert g \vert^{\rm i \, s}
\qquad \fl \, g \in C_k \, . \leqno (29)
$$
The coherence of the statement of the theorem is insured by the equality
$$
L(\wt{\chi}' ,z) = L(\wt{\chi} , z + i \, s) \qquad \fl \, z
\in \Cb \, . \leqno (30)
$$
When the zeros of $L$ have multiplicity and $\d$ is large
enough the operator $D$ is {\it not} semisimple and has a non
trivial Jordan form (cf. Appendix I). This is compatible with the almost
unitary condition (22) but not with skew symmetry for $D$.
\smallskip

\noindent The proof of theorem 1, explained in Appendix I, is based on 
the distribution theoretic 
interpretation by A. Weil [W2] of the idea of Tate and Iwasawa on the 
functional 
equation. Our construction should be compared with [Bg] and [Z]. 
\smallskip

\noindent As we expected from (C), the Polya-Hilbert space
$\Hc$ appears as a cokernel. Since we obtain the Hilbert space
$L_{\d}^2 (X)_0$ by imposing two linear conditions on $\Sc (A)$,
$$
0 \ra \Sc (A)_0 \ra \Sc (A) \build \ra_{}^{L} \Cb \oplus \Cb
(1) \ra 0 \leqno (31)
$$
we shall define $L_{\d}^2 (X)$ so that it fits in an exact sequence of  
$C_k$-modules
 
$$
0 \ra L_{\d}^2 (X)_0 \ra L_{\d}^2 (X) \ra \Cb \oplus \Cb (1) \ra
0 \, . \leqno (32)
$$
We can then use the exact sequence of $C_k$-modules
$$
0 \ra L_{\d}^2 (X)_0 \ra L_{\d}^2 (C_k) \ra \Hc \ra 0 \leqno
(33)
$$
together with Corollary 2 to compute in a formal manner what
the character of the module $L_{\d}^2 (X)$ should be. Using
(32) and (33) we obtain,
$$
\hbox{``Trace''} \ (U(h)) = \wh h (0) + \wh h (1) -
\sum_{{L(\chi,\rho)=0 \atop {\rm Re} \rho = {1\over 2}}} \wh h
(\chi,\rho) + \ify \, h(1) \leqno (34)
$$
where $\wh h (\chi,\rho)$ is defined by Corollary 2 and
$$
U(h) = \int_{C_k} U(g) \, h(g) \, d^* \, g \leqno (35)
$$
while the test function $h$ is in a suitable function space.
Note that the trace on the left hand side of (34) only makes
sense after a suitable regularisation since the left regular
representation of $C_k$ is not tra\c cable. This situation is
similar to the one encountered by Atiyah and Bott ([AB]) in
their proof of the Lefchetz formula. We shall first learn how
 to compute in a formal manner the above trace
from the fixed points of the action of $C_k$ on $X$. In section VII,
 we shall show how to regularize the trace and completely eliminate the 
parameter $\d$.
\vglue 1cm

\noindent {\bf IV The distribution trace formula for flows on
manifolds}

\medskip

In order to understand how the left hand side of III(34) should
be computed we shall first give an account of the proof of the usual 
Lefchetz formula by Atiyah-Bott ([AB])
and describe the  computation of the distribution
theoretic trace for flows on manifolds, which is a variation on the 
theme of [AB] and is due to Guillemin-Sternberg [GS]. 
We refer to Appendix III for a more detailed coordinate independent 
treatment following [GS].
\smallskip

Let us start with a diffeomorphism $\varphi$ of a smooth compact 
manifold $M$ 
and assume that the graph of $\varphi$ is transverse to the diagonal in 
$M 
\times M$. One can then easily define and compute the distribution 
theoretic 
trace of the operator $U : C^{\infty} (M) \rightarrow C^{\infty} (M)$,
$$
(U\xi) (x) = \xi (\varphi (x)) \, . \leqno (1)
$$
Indeed let $k(x,y)$ be the Schwartz distribution on $M \times M$ such 
that
$$
(U\xi) (x) = \int k(x,y) \, \xi (y) \, dy \, , \leqno (2)
$$
The distributional trace of $U$ is simply
$$
\hbox{``Trace''} \, (U) = \int k(x,x) \, dx \, , \leqno (3)
$$
\smallskip

Near the diagonal and in local coordinates one gets,
$$
k(x,y) = \delta (y - \varphi (x)) \leqno (4)
$$
where $\delta$ is the Dirac distribution.

Since, by hypothesis, the fixed points of $\varphi$ are isolated, one 
can 
compute the trace (3) as a finite sum $\build \sum_{x,\varphi (x) = 
x}^{}$ and get the 
contribution of each fixed point $x \in M $, $\varphi (x) = x$, as
$$
{1 \over \vert 1 - \varphi' (x)\vert} \  \leqno (5)
$$
where $\varphi' (x)$ is the Jacobian of $\varphi$ and $ \vert A \vert = 
\vert \det A \vert .$
One just uses the invertibility of ${\rm id} - \varphi' (x)$ to change 
variables in the 
integral,$$
\int \delta (y - \varphi (y)) \, dy \, . \leqno (6)
$$
One thus gets (cf. [AB]),
$$
\hbox{``Trace''} \, (U) = \sum_{x,\varphi (x) = x} {1 \over \vert 1 - 
\varphi' (x) \vert} \, . \leqno (7)
$$
This computation immediately extends to the action of $\varphi$ on 
sections of 
an equivariant vector bundle $E$ such as the bundle $\wdg^k T^*$ whose 
sections, $C^{\infty} (M,E)$ are the smooth forms of degree $k$. The 
alternate 
sum of the corresponding distribution theoretic traces is the {\it 
ordinary} 
trace of the action of $\varphi$ on the de Rham cohomology, thus 
yielding the 
usual Lefchetz formula,
$$
\sum (-1)^j \, \hbox{Trace} \, \varphi^* / H^j = \sum_{\varphi (x) = x} 
\, \hbox{sign} \, \det (1 - \varphi' (x)) \, . \leqno (8)
$$
Let us refer to the appendix for more pedantic notations which show that 
the distribution theoretic trace is coordinate independent.
\smallskip
We shall now write down the analogue of formula (7) in the case of a 
flow $F_t = 
\exp (tv)$ of diffeomorphisms of $M$, where $v \in C^{\infty} (M,T)$ is 
a vector 
field on $M$. We get a one parameter group of operators acting on 
$C^{\infty} 
(M)$,
$$
(U_t \, \xi) (x) = \xi (F_t (x)) \qquad \forall \, \xi \in C^{\infty} 
(M) \, , \ x \in M \, , \ t \in \Rb \, , \leqno (9)
$$

and we need the formula for,
$$
  \rho(h)=\hbox{``Trace''} \, \left( \int h(t) \, U_t \, dt \right) \ , 
\ h \in C_c^{\infty} (\Rb) \, , \ h(0) = 0 \, . \leqno (10)
$$
The condition $h(0) = 0$ is required because we cannot expect that the 
identity 
map $F_0$ is transverse to the diagonal.
In order to define $\rho $  as a distribution evaluated on the test 
function $h$, we let $f$ be the following map,
$$
f : X=M \ts \Rb \ra Y= M \ , \ f(x,t) = F_t (x) \, . \leqno (11)
$$
The graph of $f$ is the submanifold $Z$ of $X \ts Y $,
$$
Z = \{ (x,t,y) \ ; \ y = F_t (x)\} \, . \leqno (12)
$$
One lets $\vp$ be the diagonal map,
$$
\vp (x,t) = (x,t,x) \ , \ \vp : M \ts \Rb \ra X \ts Y
\leqno (13)
$$
and one assumes the transversality $\quad \vp \, {\cap \!\!\! ^{\mid}} 
\, Z \quad$  outside $\quad M\ts (0)$.

let $ \tau $ be the distribution,
$$
\tau = \vp^* (\d (y-F_t (x)) \, dy) \, , \leqno (14)
$$

\noindent and $q$  be the second projection,
$$
q(x,t) = t \in \Rb ,\leqno (15)
$$
then by definition $\rho $ is the pushforward $q_* (\tau)$ of the 
distribution
$\tau $.

\smallskip

\noindent One checks (cf. Appendix III) that $q_* (\tau)$ is a 
generalized function.

\smallskip
 Exactly as in the case of a single 
transformation, the contributions to (10) will come from the fixed 
points of 
$F_t$. The latter will come either from a {\it zero} of the vector field 
$v$, 
(i.e. $x \in M$ such that $v_x = 0$) or from a {\it periodic orbit} 
$\gamma$ of 
the flow and we call $T_{\gamma}^{\#}$ the length of such a periodic 
orbit. Under 
the above transversality 
hypothesis the formula for (10) is (cf. [GS], [G] and the Appendix III),
$$
\ldisplaylinesno{
\hbox{``Trace''} \, \left( \int h(t) \, U_t \, dt \right) = &(16) \cr
\sum_{x,v_x = 0} \ \int {h(t) \over \vert 1 - (F_t)_* \vert} \, dt + 
\sum_{\gamma} \sum_{T} T_{\gamma}^{\#} \, {1 \over \vert 1 - (F_{T /})_* 
\vert} 
\, h(T) \cr
}
$$

where in the second sum $\gamma$ is a periodic orbit with length 
$T_{\gamma}^{\#}$, and $T$ varies in $\Zb \, T_{\gamma}^{\#}$ while 
$(F_{T /})_* $
 is the Poincare return map, i.e. the restriction of the tangent map to 
the transversal of the orbit.
\smallskip

One can rewrite (16) in a better way as,
$$
\hbox{``Trace''} \, \left( \int h(t) \, U_t \, dt \right) = 
\sum_{\gamma} 
\int_{I_{\gamma}} {h(u) \over \vert 1 - (F_u)_* \vert} \, d^* u \, , 
\leqno (17)
$$
where the zeros $x \in M $, $v_x = 0$, are considered also as periodic 
orbits $\gamma$, 
while $I_{\gamma} \subset \Rb$ is the isotopy subgroup of any $x \in 
\gamma$, 
and $d^* u$ is the unique Haar measure in $I_{\gamma}$ such that the 
covolume of 
$I_{\gamma}$ is equal to $1$, i.e. such that for the unique Haar measure  
$d\mu $ of total
 mass $1$ on $\Rb / I$ and any $f \in C_c^{\infty} (\Rb)$,
$$
\int_{\Rb} f(t) \, dt = \int_{\Rb / I} \left( \int_I f(u+s) \, d^* u 
\right) \, d\mu (s) \, , \leqno (18)
$$
\smallskip
Also we still write $(F_u)_*$ for the restriction of the tangent map to 
$F_u$ to the transverse space of the orbits.

\noindent To understand what $(F_t)_*$ looks like at a zero of $v$ we 
can replace
$v(x)$ for $x$ near $x_0$ by its tangent map. For simplicity we
take the one dimensional case, with $v(x) = x \, {\part \over
\part x}$, acting on $\Rb = M$.

\smallskip

\noindent One has $F_t (x) = e^t \, x$ 
 . Since $F_t$ is linear the
tangent map $(F_t)_*$ is
$$
(F_t)_* = e^t \leqno (19)
$$
and (12) becomes
$$
\hbox{``Trace''} \, \left( \int h(t) \, U_t \, dt \right) =  \int {h(t) 
\over \vert 1 - e^t \vert} \, dt \, , \leqno (20)
$$
Thus for this flow the distribution trace formula is
$$
\hbox{``Trace''} \, (U(h)) = \int {h(u) \over \vert 1-u \vert}
\ d^* u \leqno (21)
$$
where we used the multiplicative notation so that $\Rb_+^*$
acts on $\Rb$ by multiplication,while $U(h)=\int { U(v)h(v) }\ d^* v $ 
and  $d^* v $ is the haar measure of the group $\Rb_+^*$.

One can treat in a similar way the action, by multiplication, of the 
group of non zero complex numbers on the manifold $\Cb $.
\smallskip

\noindent We shall now investigate the more general case of an
arbitrary local field.

\vglue 1cm

\noindent {\bf V The action $(\lb ,x) \ra \lb \, x$ of $K^*$ on a local 
field 
$K$.}

\medskip

We let $K$ be a local field and consider the map,
$$
f:K \ts K^* \ra K \ , \ f(x,\lb) = \lb \, x \leqno (1)
$$
together with the diagonal map,
$$
\vp : K \ts K^* \ra K \ts K^* \ts K \ , \ \vp (x, \lb) =
(x,\lb ,x) \leqno (2)
$$
as in IV (11) and (12) above.

\smallskip

\noindent When $K$ is Archimedian we are in the framework of
manifolds and we can associate to $f$ a $\d$-section with
support $Z = \ \hbox{Graph} \ (f)$,
$$
\d_Z = \d (y-\lb \, x) \, dy \, . \leqno (3)
$$
Using the projection $q (x,\lb) =\lb$ from $K \ts K^*$ to
$K^*$ we then consider as above the generalized function on
$K^*$ given by,
$$
q_* (\vp^* \, \d_Z) \, . \leqno (4)
$$
The formal computation of this generalized function of $\lb $ is
$$
\int \d (x-\lb \, x) \, dx = \int \d ((1-\lb) x) \, dx = \int
\d (y) \, d((1-\lb)^{-1} \, y)
$$
$$
= \vert 1-\lb \vert^{-1} \int \d (y) \, dy = \vert 1-\lb
\vert^{-1} \, .
$$
We want to justify it by computing the convolution of the
Fourier transforms of $\d (x-y)$ and $\d (y-\lb \, x)$ since this is the 
correct
 way of defining the product of two distributions in this local context
. Let us
first compute the Fourier transform of $\d (ax + by)$ where
$(a,b) \in K^2 (\not= 0)$. The pairing between $K^2$ and its
dual $K^2$ is given by
$$
\lgl (x,y) , (\xi ,\eta) \rgl = \a (x \, \xi + y \, \eta) \in
U(1) \, . \leqno (5)
$$
where $\a$ is a fixed nontrivial character of the additive group $K$.

\smallskip

\noindent Let $(c,d) \in K^2$ be such that $ad-bc =1$ and
consider the linear invertible transformation of $K^2$,

$$
L \left[ \matrix{ x \cr y \cr} \right] = \left[ \matrix{ a &b
\cr c &d \cr} \right] \, \left[ \matrix{x \cr y \cr} \right]
\, .\leqno (6)
$$
The Fourier transform of $\vp \circ L$ is given by
$$
(\vp \circ L)^\wedge = \vert \det L\vert^{-1} \, \wh{\vp}
\circ (L^{-1})^t \, . \leqno (7)
$$
Here one has $\det L =1$ and $(L^{-1})^t$ is
$$
(L^{-1})^t = \left[ \matrix{ d &-c \cr -b &a \cr} \right] \, .
\leqno (8)
$$

\smallskip

\noindent One first computes the Fourier transform of $\d
(x)$, the additive Haar measure $dx$ is normalized so as to
be selfdual, and in one variable $\d (x)$ and 1 are Fourier 
transforms of each other, thus
$$
(\d \ot 1)^\wedge = 1\ot \d \, . \leqno (9)
$$
Using (7) one gets that the Fourier transform of $\d (ax
+ by)$ is $\d (-b \, \xi + a \, \eta)$. Thus we have to
compute the convolution of the two generalized functions, $\d
(\xi + \eta)$ and $\d (\xi + \lb \, \eta)$. Now
$$
\int f(\xi ,\eta) \, \d (\xi + \eta) \, d\xi \, d\eta = \int
f(\xi ,-\xi ) \, d\xi
$$
and
$$
\int f (\xi ,\eta) \, \d (\xi +\lb \, \eta) \, d\xi \, d\eta =
\int f (-\lb \, \eta , \eta) \, d\eta
$$
thus we are dealing with two measures carried respectively by
two distinct lines. Their convolution evaluated on $f \in
C_c^{\ify} (K^2)$ is $\int f(\a + \b) \, d\mu (\a) \, d\nu
(\b) = \int \int f((\xi ,-\xi) + (-\lb \, \eta ,\eta)) \,
d\xi \, d\eta = \int \int f(\xi - \lb \, \eta , -\xi + \eta)
\, d\xi \, d\eta = \Bigl( \int \int f(\xi' , \eta')$ $d\xi'
\, d\eta' \Bigl) \ts \vert J \vert^{-1}$ where $J$ is the
determinant of the matrix $\left[ \matrix{ 1&-\lb \cr -1&1
\cr} \right] = L$, so that $\left[ \matrix{ \xi' \cr \eta'
\cr} \right] = J \, \left[ \matrix{ \xi \cr \eta \cr}
\right]$. One has $J = 1-\lb$ and thus the convolution of the
generalized functions $\d (\xi + \eta)$ and $\d (\xi + \lb \,
\eta)$ gives as expected the constant function
$$
\vert 1-\lb\vert^{-1} \, 1 \, . \leqno (10)
$$
Correspondingly, the product of the distribution $\d (x-y)$
and $\d (y-\lb \, x)$ gives $\vert 1-\lb \vert^{-1} \, \d_0$
so that,
$$
\int \d (x-y) \, \d (y- \lb \, x) \, dx \, dy = \vert
1-\lb\vert^{-1} \, . \leqno (11)
$$

\smallskip

\noindent In this local case the Fourier transform alone was
sufficient to make sense of the relevant product of
distributions. In fact this would continue to make sense if
we replace $\d (y-\lb \, x)$ by $\int h(\lb^{-1}) \, \d
(y-\lb \, x) \, d^* \, \lb$ where $h(1)=0$ .
\smallskip
We shall now treat in detail the more delicate general case where $h(1)$ 
is arbitrary.

We shall prove a precise general result (theorem 3) which handles 
the lack of transversality when $h(1) \not = 0$. We deal directly with 
the 
following operator in $L^2 (K)$,
$$
U(h) = \int h (\lb) \, U (\lb) \, d^* \lb \, , \leqno (12)
$$
where the scaling operator $U (\lb)$ is defined by
$$
(U (\lb) \, \xi) (x) = \xi (\lb^{-1} \, x) \qquad \fl \, x \in K 
\leqno (13)
$$
and where the multiplicative Haar measure $d^* \lb$ is normalized by,
$$
\int_{\vert \lb \vert \in [1, \L]} d^* \lb \sim \log \L \qquad 
\hbox{when} \ \L \ra \ify \, . \leqno (14)
$$
To understand the ``trace'' of $U(h)$ we shall proceed as in the 
Selberg trace formula ([Se]) and use a cutoff. For this we use the 
orthogonal projection $P_{\L}$ onto the subspace,
$$
P_{\L} = \{ \xi \in L^2 (K) \, ; \ \xi (x) = 0 \qquad \fl x \ , \ 
\vert x \vert > \L \} \, . \leqno (15)
$$

Thus, $P_{\L}$ is the 
multiplication operator by the function $\rho_{\L}$, where $\rho_{\L} 
(x) = 1$ if $\vert x \vert \leq \L$, and $\rho (x) = 0$ for $\vert x 
\vert > \L$. 
This gives 
an infrared cutoff and to get an ultraviolet cutoff we use 
$\wh{P}_{\L} = F P_{\L} F^{-1}$ where $F$ is the Fourier transform 
(which depends upon the basic character $\a$). We let
$$
R_{\L} = \wh{P}_{\L} \, P_{\L} \, . \leqno (16)
$$
The main result of this section is then,

\medskip

\noindent {\bf Theorem 3.} {\it Let $K$ be a local field with basic 
character $\a$. Let $h \in \Sc (K^*)$ have compact support. Then
$R_{\L} \, U(h)$ is a trace class operator and when 
 $\L \ra \ify$, one has 
$$
{\rm Trace} \, (R_{\L} \, U(h)) = 2h (1) \log' \L + \int' {h(u^{-1}) 
\over \vert 1-u \vert} \, d^*u + o(1)
$$
where $2 \log' \L = \int_{\lb \in K^*, \, \vert \lb \vert \in [\L^{-1}, 
\L]} d^* \lb$,
 and the principal value $\int'$ is uniquely determined by the pairing 
with the 
unique distribution on $K$ which agrees with ${du \over \vert 1-u 
\vert}$ for $u \not= 1$ and whose Fourier transform vanishes at $1$.
}

\medskip
\noindent {\it Proof.} We normalize as above the additive 
Haar measure to be the selfdual one on $K$. Let the constant $\rho > 0$ 
be determined by the equality,
$$
\int_{1 \leq \vert \lb \vert \leq \L} \, {d\lb \over \vert \lb \vert} 
\sim \rho \log \L \qquad \hbox{when} \ \L \ra \ify \, . \leqno (17)
$$
so that $d^* \lb = \, \rho^{-1} {d\lb \over \vert \lb \vert} $. 
Let $L$ be the unique distribution, extension of $ \rho^{-1}{du \over 
\vert 1-u \vert}$ whose 
Fourier transform vanishes at 1, $\wh L (1) = 0$. One then has by 
definition,
$$
\int' {h(u^{-1}) 
\over \vert 1-u \vert} \, d^* u =\,\left\lgl L ,  \, {h(u^{-1}) \over 
\vert u \vert} \right\rgl 
\, , \leqno (18)
$$
where ${h(u^{-1}) \over \vert u \vert}=0$ for $u^{-1}$ outside the 
support of $h$.
\smallskip 
Let $T = U(h)$. We can write the Schwartz 
kernel of $T$ as,
$$
k(x,y) = \int h(\lb^{-1}) \, \d (y - \lb x) \, d^* \lb \, . \leqno 
(19)
$$
Given any such kernel $k$ we introduce its symbol,
$$
\s (x,\xi) = \int k(x,x+u) \, \a (u\xi) \, du \leqno (20)
$$
as its partial Fourier transform.  The Schwartz kernel 
$r_{\L}^t (x,y)$ of the transpose $R_{\L}^t$ is given by,
$$
r_{\L}^t (x,y) = \rho_{\L} (x) \, (\wh{\rho_{\L}}) \, (x-y) \, . 
\leqno (21)
$$
Thus, the symbol $\s_{\L}$ of $R_{\L}^t$ is simply,
$$
\s_{\L} (x,\xi) = \rho_{\L} (x) \, \rho_{\L} (\xi) \, . \leqno (22)
$$
The operator $R_{\L}$ is of trace class and one has,
$$
{\rm Trace} \, (R_{\L} \, T) = \int k (x,y) \, r_{\L}^t (x,y) \, dx 
\, dy \, . \leqno (23)
$$
Using the Parseval formula we thus get,
$$
{\rm Trace} \, (R_{\L} \, T) = \int_{\vert x \vert \leq \L , \vert 
\xi \vert \leq \L} \s (x,\xi) \, dx \, d\xi \, . \leqno (24)
$$
Now the symbol $\s$ of $T$ is given by,
$$
\s (x,\xi) = \int h(\lb^{-1}) \left( \int \d (x+u-\lb x) \, \a (u\xi) 
\, du \right) \, d^* \lb \, . \leqno (25)
$$
One has,
$$
\int \d (x+u-\lb x) \, \a (u\xi) \, du = \a ((\lb - 1) \, x \xi) \, , 
\leqno (26)
$$
thus (25) gives,
$$
\s (x,\xi) = \rho^{-1} \int_K g(\lb) \, \a (\lb x \xi) \, d\lb \leqno 
(27)
$$
where,
$$
g(\lb) = h ((\lb + 1)^{-1}) \, \vert \lb + 1 \vert^{-1} \, . \leqno 
(28)
$$
Since $h$ is smooth with compact support on $K^*$ the function $g$ 
belongs to $C_c^{\ify} (K)$.

\smallskip

Thus $\s (x,\xi) = \rho^{-1} \, \wh g (x \xi)$ and,
$$
{\rm Trace} \, (R_{\L} \, T) = \rho^{-1} \int_{\vert x \vert \leq \L 
, \vert \xi \vert \leq \L} \wh g (x\xi) \, dx \, d\xi \, . \leqno 
(29)
$$
With $u = x \xi$ one has $dx \, d\xi = du \, {dx \over \vert x 
\vert}$ and, for $\vert u \vert \leq \L^2$,
$$
\rho^{-1} \int_{{\vert u 
\vert \over \L } \leq \vert x \vert  \leq \L  } \, 
{dx \over \vert x \vert} = 2 \log' \L - \log \vert u \vert \leqno (30)
$$
(using the precise definition of $\log' \L$ to handle the boundary 
terms). Thus we can rewrite (29) as,

$$
{\rm Trace} \, (R_{\L} \, T) = \int_{\vert u \vert \leq \L^2} \, \wh 
g (u) \, (2\log' \L - \log \vert u \vert) \, du \, \leqno (31)
$$
Since $g \in C_c^{\ify} (K)$ one has,
$$
\int_{\vert u \vert \geq \L^2} \, \vert \wh g (u) \vert \, du = O 
(\L^{-N}) \qquad \fl \, N \leqno (32)
$$
and similarly for $\vert \wh g (u) \, \log \vert u \vert \vert$. Thus
$$
{\rm Trace} \, (R_{\L} \, T) = 2 \, g(0) \, \log' \L - \int \wh g (u) 
\, \log \vert u \vert \, du \, + o(1) . \leqno (33)
$$
Now for any local field $K$ and basic character $\a$, if we take for 
the Haar measure $da$ the selfdual one, the Fourier transform of the 
distribution $\vp (u) = - \log \vert u \vert$ is given outside 0 by
$$
\wh{\vp} (a) = \rho^{-1} \, {1 \over \vert a \vert} \, , \leqno (34)
$$
with $\rho$ determined by (17). To see this one lets $P$ be the 
distribution on $K$ given by,
$$
P (f) = \lim_{\ve \ra 0 \atop \ve \in {\rm Mod} (K)} \left( 
\int_{\vert x \vert \geq \ve} f(x) \, d^* x + f(0) \log \ve \right) 
\, . \leqno (35)
$$
One has $P(f_a) = P(f) - \log \vert a \vert \, f(0)$ which is enough 
to show that the function $\wh P (x)$ is equal to $-\log \vert x 
\vert + $ cst, and $\wh{\vp}$ differs from $P$ by a multiple of 
$\d_0$.
\smallskip
 Thus the Parseval formula gives, with the convention of 
theorem 3,
$$
-\int \wh g (u) \, \log \vert u \vert \, du = {1 \over \rho} \int' 
g(a) \, {da \over \vert a \vert} \, . \leqno (36)
$$
Replacing $a$ by $\lb - 1$ and applying (28) gives the desired 
result.~\xx

\medskip

We shall show in appendix II that the privileged principal value, 
which depends upon the basic character $\a$, is the same as in Weil's 
explicit formulas.
 \vglue 1cm

\noindent {\bf VI The global case, and the formal trace computation.}

\medskip

We shall now consider the action of $C_k$ on $X$ and write
down the analogue of IV (17) for the distribution trace
formula.
\smallskip
Both $X$ and $C_k$ are defined as quotients and we let

$$\pi : A \to X\, ,\ c :{\rm GL}_1 (A) \to C_k \leqno (1)$$
be the corresponding quotient maps.
\smallskip
As above we consider the graph $Z$ of the action
$$ f : X \times C_k \to X\, , \ f(x,\lb) = \lb x\leqno (2)$$
and the diagonal map
$$ \vp : X \times  C_k \to X \times C_k \times X\, \qquad
\vp(x,\lb) = (x,\lb,x)\,. \leqno (3)$$

We first investigate the fixed points, $\vp^{-1}(Z)$, i.e.~the pairs
$(x,\lb) \in X \times C_k$ such that $\lb x = x$. Let $x =
\pi (\tilde x)$ and $\lb = c(j)$. Then the equality $\lb x =
x$ means that $\pi (j \tilde x) = \pi(\tilde x)$ thus
there exists $q \in k^*$ such that with $ \tilde j = qj$, one has
$$ \tilde j \tilde x = \tilde x\,  . \leqno(4)
$$
Recall now that $A$ is the restricted direct product $A =
\build{\Pi}_{\rm res}^{}\, k_v$ of the local fields $k_v$
obtained by completion of $k$ with respect to the place $v$.
The equality (4) means that $\tilde j_v \tilde x_v = \tilde
x_v$, thus, if $\tilde x_v \not = 0$ for all $v$ it follows
that $\tilde j_v = 1\ \forall v$ and $\tilde j = 1$. This
shows that the projection of $\vp^{-1} (Z) \cap C_k
\backslash \{1\}$ on $X$ is the union of the hyperplanes
$$\cup H_v\, ;\ H_v = \pi (\tilde H_v)\, ,\ \tilde H_v = \{ x \, ;\ x_v 
=
0\}\, . \leqno (5)$$
Each $\tilde H_v$ is closed in $A$ and is invariant under
multiplication by elements of $k^*$. Thus each $H_v$ is a closed
subset of $X$ and one checks that it is the closure of the orbit under
$C_k$ of any of its generic points
$$x\, ,\ x_u = 0 \quad \iff \quad u = v \, .\leqno
(6)  $$ For any such point $x$, the isotropy group $I_x$ is the image in
$C_k$ of the multiplicative group $k^*_v$,
$$ I_x = k^*_v \leqno (7)$$
by the map $\lb \in k^*_v \to (1,\ldots,1,\lb,1,\ldots)$. This map
already occurs in class field theory (cf [W1]) to relate the local 
Galois
theory to the global one.
\smallskip
Both groups $k^*_v$ and $C_k$ are commensurable to $\Rb^*_+$ by the
module homomorphism, which is proper with cocompact range,
$$ G \build{\longrightarrow}_{}^{|\  |} \Rb^*_+\,. \leqno (8)$$
Since the restriction to $k^*_v$ of the module of $C_k$ is the module of
$k^*_v$, it follows that
$$I_x \ \hbox{\rm is a cocompact subgroup of}\ C_k\, . \leqno (9) $$
This allows to normalize the respective Haar measures in such a way
that the covolume of $I_x$ is 1. This is in fact insured by the
canonical normalisation of the Haar measures of modulated groups ([W3
]),
$$
\int_{\vert g \vert \in [1,\L]} \, d^* g \sim \log \L \quad
\hbox{when} \quad  \L \ra +\ify \, . \leqno (10)
$$
It is important to note that though $I_x$ is cocompact in $C_k$, the
orbit of $x$ is not closed and one needs to close it, the result being
$H_v $. We shall learn how to justify this point later in section VII, 
in the similar situation of the action 
of $C_S$ on $X_S$.
We can now in view of the results of the two preceding sections, write
down the contribution of each $H_v$ to the distributional trace;
\medskip
Since $\tilde H_v$ is a hyperplane, we can identify the transverse
space $N_x$ to $H_v$ at $x$ with the quotient
$$N_x = A/\tilde H_v = k_v \leqno (11)$$
namely the additive group of the local field $k_v$. Given $j \in I_x$
one has $j_u = 1 \ \forall u \not = v$, and $j_v = \lb \in k^*_v$. The
action of $j$ on $A$ is linear and fixes $x$, thus the action on the
transverse space $N_x$ is given by
$$ (\lb,a) \to \lb a \,\quad \forall a\in k_v. \leqno (12)$$
We can thus proceed with some faith and write down the contribution of
$H_v$ to the distributional trace in the form,
$$\int_{k^*_{v}} \, {h(\lb) \over |1-\lb|}\, d^*\lb \leqno (13)$$
where $h$ is a test function on $C_k$ which vanishes at 1. We now have
to take care of a discrepancy in notation with the third section
(formula 9), where we used the symbol $U(j)$ for the operation
$$\bigl (U(j)f \bigr) (x) = f(j^{-1} x) \leqno (14)$$
whereas we use $j$ in the above discussion. This amounts to replace the
test function $h(u)$ by $h(u^{-1})$ and we thus obtain as a formal
analogue of III(17) the following expression for the distributional 
trace
  $$ \hbox{\rm ``Trace"}\, (U(h)) = \sum_{v} \int_{k^*_v} \, {h(u^{-1})
\over |1-u|}\, d^* u \, . \leqno (15)$$
Now the right-hand side of (15) is, when restricted to the hyperplane
$h(1)=0$, the distribution obtained by Andr\'e Weil [W3] as the 
synthesis of
the explicit formulas of number theory for all $L$-functions with
Gr\"ossencharakter. In particular we can rewrite it as
$$\hat h(0) + \hat h(1) - \sum_{L(\chi,\rho)=0} \hat h (\chi,\rho) +
\infty\  h(1) \leqno (16)$$
where this time the restriction Re$(\rho) = {1 \over 2}$ has been
eliminated.
\smallskip
Thus, equating (34) of section III and (16) for $h(1)=0$ would yield the 
desired
information on the zeros. Of course, this does require first eliminating 
the role
 of $\d$, and (as in [AB]) to prove that
the distributional trace coincides with the ordinary operator theoretic
trace on the cokernel of $E$. This is achieved for the usual set-up of
the Lefchetz fixed point theorem by the use of families.
\smallskip
A very important property of the right hand side of (15) (and of IV (17) 
in
general) is that if the test function $h, h(1)=0 $ is positive,
$$ h(u) \ge 0 \quad \forall \ u \in C_k \leqno (17)$$
then the right-hand side is {\it positive\/}. This indicated from the 
very
start that in order to obtain the Polya-Hilbert space from the Riemann
flow, it is {\it not\/} quantization that should be involved but simply
the passage to the $L^2$ space, $X \to L^2(X)$. Indeed the positivity of 
IV
(17) is typical of {\it permutation matrices\/} rather than of
quantization. This distinction plays a crucial role in the above 
discussion of the trace formula,
 in particular the expected trace formula is not a semi-classical 
formula but a Lefchetz formula 
in the spirit of [AB]. 
\smallskip

The above discussion is {\it not\/} a rigorous justification of
this formula. The first obvious obstacle is that the distributional 
trace
 is only formal and to give it a rigorous meaning tied up to Hilbert 
space operators, one needs as in section V, to perform a cutoff.
 The second difficulty comes from the presence of the parameter
 $\d$ as a label for the 
Hilbert space, while $\d$ does not appear in the trace formula.
 As we shall see in the next two sections the cutoff will 
completely eliminate the role of $\d$, and we shall
 nevertheless show (by proving positivity of the Weil distribution)
 that the validity of the ($\d$ independent) trace formula
 is equivalent to the Riemann Hypothesis for all Gr\"ossencharakters of 
$k$.
 
\vglue 1cm
\noindent {\bf VII Proof of the trace formula in the $S$-local case.}

\medskip
In the formal trace computation of section VI, we skiped over
 the difficulties
inherent to the tricky structure of the space $X$.
 In order to understand how to handle trace formulas on such spaces 
we shall consider the slightly simpler situation which arises when
 one only considers a finite set $S$ of places of $k$.
As soon as the cardinality of $S$ is larger than $3$,
 the corresponding space $X_S$ does share most of the tricky features
of the space $X$. In particular it is no longer of type I
 in the sense of Noncommutative Geometry.
\smallskip
 We shall nevertheless prove a precise general result (theorem 4)
 which shows that the above handling of periodic orbits and
 of their contribution to the trace is the correct one.
 It will in particular show why the orbit of the fixed
 point $0$, or of elements $x \in A$, such that $x_v$ vanishes for at 
least two places 
 do not contribute to the trace formula.
\smallskip
At the same time, we shall handle as in section V, 
the lack of transversality when $h(1) \not = 0$.
\smallskip
Let us first describe the reduced framework for the trace formula.
 We let $k$ be a global field and $S$ a finite set of places of $k$
 containing all infinite places. The group $O_S^*$ of $S$-units
 is defined as the subgroup of $k^*$, 
$$
O_S^*=\{q \in k^*, \vert q_v \vert =1, v \notin S \} \leqno (1)
$$
It is cocompact in $J^1_S$ where,
$$
J_S=\prod_{v\in S} k^*_v \leqno (2)
$$
and,
$$
J_S^1=\{j \in J_S, \vert j \vert=1 \}.\leqno (3)
$$
Thus the quotient group $C_S=J_S/O_S^*$ plays the same role as $C_k$, 
and acts on the quotient
 $X_S$ of $A_S=\prod_{v\in S} k_v$ by $O_S^*$.
\smallskip
 To keep in mind a simple example, one can take $k= \Qb$,
 while $S$ consists of the three places $2$, $3$, and $\ify$.
 One checks in this example that the topology of $X_S$ is not
 of type I since for instance the group 
$O_S^*= \{ \pm 2^n3^m; \,   n, m \in \Zb \}$
 acts ergodically on $ \{0 \} \ts \Rb \sbs A_S$.
\smallskip
We normalize the multiplicative Haar measure $d^* \lb$ of $C_S$ by,
$$
\int_{\vert \lb \vert \in [1, \L]} d^* \lb \sim \log \L \qquad 
\hbox{when} \ \L \ra \ify \, , \leqno (4)
$$ 
and normalize the multiplicative Haar measure $d^* \lb$ of $J_S$
 so that it agrees with the above on 
a fundamental domain $D$ for the action of $O_S^*$ on $J_S$.
\smallskip
There is no difficulty in defining the Hilbert space $L^2 (X_S)$
 of square integrable functions on $X_S$.
 We proceed as in section III (without the $\d$),
 and complete (and separate) the Schwartz space $\Sc(A_S)$
 for the pre-Hilbert structure given by,
$$
\Vert f \Vert^2 = \int \Bigl\vert \sum_{q\in O_S^*} f (qx)
\Bigl\vert^2 \,  \vert x \vert \,
d^* x \leqno (5)
$$
where the integral is performed on $C_S$ or equivalently
 on a fundamental domain $D$ for the action of $O_S^*$ on $J_S$.
To show that (5) makes sense, one proves that for $f \in \Sc(A_S)$,
 the function $E_0(f)(x)= \sum_{q\in O_S^*} f (qx)$
is bounded above by a power of $ Log \vert x \vert$ when $\vert x \vert$
 tends to zero. To see this when $f$ is the characteristic function of
 $\{x \in A_S, \vert x_v \vert \leq 1 \, , \fl v \in S \}$,
 one uses the cocompactness of $O_S^*$ in $J^1_S$, to replace
 the sum by an integral. The latter is then comparable to, 
$$
 \int_{u_i \geq 0, \sum u_i = -Log  \vert x \vert} \prod du_i , \leqno 
(6)
$$
where the index $i$ varies in $S$. The general case follows. 
\smallskip
The scaling operator $U (\lb)$ is defined by,
$$
(U (\lb) \, \xi) (x) = \xi (\lb^{-1} \, x) \qquad \fl \, x \in A_S
\leqno (7)
$$
and the same formula, with $x \in X_S$ defines its action on $L^2 
(X_S)$.
Given a smooth compactly supported function 
$h$ on $C_S$, $U(h)= \int h(g)U(g) dg $ makes sense as an operator 
acting on $L^2 (X_S)$.
\smallskip
We shall now show that the Fourier transform $F$ on
 $\Sc(A_S)$ does extend to a unitary operator
 on the Hilbert space $L^2 (X_S)$.
\medskip
\noindent {\bf Lemma 1.} {\it a) For any $f_i \in \Sc(A_S)$ 
the series $\sum_{O_S^*}  \left< f_1 , U(q)\, f_2 \right>_A $
of inner products in  $L^2 (A_S)$
 converges geometrically on the abelian finitely generated group 
$O_S^*$.
 Moreover its sum is equal to the inner product of $f_1$
 and $f_2$ in the Hilbert space $L^2 (X_S)$.}
\smallskip
{\it b) Let $\a =\prod \a_v$ be a  basic character of
 the additive group $A_S$ and $F$ the corresponding Fourier 
transformation.
 The map $ f \rightarrow F(f)$, $ f \in \Sc(A_S)$
 extends uniquely to a unitary operator in
 the Hilbert space $L^2 (X_S)$.} 
\medskip

\noindent {\it Proof.} The map $L : O_S^* \ra \Rb^S $, given by 
$L(u)_v= \, \log \vert u_v \vert$
has a finite kernel and its range is a lattice in the
 hyperplane $H= \, \{(y_v), \sum_S y_v = 0 \}$. On $H$ one has
$ Sup_S y_v \geq 1/2n \sum \vert y_v \vert $, where $n= card(S)$.
 Thus one has the inequality,
$$
 Sup_S \vert q_v \vert \geq exp(d(q,1)  \qquad \fl q \in  O_S^* \leqno 
(8)
$$
for a suitable word metric $d$ on $O_S^*$.
\smallskip
Let $K_n =\,  \{ x \in A_S ; \, \vert x_v \vert  \leq n ,\,\,\,
\fl v \in S \}$ and $k_n$ be the characteristic function of $K_n $.
 Let $( \lb_n)$ be a sequence of rapid decay such that,
$$
 \vert f_i(x) \vert  \leq  \sum  \lb_n \, k_n(x) \qquad \fl x \in A_S .  
\leqno (9)
$$
 One has for a suitable constant $c$, 
$$
 \vert \left< k_n , U(q^{-1})\, k_n \right> \vert \leq c \,
 n^m (Sup_S \vert q_v \vert)^{-1}  \leqno (10)
$$
where $m= \, Card(S)$.
\smallskip
Using (9) we thus see that $ \left< f_1 , U(q)\, f_2 \right>_A$
 decays exponentially on $O_S^*$. Applying Fubini's theorem
yields the equality,
$$
\int \Bigl\vert \sum_{q\in O_S^*} f (qx)
\Bigl\vert^2 \,  \vert x \vert \,
d^* x = \,    \sum_{O_S^*}  \left< f , U(q)\, f \right>_A .\leqno (11)
$$
This proves a). To prove b), one
 just uses (11) and the equalities $\left< Ff , \, Ff \right>_A =
\, \left< f , \, f \right>_A$ and $F( U(q)\, f)= \, U(q^{-1})F(f)$.~\xx

\smallskip 

Now exactly as above for the case of local fields (theorem V.3),
 we need to use a cutoff. For this we use the 
orthogonal projection $P_{\L}$ onto the subspace,
$$
P_{\L} = \{ \xi \in L^2 (X_S) \, ; \ \xi (x) = 0 \qquad \fl x \ , \ 
\vert x \vert > \L \} \, . \leqno (12)
$$
Thus, $P_{\L}$ is the 
multiplication operator by the function $\rho_{\L}$, where $\rho_{\L} 
(x) = 1$ if $\vert x \vert \leq \L$, and $\rho (x) = 0$
 for $\vert x \vert > \L$. This gives 
an infrared cutoff and to get an ultraviolet cutoff we use 
$\wh{P}_{\L} = F P_{\L} F^{-1}$ where $F$ is the Fourier transform 
(lemma 1)
which depends upon the choice of the basic character $\a =\prod \a_v$. 
We let
$$
R_{\L} = \wh{P}_{\L} \, P_{\L} \, . \leqno (13)
$$
The main result of this section is then,

\medskip

\noindent {\bf Theorem 4.} {\it Let $A_S$ be as above, with basic 
character $\a= \prod \a_v$. Let $h \in \Sc (C_S)$ have compact support. 
Then
  when 
 $\L \ra \ify$, one has
$$
{\rm Trace} \, (R_{\L} \, U(h)) = 2h (1) \log' \L + \sum_{v \in S} 
\int'_{k^*_v} {h(u^{-1}) 
\over \vert 1-u \vert} \, d^* u + o(1)
$$
where $2 \log' \L = \int_{\lb \in C_S, \, \vert \lb \vert \in
 [\L^{-1}, \L]} d^* \lb$, each $k^*_v$ is embedded in $C_S$ 
by the map $u \rightarrow (1,1,...,u,...,1)$ and the principal
 value $\int'$ is uniquely determined by the pairing with the 
unique distribution on $k_v$ which agrees with ${du \over \vert 1-u 
\vert}$ for $u \not= 1$ and whose Fourier transform relative to $\a_v$ 
vanishes at $1$.
}

\medskip

\noindent {\it Proof.} We normalize as above the additive 
Haar measure $dx$ to be the selfdual one on the abelian group $A_S$.
 Let the constant $\rho > 0$ be determined by the equality,
 (where the fundamental domain $D$ is as above),
$$
\int_{\lb \in D, \, 1 \leq \vert \lb \vert \leq \L} \, {d\lb \over \vert 
\lb \vert} 
\sim \rho \log \L \qquad \hbox{when} \ \L \ra \ify \, . 
$$
so that $d^* \lb = \, \rho^{-1} {d\lb \over \vert \lb \vert} $. 
\smallskip 
We let $f$ be a smooth compactly supported function on $J_S$
such that
$$
\sum_{q\in O_S^*} f(qg)=h(g)
\qquad \fl \, g \in C_S .\leqno (14)
$$
The existence of such an $f$ follows from the discreteness of $O_S^*$ in 
$J_S$. We then
 have the equality $ U(f) =\, U(h)$, where
$$
U(f) = \int f (\lb) \, U (\lb) \, d^* \lb \, , \leqno (15)
$$
To compute the trace of $U(h)$ acting on functions on the quotient space 
$X_S$, 
we shall proceed as in the 
Selberg trace formula ([Se]). Thus for an operator $T$, acting
 on functions on $A_S$, which commutes with the action of $O_S^*$
 and is represented by an integral kernel,
$$
T(\xi) = \int k(x,y)\xi(y) \,dy , \leqno (16)
$$
the trace of its action on $L^2 (X_S)$ is given by,
$$
Tr(T)= \sum_{q\in O_S^*}\int_{D} k(x,qx)dx . \leqno (17)
$$
where $D$ is as above a fundamental domain for the action of $O_S^*$ on 
the subset
$J_S$ of $A_S$, whose complement is negligible.
Let $T = U(f)$. We can write the Schwartz 
kernel of $T$ as,
$$
k(x,y) = \int f(\lb^{-1}) \, \d (y - \lb x) \, d^* \lb \, . \leqno 
(18)
$$
by construction one has,
$$
k(qx,qy)=k(x,y)\,\qquad q \in O_S^*.\leqno (19)
$$
For any $q \in O_S^*$, we shall evaluate the integral,
$$
I_q=\int_{x \in D} \,  k(qx,y) r_{\L}^t (x,y) dy dx \leqno (20)
$$
where the Schwartz kernel 
$r_{\L}^t (x,y)$ for the transpose $R_{\L}^t$ is given by,
$$
r_{\L}^t (x,y) = \rho_{\L} (x) \, (\wh{\rho_{\L}}) \, (x-y) \, . 
\leqno (21)
$$
To evaluate the above integral, we let $y=x+a$ and perform a Fourier 
transform in $a$.
For the Fourier transform in $a$ of  $r_{\L}^t (x,x+a)$, one gets,
$$
\s_{\L} (x,\xi) = \rho_{\L} (x) \, \rho_{\L} (\xi) \, . \leqno (22)
$$
For the Fourier transform in $a$ of $k(qx,x+a)$, one gets,
$$
\s (x,\xi) = \int f(\lb^{-1}) \left( \int \d (x+a-\lb qx) \, \a (a\xi) 
\, da \right) \, d^* \lb \, . \leqno (23)
$$
One has,
$$
\int \d (x+a-\lb qx) \, \a (a\xi) \, da = \a ((\lb q - 1) \, x \xi) \, , 
\leqno (24)
$$
thus (23) gives,
$$
\s (x,\xi) = \rho^{-1} \int_{A_S} g_q(u) \, \a (u x \xi) \, du \leqno 
(25)
$$
where,
$$
g_q(u) = f (q(u + 1)^{-1}) \, \vert u + 1 \vert^{-1} \, . \leqno 
(26)
$$
Since $f$ is smooth with compact support on $A_S^*$ the function $g_q$ 
belongs to $C_c^{\ify} (A_S)$.

\smallskip

Thus $\s (x,\xi) = \rho^{-1} \, \wh g_q (x \xi)$ and, using the Parseval 
formula we get,
$$
I_q= \int_{ x \in D , \,\vert x \vert \leq \L , \vert 
\xi \vert \leq \L} \s (x,\xi) \, dx \, d\xi \, . \leqno (27)
$$
 This gives,
$$
I_q= \rho^{-1} \int_{x \in D , \,\vert x \vert \leq \L 
, \vert \xi \vert \leq \L} \wh g_q (x\xi) \, dx \, d\xi \, . \leqno 
(28)
$$
With $u = x \xi$ one has $dx \, d\xi = du \, {dx \over \vert x 
\vert}$ and, for $\vert u \vert \leq \L^2$,
$$
\rho^{-1} \int_{x \in D , \,{\vert u 
\vert \over \L } \leq \vert x \vert  \leq \L  } \,
 {dx \over \vert x \vert} = 2 \log' \L - \log \vert u \vert \leqno (29)
$$

(using the precise definition of $\log' \L$ to handle the boundary 
terms). Thus we can rewrite (28) as,

$$
{\rm Trace} \, (R_{\L} \, T) = \sum_{q \in O_S^*} \, \int_{\vert u \vert 
\leq \L^2} \, \wh 
g_q (u) \, (2\log' \L - \log \vert u \vert) \, du \,\leqno (30) 
 $$
Now $\log \vert u \vert= \sum_{v \in S}\log \vert u_v \vert$, and we 
shall first prove that,
$$
 \sum_{q \in O_S^*} \, \int \, \wh 
g_q (u) \, du =\, h(1)  ,\leqno (31) 
$$
while for any $v \in S$,
$$
 \sum_{q \in O_S^*} \, \int \, \wh 
g_q (u) \,(- \log \vert u_v \vert) \, du =  \int'_{k^*_v} {h(u^{-1}) 
\over \vert 1-u \vert} \, d^* u \, .  \leqno (32) 
$$
In fact all the sums in $q$ will have only finitely many non zero terms.
It will then remain to control the error term, namely to show that,
$$
\sum_{q \in O_S^*} \, \int \, \wh 
g_q (u) \, ( \log \vert u \vert- 2 \log'\L)^{+} \, du= 0( \L^{-N}) 
\,\leqno (33) 
 $$
for any $N$, where we used the notation $ x^{+}=0$ if $x \leq 0$ and $ 
x^{+}=x$ if $x>0$.
Now recall that,
$$
g_q(u) = f (q(u + 1)^{-1}) \, \vert u + 1 \vert^{-1} \, , 
$$
so that $\int \, \wh 
g_q (u) \, du =\, g_q(0) =\,f(q)$. Since $f$ has compact
 support in $A^*_S$, the intersection of $O_S^*$
with the support of $f$ is finite and by (14) we get the equality (31).
\smallskip
To prove (32), we consider the natural projection 
$pr_v$ from $ \prod_{l \in S} k^*_l$ to $ \prod_{l \neq v} k^*_l$.
The image $pr_v( O_S^*)$ is still a discrete subgroup of 
$ \prod_{l \neq v} k^*_l$, (since $k_v^*$ is cocompact in $C_S$),
 thus there are only finitely many $q \in O_S^*$ such that $k_v^*$
 meets the support of $f_q$, where $f_q(a)= \, f(qa)$ for all $a$.
\smallskip
For each $q \in O_S^*$ one has, as in section V,
$$
 \int \, \wh 
g_q (u) \,(- \log \vert u_v \vert) \, du =  \int'_{k^*_v} {f_q(u^{-1}) 
\over \vert 1-u \vert} \, d^* u \, , \,\leqno (34)
$$
and from what we have just seen, this vanishes except for finitely many 
$q's$,
 so that by (14) we get the equality (32).
Let us prove (33). Let $ \ve_{\L}(u)=\, ( \log \vert u \vert- 2 
\log'\L)^{+}$, and,
$$
\delta_q(\L)= \, \int \, \wh 
g_q (u) \,\ve_{\L}(u)  \, du \,\leqno (35)
$$
be the error term. We shall prove,
\medskip
\noindent {\bf Lemma 2.} {\it  For any $\L$ the series 
$\sum_{O_S^*} \vert \delta_q(\L) \vert$ converges geometrically
 on the abelian finitely generated group $O_S^*$.
 Moreover its sum $ \s(\L)$ is $O(\L^{-N})$ for any $N$.}
\medskip
\noindent {\it Proof.} Let (cf. (8)),
 $d$ be a suitable word metric on $O_S^*$ such that,
$$
 Sup_S \vert q_v \vert \geq exp(d(q,1))  \qquad \fl q \in  O_S^* \leqno 
(36)
$$
\smallskip
Let $\xi \in \Sc(A_S)$ be defined by $\xi(x)= \,
 f(x^{-1})\vert x^{-1} \vert$ for all $x \in A^*_S$
 and extended by $0$ elsewhere. One has $g_q(x)= \, 
\xi(q^{-1}(1+x))$ for all $x \in A_S$, so that $\wh 
g_q (u) \,= \int g_q(x) \a (u \, x) \, dx = \, \a(-u)\, \wh \xi(q \, 
u)$.
Now, $ \delta_q(\L)= \, \int \, \wh 
g_q (u) \,\ve_{\L}(u)  \, du \,= \int \, \wh 
\xi(q \, u) \, \a(-u)\,\,\ve_{\L}(u)  \, du \,=  \int \, \wh 
\xi(y) \, \a(-q^{-1} \, y)\,\ve_{\L}(y)  \, dy $,
 since $\ve_{\L}(q \, u)\,=\ve_{\L}(u)$ for all $u$.
\smallskip
Thus we get, using the symbol $ \ov{F} \eta$ for the
 inverse Fourier transform of $ \eta$, the equality,
$$
\delta_q(\L)= \,\ov{F}(\ve_{\L}\wh \xi) (q^{-1}) . \leqno (37)
$$
Let $ \a \in ]0,1/2[$ and consider the norm,
$$
\Vert \eta \Vert = \, Sup_{x \in A_S} \vert F( \eta )(x) \, Sup_S \vert
x_v \vert^{ \a} \vert . \leqno (38)
$$
In order to estimate (38), we fix a smooth function $ \psi$ on $ \Rb$,
 equal to $1$ in a neighborhood of $0$ and with support in $[-1,1]$,
 and introduce the convolution operators,
$$
 (C_{ \a,v}* \eta)(x)= \, \int_{k_v} \, \psi( \vert \ve \vert)( \eta(x+ 
\ve)- \eta(x))
{d \ve \over \vert \ve \vert^{1+ \a}} \, \, ,\leqno (39)
$$
and the norms,
$$
\Vert \eta \Vert_{(1, \a, v)} = \, \Vert C_{ \a,v}* \eta \Vert_1 \, \, , 
\leqno (40)
$$
where $ \Vert \, \, \, \Vert_1$ is the $L^1$ norm.
\smallskip
The Fourier transform on $k_v$ of the distribution 
$ C_{ \a,v}$ behaves like $ \vert x_v \vert^{ \a}$ for $\vert x_v \vert 
\  \ra \ify$ . Thus, 
using the equality $F(C_{ \a,v}* \eta)=\, F(C_{ \a,v})\,F(\eta)$,
 and the control of the sup norm of $F(g)$ by the $L^1$ 
norm of $g$, we get an inequality of the form,
$$
Sup_{x \in A_S} \vert F( \eta )(x) \, Sup_S \vert
x_v \vert^{ \a} \vert \leq c_{ \a} \sum_S \Vert \eta \Vert_{(1, \a, v)} 
\, . \leqno (41)
$$
Let us now show that for any $\eta \in \Sc (A_S)$, and $ \a < 1/2 $, one 
has,
$$
 \Vert  \ve_{\L} \, \eta \Vert_{(1, \a, v)}= O( \L^{-N})  \, , \leqno 
(42)
$$
for any $N$.
\smallskip
One has $ \vert (\ve_{\L}(x+ \ve) \eta(x+ \ve)
 -\ve_{\L}(x) \eta(x))- \ve_{\L}(x)( \eta(x+ \ve)-\eta(x)) \vert \leq
\vert (\ve_{\L}(x+ \ve)- \ve_{\L}(x)) \vert \vert  \eta(x+ \ve) \vert$.
\smallskip
Moreover using the inequality,
$$
 \vert a^{+}- b^{+} \vert \leq \vert a - b \vert    \, , \leqno (43)
$$
we see that $\vert (\ve_{\L}(x+ \ve)- \ve_{\L}(x))
 \vert \leq \vert \log \vert x_v+ \ve \vert \,
 - \log \vert x_v \vert \vert $, for $ \ve \in k_v$. Let then,
$$
c'_{ \a }= \, \int_{k_v} \log \vert 1+y \vert {dy \over \vert y 
\vert^{1+ \a}} \, \, .\leqno (44)
$$
It is finite for all places $v \in S $
provided $ \a < 1/2 $, and one has,
$$
 \int_{k_v} \, \psi( \vert \ve \vert)
( \vert \log \vert x+ \ve \vert- \log \vert x \vert \vert)
{d\ve \over \vert \ve \vert^{1+ \a}} \,
 \leq c'_{ \a } \vert x \vert^{- \a} \, . \leqno (45)
$$
Thus one obtains the inequality,
$$
 \vert C_{ \a,v}* \ve_{\L}\, \eta  - \ve_{\L}\,(C_{ \a,v}* \eta) \vert 
(x) \leq 
c'_{ \a } \vert x_v \vert^{- \a} \, 
Sup_{ \ve \in k_v, \vert \ve \vert \leq 1} \, \vert \eta (x+ \ve) \vert 
. \leqno (46)
$$
Since the function $ \vert x_v \vert^{- \a}$ is locally integrable,
 for $\a < 1$, one has for $ \eta \in \Sc (A_S)$, and any $N$,
$$
\int_{X_{ \L}} \, \vert x_v \vert^{- \a} \,
 Sup_{  \ve \in k_v, \vert \ve \vert \leq 1} \, \vert \eta (x+ \ve) 
\vert dx \, =
O( \L^{-N}) \,, \leqno (47)
$$
where $X_{ \L} =\, \{\,y+ \ve;\, \vert y \vert \geq \L,\,
 \ve \in k_v, \,\vert \ve \vert \leq 1 \, \}$.
\smallskip
Moreover one has for any $N$,
$$
  \Vert \ve_{\L}\,(C_{ \a,v}* \eta) \Vert_1 = \, O( \L^{-N}). \leqno 
(48)
$$
Thus, using (46), we obtain the inequality (42).
\smallskip
Taking $ \eta = \wh \xi$ and using (41), we thus get
 numbers $ \d_{\L}$, such that $ \d_{\L}\, =
O( \L^{-N}) $ for all $N$ and that,
$$
 \vert \ov{F}(\ve_{\L}\wh \xi) Sup_S \vert
x_v \vert^{ \a} \vert \vert \leq  \d_{\L} 
 \qquad \fl x \in A_S \, \fl \L  . \leqno (49)
$$
Taking $x=q \in O_S^*$, and using (36) and (37), we thus get,
$$
 \vert \delta_q(\L) \vert \leq  
 \d_{\L} exp(-d(q,1))  \qquad \fl q \in  O_S^*   , \leqno (50)
$$
which is the desired inequality.

~\xx

\vglue 1cm
\noindent {\bf VIII The trace formula in the global case, and 
elimination of $\d$.}
\smallskip
 The main difficulty created by the parameter $\d$ in Theorem 1
 is that the formal trace computation of section VI is independent
 of $\d$, and thus cannot give in general the expected value of the
 trace of theorem  1, since in the latter each critical zero $ \rho$
 is counted with a multiplicity equal
to the largest integer  $ n < {1+\d \over 2}$ , $ n\leq $ multiplicity 
of
$ \rho$ as a zero of $L$. In particular for $L$ functions with multiple 
zeros,
 the $\d$-dependence of the spectral side is nontrivial. It is also 
clear that
 the function space $L_{\d}^2 (X)$ artificially eliminates the 
non-critical zeros
 by the introduction of the $\d$.
\smallskip
 As we shall see, all these problems are eliminated by the cutoff.
 The latter will be performed directly on the Hilbert space $L^2 (X)$
 so that the only value of $\d$ that we shall use is $\d =0$.
 All zeros will play a role in the spectral side of the trace formula,
 but while the critical zeros will appear per-se, the non critical ones
 will appear as resonances and 
enter in the trace formula through their harmonic potential with respect
 to the critical line. Thus the spectral side is entirely canonical and
 independent of $\d$, and by proving positivity of the Weil 
distribution,
 we shall show that its equality with the geometric side, i.e. the 
global
 analogue of Theorem 4, is equivalent to the Riemann Hypothesis for all
 $L$-functions with Gr\"ossencharakter.  
\smallskip
 The Abelian group $A$ of Adeles of $k$ is its own Pontrjagin
dual by means of the pairing
$$
\langle a,b \rangle = \a (ab) \leqno (1)
$$
where $\a : A \to U(1)$ is a nontrivial character which vanishes on $k
\sbs A$. Note that such a character is {\it not canonical\/}, but that 
any
two such characters $ \a $ and $\a^{'}$ are related by $k^*$,
$$ \a^{'}(a) = \a(qa) \quad \forall a \in A\,. \leqno (2)$$
It follows that the corresponding Fourier transformations on $A$ are
related by
$$\hat f^{^ {'}} = \hat f_q \,. \leqno (3)$$
This is yet another reason why it is natural to mod out by functions of
the form $f-f_q$, i.e.$\,$to consider the quotient space $X$.
\smallskip
We fix the additive character $\a$ as above, $\a= \prod \a_v$ and let
 $d$ be a differental idele,
$$
 \a (x) = \, \a_0 (d \, x)  \quad \forall x \in A\, , \leqno (4)
$$ 
where $\a_0= \prod \a_{0,v}$ is the product of the local
 normalized additive characters (cf [W1]).
We let $S_0$ be the finite set of places where $\a_v$ is ramified. 
\smallskip
 We shall first concentrate on the case of positive characteristic,
 i.e. of function fields, both because it is technically simpler
 and also because it allows to keep track of the geometric significance
 of the construction (cf. section II).
\smallskip
In order to understand how to perform in the global case,
 the cutoff $ R_{\L} = \wh{P}_{\L} \, P_{\L}$ of section VII,
 we shall first analyze the relative position of the pair of
 projections $\wh{P}_{\L}$, $P_{\L}$ when 
 $\L \ra \ify$. Thus, we let $S \supset S_0$ be a finite set
 of places of $k$, large enough so that $mod(C_S)= \, mod(C_k)
 =  \, q^{ \Zb}$ and that for any fundamental domain $D$ for
 the action of $O_S^*$ on $J_S$, the product $D \ts \prod R^*_v $
 is a fundamental domain for the action of $k^*$ on $J_k$.

Both $\wh{P}_{\L}$ and $ P_{\L}$ commute with the decomposition
 of $L^2 ( X_S )$ as the direct sum of the 
subspaces, indexed by characters $\chi_0$ of $C_{S,1}$, 
$$
L_{\chi_0}^2 = \{ \xi \in L^2 ( X_S) \, ; \, \xi (a^{-1} x)
 = \chi_0 (a) \, \xi (x),\, \,\, \fl  \, x \in  X_S \, , a \in C_{S,1} 
\} \leqno (5)
$$
which corresponds to the projections $ P_{\chi_0} = \int \ov{\chi_0} (a)
 \, U(a) \, d_1 \, a $, where $ d_1 \, a $ is the Haar 
measure of total mass $1$ on $C_{S,1}$.

\smallskip

\noindent {\bf Lemma 1.} {\it Let $\chi_0$ be a character 
of $C_{S,1}$, then for $\L$ large enough $\wh{P}_{\L}$ and
 $ P_{\L}$ commute on the Hilbert space $L_{\chi_0}^2 $. }
\medskip

\noindent {\it Proof.} Let $ \Uc_S$ be the image in $C_S$
 of the open subgroup $ \prod R^*_v $. It is a subgroup of finite
 index $l$ in $C_{S,1}$. Let us fix a character $ \chi$ of $ \Uc_S$ 
and consider the finite direct sum of the Hilbert spaces $L_{\chi_0}^2 $ 
where $ \chi_0$ varies among the characters of $C_{S,1}$
 whose restriction to $ \Uc_S$ is equal to $ \chi$,
$$
L^2 ( X_S)_{\chi} = \{ \xi \in L^2 ( X_S) \, ; \, \xi (a^{-1} x)
 = \chi (a) \, \xi (x) ,\, \,\,\fl  \, x \in  X_S \, , a \in \Uc_S \} 
\leqno (6)
$$
The corresponding orthogonal projection is $U(h_{\chi})$, 
where $h_{\chi} \in \Sc(C_S)$ is such that,
$$
Supp(h_{\chi})=\, \Uc_S  \qquad  h_{\chi}(x) = \, \lb 
\,  \ov{\chi}(x)  \qquad \fl  \, x \in  \Uc_S  \leqno (7)
$$
and the constant $ \lb = l / \log(q)$ corresponds
 to our standard normalization of the Haar measure on $C_S$.
Let as in section VII, $f \in \Sc(J_S)$ with support 
$ \prod R^*_v $ be such that $U(f)= \, U(h)$ and let
 $\xi \in \Sc(A_S)$ be defined by $\xi(x)= \, 
f(x^{-1})\vert x^{-1} \vert$ for all $x \in A^*_S$ and extended by $0$ 
elsewhere.
\smallskip
Since $\xi$ is locally constant, its Fourier transform has
 compact support and the equality (37) of section VII shows 
that for $\L$ large enough one has the equality,
$$
{\rm Trace} \, (\wh{P}_{\L} \, P_{\L} \, U(h_{\chi})) = 
2h_{\chi} (1) \log' \L + \sum_{v \in S} \int'_{k^*_v} {h_{\chi}(u^{-1})
\over \vert 1-u \vert} \, d^* u  \leqno (8)
$$
With $ \L= \, q^N$, one has $2  \log' \L = \, (2N+1) \log(q)$ so that,
$$
2h_{\chi} (1) \log' \L = \, (2N+1)l  \leqno (9)
$$
The character $ \chi$ of $ \prod R^*_v $ is a product,
 $ \chi= \, \prod \chi_v$ and if one uses the standard additive 
character
$ \a_0$ to take the principal value one has, (cf. [W1] Appendix IV),
$$ 
\int'_{R*_v} {\chi_v(u) \over \vert 1-u \vert} \, d^* u = \, -f_v 
\log(q_v)  \leqno (10)
$$
where $f_v$ is the order of ramification of $ \chi_v$. We thus get,
$$
\int'_{k^*_v} {h_{\chi}(u^{-1})
\over \vert 1-u \vert} \, d^* u  = \, -f_v \, deg(v) \, l + 
l \, {\log( \vert d_v \vert) \over \log(q) } \leqno (11)
$$
where $q_v = \, q^{deg(v)}$, and since we use the additive 
character $ \a_v$, we had to take into account the shift 
$\log( \vert d_v \vert) \, h_{\chi} (1)$ in the principal value.

\smallskip
Now one has $\vert d \vert = \prod \vert d_v \vert = q^{2-2g}$,
 where $g$ is the genus of the curve. Thus we get,
$$ 
{\rm Trace} \, (\wh{P}_{\L} \, P_{\L} \, U(h_{\chi}))= \,
 (2N+1)l \, - f \, l + (2-2g) \, l  \leqno (12)
$$
where $f = \sum_S f_v \, deg(v) $ is the order of
 ramification of $ \chi$, i.e. the degree of its conductor.
\smallskip
Let $B_{ \L} = \, Im(P_{\L}) \cap Im( \wh{P}_{\L}) $ be 
the intersection of the ranges of the projections $P_{\L}$ 
and $\wh{P}_{\L}$, and $B_{ \L}^{ \chi}$ be its intersection with
$L^2 ( X_S)_{\chi} $. We shall exhibit for each character $ \chi$ of $ 
\Uc_S$ a vector 
$ \eta_{ \chi} \in L^2 ( X_S)_{\chi} $ such that,
$$
U(g)(\eta_{ \chi}) \in \, B_{ \L} \qquad \fl g \in C_S , \vert g \vert 
\leq \L , 
\vert g^{-1} \vert \leq q^{2-2g-f} \, \L , \leqno (13)
$$
while the vectors $U(g)(\eta_{ \chi})$ are linearly independent 
for $g \in D_S$, where $D_S$ is the quotient of $C_S$ by the 
open subgroup $ \Uc_S$. 
\smallskip
With $ \L= \, q^N$ as above, the number of elements $g$ 
of $D_S$ such that $  \vert g \vert \leq \L , \vert g^{-1} 
\vert \leq q^{2-2g-f} \, \L $ is precisely equal to $(2N+1)l 
\, - f \, l + (2-2g) \, l$, which allows to conclude that the
 projections $\wh{P}_{\L}$ and$ P_{\L}$ commute in $L^2 ( X_S)_{\chi} $
 and that the subspace $B_{ \L}^{ \chi}$ is the linear span of the 
$U(g)(\eta_{ \chi}) $.
\smallskip
Let us now construct the vectors $ \eta_{ \chi} \in L^2 ( X_S)_{\chi} $. 
With the notations of
[W1] Proposition  VII.13, we let,
$$
 \eta_{ \chi} = \, \prod_S \phi_v   \leqno (14)
$$
be the standard function associated to $ \chi= \, \prod \chi_v$ 
so that for unramified $v$, $\phi_v$ is the characteristic function 
of $R_v$, while for ramified $v$ it vanishes outside
$R^*_v$ and agrees with $ \ov \chi_v$ on $R^*_v$. By construction 
the support of $ \eta_{ \chi}$ is contained in $R= \, \prod R_v$, thus 
one has 
$ U(g)(\eta_{ \chi}) \in \, Im(P_{\L})$ if $\vert g \vert \leq \L $.
 Similarly by [W1] Proposition  VII.13, we get that
 $ U(g)(\eta_{ \chi}) \in \, Im(\wh{P}_{\L})$ as soon
 as $\vert g^{-1} \vert \leq q^{2-2g-f} \, \L $. 
This shows that $\eta_{ \chi}$ satisfies (13) and it remains to show 
that
the vectors $U(g)(\eta_{ \chi})$ are linearly independent for $g \in 
D_S$. 
\smallskip
Let us start with a non trivial relation of the form,
$$
\Vert  \sum \lb_g U(g)(\eta_{ \chi})  \Vert  = \, 0 \leqno (15)
$$
where the norm is taken in $L^2 ( X_S)$, (cf. VII. 5). 
Let then $ \xi_{ \chi} = \, \prod_S \phi_v  \ot 1_R $ 
where $ R = \, \prod_{v \notin S} R_v$.
Let us assume first that $  \chi \neq 1$, then $ \xi_{ \chi}$ 
gives an element of $L_{\d}^2 (X)_0$
which is cyclic for the representation $U$ of $C_k$ 
in the direct sum of the subspaces $L_{\d,\chi_0}^2 (X)_0$ 
where $ \chi_0$ varies among the characters of $C_{k,1}$
 whose restriction to $ \Uc$ is equal to $ \chi$.
\smallskip
Now (15) implies that in $L_{\d}^2 (X)_0$
one has $\sum \lb_g U(g)(\xi_{ \chi})  = \, 0$.
 By the cyclicity of $ \xi_{ \chi}$ one then gets 
$\sum \lb_g U(g)  = \, 0$ on any $L_{\d,\chi_0}^2 (X)_0$
 which gives a contradiction (cf. Appendix 1, Lemma 3).
\smallskip
The proof for $  \chi = 1$ is similar but requires more care since $ 1_R 
\notin \Sc_0(A)$.~\xx

\smallskip
We can thus rewrite Theorem 4 in the case of positive characteristic as,

\medskip

\noindent {\bf Corollary 2.} {\it Let $Q_{\L}$ be the orthogonal
 projection on the subspace of $L^2 ( X_S)$ spanned by the
$f \in \Sc(A_S)$ which vanish as well as their Fourier transform
 for $ \vert x \vert > \L$. Let $h \in \Sc (C_S)$ have compact support. 
Then
  when 
 $\L \ra \ify$, one has
$$
{\rm Trace} \, (Q_{\L} \, U(h)) = 2h (1) \log' \L + \sum_{v \in S} 
\int'_{k^*_v} {h(u^{-1}) 
\over \vert 1-u \vert} \, d^* u + o(1)
$$
where $2 \log' \L = \int_{\lb \in C_S, \, \vert \lb \vert \in [\L^{-1}, 
\L]} d^* \lb$, 
and the other notations are as in Theorem VII.4.}

\medskip

In fact the proof of lemma 1 shows that the subspaces 
$B_{ \L}$ stabilize very quickly, so that the natural map 
$ \xi \rightarrow \xi \ot 1_R $ from $L^2 ( X_S)$ to $L^2 ( X_S')$ for 
$S \subset 
S'$ maps $B^S_{ \L}$ onto $B^{S'}_{ \L}$.
\smallskip
We thus get from corollary 2 an $S$-independent global formulation
 of the cutoff and of the trace formula. We let $L^2 ( X)$ be the
 Hilbert space $L_{\d}^2 ( X)$ of section III for the trivial value
 $ \d = \, 0$ which of course eliminates the unpleasant term from the
 inner product, and we let $Q_{\L}$ be the orthogonal projection on the 
subspace $B_{ \L}$ of $L^2 ( X)$ spanned by the
$f \in \Sc(A)$ which vanish as well as their Fourier transform for
 $ \vert x \vert > \L$. As we mentionned earlier, the 
proof of lemma 1 shows that for $S$ and $\L$ large enough 
(and fixed character $ \chi$), the natural map $ \xi \rightarrow 
\xi \ot 1_R $ from $L^2 ( X_S)_{\chi}$ 
to $L^2 ( X)_{\chi}$ maps $B^S_{ \L}$ onto $B_{ \L}$.
\smallskip
It is thus natural to expect that the following global
 analogue of the trace formula of corollary 2 actually holds, i.e. that 
when 
 $\L \ra \ify$, one has,
$$
{\rm Trace} \, (Q_{\L} \, U(h)) = 2h (1) \log' \L + \sum_{v} 
\int'_{k^*_v} {h(u^{-1}) 
\over \vert 1-u \vert} \, d^* u + o(1) \leqno (16)
$$
where $2 \log' \L = \int_{\lb \in C_k, \, \vert \lb \vert \in [\L^{-1}, 
\L]} d^* \lb$,
 and the other notations are as in Theorem VII.4.
\smallskip
We can prove directly that (16) holds when $h$ is supported by $C_{k,1}$ 
but 
are not able to prove (16) directly for arbitrary $h$
(even though the right hand side of the formula only 
contains finitely many nonzero terms since $ h \in \Sc (C_k)$ has 
compact support).
 What we shall show however is that the trace formula (16) 
implies the positivity of the Weil distribution, and hence 
the validity of RH for $k$. Remember that we are still in 
positive characteristic where RH is actually a theorem of A.Weil. 
It will thus be important to check the actual equivalence between
 the validity of RH and the formula (16). This is achieved by,

\medskip

\noindent {\bf Theorem 5.} {\it Let $k$ be a global field
 of positive characteritic and $Q_{\L}$ be the orthogonal
 projection on the subspace of $L^2(X)$ spanned by the
 $f \in \Sc(A)$ such that $f(x)$ and $\wh f(x)$ vanish
 for $\vert x \vert > \L$ . Let $h \in \Sc (C_k)$ have
 compact support. Then the following conditions are equivalent, 
\smallskip

a) When 
 $\L \ra \ify$, one has
$$
{\rm Trace} \, (Q_{\L} \, U(h)) = 2h (1) \log' \L + \sum_{v } 
\int'_{k^*_v} {h(u^{-1}) 
\over \vert 1-u \vert} \, d^* u + o(1)
$$  }
\smallskip
{\it b) All $L$ functions with Gr\"ossencharakter on $k$ satisfy the 
Riemann Hypothesis.}
\medskip

\noindent {\it Proof.} To prove that a) implies b), we shall prove 
(assuming a)) the positivity of the Weil distribution
 (cf. Appendix 2),
$$
\D = \log \vert d^{-1} \vert \, \d_1 + D - \sum_v D_v \, . \leqno (17)
$$
First, by theorem III.1 applied for $\d= \, 0$, the map $E$,
$$
E(f) \, (g) = \vert g \vert^{1/2} \sum_{q\in k^*} f(qg)
\qquad \fl \, g \in C_k \, , \leqno (18)
$$
defines a surjective isometry from $L^2(X)_0$ to $L^2(C_k)$ such that,
$$
E \, U (a) = \vert a \vert^{1/2} \, V(a) \, E \, , \leqno (19)
$$
where the left regular representation $V$ of $C_k$ on $L^2
(C_k)$ is given by,
$$
(V(a) \, \xi) \, (g) = \xi (a^{-1} \, g) \qquad \fl \, g,a
\in C_k \, . \leqno (20)
$$
Let $S_{ \L}$ be the subspace of $L^2
(C_k)$ given by,
$$
S_{\L} = \{ \xi \in L^2 (C_k) \, ; \ \xi (g) = 0 , \,\,\, \fl g \ , \ 
\vert g \vert \notin [ \L^{-1}, \,  \L] \} \, . \leqno (21)
$$
We shall denote by the same letter the corresponding orthogonal 
projection.
\smallskip
Let $B_{ \L,0}$ be the subspace of $L^2(X)_0$ spanned by the 
$f \in \Sc(A)_0$ such that $f(x)$ and $\wh f(x)$ vanish for
 $\vert x \vert > \L$ and $Q_{ \L,0}$ be the corresponding orthogonal
 projection. Let $f \in \Sc(A)_0$ be such that $f(x)$ and $\wh f(x)$ 
vanish for $\vert x \vert > \L$, then $E(f) \, (g)$ vanishes for
 $\vert g \vert > \L$, and the equality (Appendix 1),
$$
E(f) (g) = E (\wh f ) \left( {1 \over g} \right) \,\qquad  f \in \Sc 
(A)_0, \leqno (22)
$$
shows that $E(f) \, (g)$ vanishes for $\vert g \vert < \L^{-1}$.
\smallskip
This shows that $E(B_{ \L,0}) \subset S_{ \L}$, so that if
 we let $Q'_{ \L,0}= \, E \,Q_{ \L,0} \, E^{-1} $, we get the 
inequality,
$$
Q'_{ \L,0} \leq S_{ \L}   \leqno (23)
$$
and for any $ \L$ the following distribution on $C_k$ is of positive 
type,
$$
  \D_{ \L}(f) = \,  {\rm Trace} \, (( S_{ \L}- \, Q'_{ \L,0}) \, V(f))  
, \leqno (24)
$$
i.e. one has,
$$
\D_{ \L}(f * f^*) \geq \, 0 , \leqno (25)
$$
where $ f^*(g)= \, \ov f(g^{-1})$ for all $g \in C_k$.
\smallskip
 Let then $ f(g)= \,  \vert g \vert^{-1/2} \, h (g^{-1})$,
 so that by (19) one has $ E \, U(h)= \, V( \tilde f) E$ where 
$ \tilde f(g)= \, f(g^{-1})$ for all $g \in C_k$. By lemma 3 of Appendix 
2 one has,
$$
  \sum_v  D_v(f)-  \, \log \vert d^{-1} \vert = \, \sum_{v } 
\int'_{k^*_v} {h(u^{-1}) 
\over \vert 1-u \vert} \, d^* u . \leqno (26)
$$  
One has $ {\rm Trace} \, ( S_{ \L} \, V(f))= \, 2f (1) \log' \L$, 
thus using a) we see that the limit of  $\D_{ \L}$ when
 $\L \ra \ify$ is the Weil distribution $\D$ (cf.(17)).
 The term $D$ in the latter comes from the nuance between 
the subspaces $B_{ \L}$  and $B_{ \L,0}$. 
This shows using (24), that the distribution $\D$ is of positive type so 
that b) holds (cf. [W3]).
\smallskip
Let us now show that b) implies a).
 We shall compute from the zeros of $L$-functions and independently of 
any hypothesis
 the limit of the distributions $\D_{ \L}$ when $\L \ra \ify$.
\smallskip
\noindent We choose (non canonically) an isomorphism
$$
C_k \sm C_{k,1} \ts N \, . \leqno (27)
$$
where $N = {\rm
range} \ \vert \ \vert \sbs \Rb_+^* $, $ N \sm \Zb$ is the subgroup 
$q^{\Zb} \sbs
\Rb_+^*$ .
\smallskip
For $\rho \in \Cb$ we let $d\mu_{\rho}(z)
$ be the harmonic measure of $\rho$ with respect 
to the line $ i \, \Rb \subset \Cb$. It is a probability measure on 
the line $ i \, \Rb $ and coincides with the Dirac mass at $\rho$ when 
$\rho$ is on the line.
\smallskip
The implication b)$\Rightarrow$a) 
follows immediately from the explicit formulas (Appendix 2) and the 
following lemma,
\bigskip
\noindent {\bf Lemma 3.} {\it The limit of the distributions $\D_{ \L}$
 when $\L \ra \ify$ is given by,
$$
\D_{ \infty}(f) = \sum_{{L\left(\wt{\chi} , {1\over 2} + \rho
\right) =0 \atop \rho \in { B / N^{ \bot}}}} N(\wt{\chi},{1\over 2} + 
\rho) 
\,  \int_{z \in i \, \Rb} \wh f (\wt{\chi} ,z) d\mu_{\rho}(z)
$$
where $B$ is the open strip $B= \,  \{ \rho \in \Cb ; 
Re(\rho) \in ]{-1 \over 2}, {1 \over 2}[ \}$,
$N(\wt{\chi}, {1\over 2} + \rho)$ is the multiplicity 
of the zero, $d\mu_{\rho}(z)
$ is the harmonic measure of $\rho$ with respect to
 the line $ i \, \Rb \subset \Cb$, and the Fourier transform $\wh f$ of 
$f$ is defined by,
$$
\wh f (\wt{\chi} ,\rho) = \int_{C_k} f(u) \, \wt{\chi} (u) \, \vert u
\vert^\rho \, d^* \, u \, . $$ }
\medskip

\noindent {\it Proof.} Let $\L = \, q^N$. The proof of Lemma 1
 gives the lower bound $(2N+1) \, - f  + (2-2g)$ for the 
dimension of $B_{\L ,\chi}$ in terms of the order of ramification $f$ of 
the 
character $\chi$ of $C_{k,1}$, where we assume first that $\chi \not= 
1$.
 We have seen moreover
that $E(B_{\L ,\chi}) \subset S_{\L ,\chi}$ while the dimension
 of $S_{\L ,\chi}$ is $2N+1$.
\smallskip
Now by Lemma 3 of Appendix 1, every element $ \xi \in E(B_{\L ,\chi})$
 satisfies the  conditions,
$$
\int \xi(x) \, \chi (x) \, \vert x \vert^{\rho } \, d^* x = 0  \qquad  
\fl 
 \rho \in { B / N^{ \bot}} , \, L\left(\chi , {1\over 2} + \rho
\right) =0   . \leqno (28)
$$
This gives $2g-2+f$ linearly independent conditions 
(for $N$ large enough), using [W1] Theorem VII.6, and shows that they 
actually characterize the subspace $E(B_{\L ,\chi})$ of $S_{\L ,\chi}$.
\smallskip
This reduces the proof of the lemma to the following simple computation:
 One lets $F$ be a finite subset (possibly with 
multiplicity) of $\Cb^*$ and $E_N$ the subspace of 
$S_N= \, \{ \xi \in l^2( \Zb); \,\,\,  \xi(n)= \, 0 \,\, \fl n>N \}$
 defined by the conditions $ \sum \xi(n) \, z^n= \, 0 \,\, \fl z \in F$.
 One then has to compute the limit when 
$N \ra \ify$ of ${\rm Trace} \,((S_N \, -E_N) \, V(f))$ where $V$
 is the regular representation of $\Zb$ (so that
 $V(f)= \, \sum f_k \, V^k$ where $V$ is the shift, $V(\xi)_n= \, 
\xi_{n-1}$).
One then checks that the unit vectors
 $\eta_z \in S_N, \, z \in F, \,\,\, \eta_z(n)= \ov z^n \,
( \vert z^{2N+1} \vert -  \vert z^{-(2N+1)} \vert)^{-{1 \over 2}} 
( \vert z \vert -  \vert z^{-1} \vert)^{{1 \over 2}} \,\,  \, 
\fl n \in [-N, \, N] $, are asymptotically orthogonal 
and span $(S_N \, -E_N)$ (when $F$ has multiplicity one
 has to be more careful). The conclusion then follows from,
$$
{\rm Lim_{N \ra \ify }} \lgl V(f)\eta_z, \, \eta_z \rgl =
 \, \int_{ \vert u \vert=1} P_z(u) \, \wh f(u) \, du ,\leqno (29)
$$ 
where $P_z(u)$ is the Poisson kernel, and $\wh f$ the Fourier
 transform of $f$. ~\xx
\smallskip
One should compare this lemma with Corollary 2 of Theorem III.1.
 In the latter only the critical zeros were coming into play
 and with a multiplicity controlled by $\d$. In the above lemma,
 all zeros do appear and with their full multiplicity, but while
 the critical zeros appear per-se, the non-critical ones play 
the role of resonances as in the Fermi theory.
\medskip
Let us now explain how the above results extend to number fields $k$.
 We first need to analyze, as above, the relative position of 
the projections $P_{\L}$ and $\wh{P}_{\L}$.
Let us first remind the reader of the well known geometry 
of pairs of projectors. Recall that a pair of orthogonal 
projections $P_i$ in Hilbert space is the same thing as a
 unitary representation of the dihedral group
 $\Gamma = {\Zb / 2}*{\Zb / 2}$. To the generators 
$U_i$ of $\Gamma$ correspond the operators $2 P_i -1$. 
The group $ \Gamma$ is the semidirect product of the subgroup
 generated by $U= \, U_1 \, U_2$ by the group $\Zb / 2$, 
acting by $U \mapsto U^{-1}$. Its irreducible unitary 
representations are parametrized by an angle $ \theta \in [0, {\pi
\over 2}]$, the corresponding orthogonal projections $P_i$ being 
associated to the one dimensional subspaces $y= 0$ and $
y=  x \, tg(\theta)$ in the Euclidean $x,y$ plane. 
In particular these representations are at most two
 dimensional. A general unitary representation is
 characterized by the operator $ \Theta$ whose value
 is the above angle $ \theta$ in the irreducible case.
 It is uniquely defined by the equality,
$$
{\rm Sin( \Theta)} = \, \vert P_1 - \, P_2 \vert , \leqno (30)
$$
and commutes with $P_i$.
\smallskip

 The first obvious difficulty is that when $v$ is an 
Archimedian place there exists no non-zero function on $k_v$ which 
vanishes
 as well as its Fourier transform for $ \vert x \vert > \L$.
 This would be a difficult obstacle were it not for the work of 
Landau, Pollak and Slepian ([LPS]) in the early sixties, 
motivated by problems of electrical engineering, 
which allows to overcome it by showing that though
 the projections $P_{\L}$ and $\wh{P}_{\L}$ do not commute exactly 
even for large $ \L$, their angle is sufficiently well behaved so
 that the subspace $B_{ \L}$ makes good sense.
\smallskip
For simplicity we shall take $k=  \, \Qb$, so that the
 only infinite place is real. Let $P_{\L}$ be the
orthogonal projection onto the subspace,
$$
P_{\L} = \{ \xi \in L^2 (\Rb) \, ; \ \xi (x) = 0 ,\,\,\, \fl x \ , \ 
\vert x \vert > \L \} \, . \leqno (31)
$$
and $\wh{P}_{\L} = F P_{\L} F^{-1}$ where $F$ is the 
Fourier transform associated to the basic character $\a (x)= \,
 e^{-2\pi i x}$. What the above authors have done is
 to analyze the relative position of the projections $P_{\L}$, 
$\wh{P}_{\L}$ for $\L \ra \ify$ in order to account 
for the obvious existence of signals (a recorded music piece for 
instance) which 
for all practical purposes have finite support both 
in the time variable and the dual frequency variable. 
\smallskip
The key observation of ([LPS]) is that the following second
 order differential operator on $\Rb$ actually commutes with
the projections $P_{\L}$, 
$\wh{P}_{\L}$,
$$
H_{ \L}\psi(x) = \,- \part (( \L^2- \,x^2) \,  \part) \psi(x) \, + (2 
\pi \L x)^2 \, \psi(x) , \leqno (32)
$$
where $  \part$ is ordinary differentiation in one variable. 
Exactly as the generator $x \, \part$ of scaling commutes
with the orthogonal projection on the space of functions with positive 
support, 
the operator $\part (( \L^2- \,x^2) \,  \part)$ commutes with $P_{\L}$. 
Moreover $H_{ \L}$ commutes with Fourier transform $F$, and the 
commutativity of 
$H_{ \L}$ with $\wh{P}_{\L}$ thus follows.
\smallskip
If one sticks to functions with support in $[- \L, \L]$, the operator 
$H_{ \L}$ has discrete 
simple spectrum, and was studied long before the work of [LPS].
 It appears from the factorization
of the Helmoltz equation $\D \, \psi+ \,k^2 \, \psi= \,0$ in one 
of the few separable coordinate systems in Euclidean $3$-space,
 called the prolate spheroidal coordinates. Its eigenvalues
 $\chi_n( \L), n \geq 0$ are simple, positive and of the order
 of $n^2$ for $n \ra \ify$. The corresponding eigenfunctions 
$ \psi_n$ are called the prolate spheroidal wave functions and
 since $P_{\L}\,\wh{P}_{\L}
\,P_{\L}$ commutes with $H_{ \L}$, they are the
 eigenfunctions of $P_{\L}\,\wh{P}_{\L}
\,P_{\L}$. A lot is known about them, in particular 
one can take them to be real valued, and they are 
 even for $n$ even and odd for $n$ odd. The key result 
of [LPS] is that the corresponding eigenvalues $ \lb_n$
 of the operator $P_{\L}\,\wh{P}_{\L}
\,P_{\L}$ are decreasing very slowly from $\lb_0 \, \sm 1$
 until the value $n \sm 4 \L^2$ of the index $n$, they then
 decrease from $\sm 1$ to $\sm 0$ in an interval of length 
$\sm \log(\L)$ and then stay close to $0$. Of course this 
gives the eigenvalues of $\Theta$, it dictates the analogue
 of the subspace $B_{ \L}$ of lemma 1, as the linear span 
of the $  \psi_n, \, n \leq 4 \L^2 $, and it gives the 
justification of the semi-classical counting of the number
 of quantum mechanical states which are localized in the 
interval $[- \L, \L]$ as well as their Fourier transform 
as the area of the corresponding square in phase space.
\smallskip
We now know what is the subspace $B_{ \L}$ for the single
 place $ \ify$, and to obtain it for an arbitrary set of
 places (containing the infinite one), we just use the 
same rule as in the case of function fields, i.e. we consider the map,
$$
 \psi \mapsto \psi \ot 1_R ,  \leqno (33)
$$
which suffices when we deal with the Riemann zeta function.
 Note also that in that case we restrict ourselves to even 
functions on $ \Rb$. This gives the analogue of Lemma 1, Theorem 5, and 
Lemma 3. 
\smallskip

\noindent To end this section we shall come back to our original 
motivation 
of section I and show how the formula for the number of zeros
$$
N(E) \sim (E/2\pi) (\log (E/2\pi) -1) + 7/8 + o(1) + N_{osc} (E)  \leqno 
(34)
$$
appears from our spectral interpretation.

\smallskip

\noindent Let us first do a semiclassical computation for the number of 
quantum mechanical states in one degree of freedom which fulfill the 
following conditions,
$$
\vert q \vert \leq \Lambda , \vert p \vert \leq \Lambda , \vert H \vert 
\leq E \, , \leqno (35)
$$
where $H = qp$ is the Hamiltonian which generates the group of scaling 
transformations,
$$
(U(\lambda)\xi) (x) = \xi (\lambda^{-1} x) \qquad \lambda \in \Rb_+^* , 
x 
\in \Rb , \xi \in L^2 (\Rb) \, , \leqno (36)
$$
as in our general framework.

\smallskip

\noindent To comply with our analysis of section III, we have to 
restrict 
ourselves to even functions so that we exclude the region $pq \leq 0$ of 
the semiclassical $(p,q)$ plane.

\smallskip

\noindent Now the region given by the above condition is equal to $ D = 
D_+ \cup (-D_+) $ where,
$$
D_+ = \{ (p,q) \in \Rb_+ \times \Rb_+ , p \leq 
\Lambda , q \leq \Lambda , pq \leq E \} \, , \leqno (37)
$$
Let us compute the area of $D_+$ 
for the canonical symplectic form,
$$
\omega = {1 \over 2\pi} \, dp \wedge dq \, .  \leqno (38)
$$
By construction $D_+$ is the union of a rectangle with sides 
$E/\Lambda$, 
$\Lambda$ with the subgraph, from $q=E/\Lambda$ to $q=\Lambda$, of the 
hyperbola $pq = E$. Thus,
$$
\int_{D_+} \omega = {1 \over 2\pi} \, E/\Lambda \times \Lambda + {1 
\over 
2\pi} \, \int_{E/\Lambda}^{\Lambda} {E \, dq \over q} \, = {E \over 
2\pi} + 
{2E \over 2\pi} \, \log \Lambda - {E \over 2\pi} \, \log E \, .  \leqno 
(39)
$$
Now the above computation corresponds to the standard normalization of 
the 
Fourier transform with basic character of $\Rb$ given by
$$
\alpha (x) = \exp (i x) \, .  \leqno (40)
$$
But we need to comply with the natural normalization at the infinite 
place,
$$
\alpha_0 (x) = \exp (-2\pi ix) \, .  \leqno (41)
$$
We thus need to perform the transformation,
$$
P = p/2\pi \ , \ Q = q \, .  \leqno (41)
$$
The symplectic form is now $dP \wedge dQ$ and the domain,
$$
 D' = \{ (P,Q) ; \vert Q \vert \leq \Lambda , \vert P \vert \leq \Lambda 
, \vert 
PQ \vert \leq E/2\pi \} . \leqno (42)
$$
The computation is similar and yields the following result,
$$
\int_{D'_+} \omega = {2E \over 2\pi} \, \log \Lambda - {E \over 2\pi} 
\left( \log \, {E \over 2\pi} - 1 \right) \, .\leqno (43)
$$
In this formula we thus see the overall term $\langle N(E) \rangle$ 
which 
appears with a {\it minus} sign which shows that the number of quantum 
mechanical states corresponding to $D'$ is less than ${4E \over 2\pi} \, 
\log \Lambda$ by the first approximation to the number of zeros of zeta 
whose imaginary part is less than $E$ in absolute value 
(one just multiplies by $2$ the equality (43) since  $ D' = D'_+ \cup 
(-D'_+) $.
Now ${1 \over 
2\pi} \, (2E) (2\log \Lambda)$ is the number of quantum states in the 
Hilbert space $L^2 (\Rb_+^* , d^* x)$ which are localized in $\Rb_+^*$ 
between $\Lambda^{-1}$ and $\Lambda$ and are localized in the dual group 
$\Rb$ (for the pairing $\langle \lambda , t \rangle = \lambda^{it}$) 
between $-E$ and $E$. Thus we see clearly that the first approximation 
to 
$N(E)$ appears as the lack of surjectivity of the map which associates 
to 
quantum states $\xi$ belonging to $D'$ the function on $\Rb_+^*$,
$$
E(\xi) (x) = \vert x \vert^{1/2} \sum_{n\in \Zb} \xi (nx) \leqno (44)
$$
where we assume the additional conditions $\xi (0) = \int \xi (x) dx = 
0$.
\smallskip
A finer analysis, which is just what the trace formula is doing,
 would yield the additional terms $ 7/8 + o(1) + N_{osc} (E)$.
 The above discussion yields an explicit construction of a large
 matrix whose spectrum approaches the zeros of zeta as $ \L \ra \ify$.
\smallskip
It is quite remarkable that the eigenvalues of the angle operator $ 
\Theta$
which we discussed above, also play a key role in the theory of random 
unitary matrices.
To be more specific, let $E(n,s)$ be the large $N$ limit of the 
probability that there
are exactly $n$ eigenvalues of a random Hermitian $N \ts N$ matrix in 
the
interval $[-{\pi \over \sqrt{2N}}t,{\pi \over \sqrt{2N}}t]$, $t=s/2$.
Clearly $\sum \limits_n E(n,s)=1$.
Let $P_{t}$ be as above the operator of multiplication by $1_{[-t,t]}$ 
--
characteristic function of the interval $[-t,t]$ in the Hilbert space 
$L^2(\Rb)$.
In general (cf [Me]), $E(n,s)$ is $(-1)^n$ times the $n$-th coefficient 
of the Taylor
expansion at $z=1$ of $\zeta_s(z) =\prod 
\limits_{1}^{\infty}(1-z\lb_{j}(s))$,
 where
$\lb_{j}(s)$ are the eigenvalues of the operator
$\widehat{P_{\pi}}P_t$. (Here we denote by
$\widehat{P_{\lb}}={\cal F}P_{\lb}{\cal F}^{-1}$, and ${\cal F}$
denotes the Fourier transform,
${\cal F} \xi(u)=\int e^{ixu}\xi(x)dx$. Note also that the eigenvalues 
of $\widehat{P_a}P_b$ only depend upon 
the product $ab$ so that the relation with the eigenvalues of $\Theta$ 
should be clear.)

\medskip

\noindent {\bf  General remarks.}

a) There is a close analogy between the construction 
of the Hilbert space $L^2(X)$ in section III,
and the construction of the physical Hilbert space ([S] theorem 2.1)
in constructive quantum field theory, in the case of gauge theories.
 In both cases the action of the invariance group 
(the group $k^* =\, GL_1(k)$ in our case, the gauge group in the case of 
gauge theories) 
is wiped out by the very definition of the inner product. 
Compare the comment after III (9) with ([S]) top of page 17. 
\smallskip
b) For global fields of zero
 characteristic, the Idele class group has a non trivial 
connected component of the identity and this connected 
group has so far received no interpretation from Galois theory (cf 
[W4]).
The occurence of type III factors in [BC] indicates that 
the classification of hyperfinite type III  factors [C]
should be viewed as a refinement of local class field theory 
for Archimedean places, and provide the missing interpretation
 of the connected component of the identity in the Idele class group.
 In particular hyperfinite type III factors are classified by closed 
(virtual) subgroups 
of $\Rb_+^*$ (cf [C]) and they all appear as "unramified" extensions 
of the hyperfinite factor of type $III_1$.   
\smallskip
 c) Our construction of the Polya-Hilbert space bears some resemblance 
to [Z]
and in fact one should clarify their relation, as well as the relation
 of the space $X$ of Adele classes with the hypothetical arithmetic site 
of Deninger [D].
 Note that the division of $A$ by $k^*$ eliminates the linear structure 
of $A$ 
and that it transforms drastically the formulas for dimensions of 
function spaces,
 replacing products by sums (cf. theorem 4 of section VII). It should be 
clear to 
the reader that the action of the Idele class group on the space of 
Adele classes 
is the analogue (through the usual dictionnary of class field theory) of 
the action
 of the Frobenius on the curve. (To be more specific one needs to divide 
first the 
space $X$ of Adele classes by the action of the maximal compact subgroup 
of the Idele class group).
 \smallskip
 d) There is a superficial resemblance between the way the function 
$N(E)$ appears
 in the last computation of section VIII and the discussion in [BK],
 directly inspired from [Co]. It is amusing to note that the computation 
of [BK]
 is actually coincidental, the two rectangles are eliminated for no 
reason,
 which changes appropriately the sign of the term in $E$. What [BK] had 
not 
taken into account is that the spectral interpretation of [Co] is as an 
{ \it absorption spectrum } rather than an emission spectrum.
\smallskip 
 d) There is an even more superficial resemblance with the work of D. 
Goldfeld [G],
 in the latter the Weil distribution is used to define a corresponding 
inner product
 on a function space on the Idele class group. The positivity of the 
inner product 
is of course equivalent to the positivity of the Weil distribution (and 
by the 
result of A. Weil to RH) but this does not give any clue how to prove 
this 
positivity, nor does it give any explanation (except for a nice 
observation
 at the Archimedean place) for what the Weil distribution is, since it 
is 
introduced by hands in the formula for the inner product.
\smallskip
e) The above framework extends naturally from the case of $ GL(1)$ 
to the case of $GL(n)$ where the Adele class space is replaced by 
the quotient of $M_n(A)$ by the action on the right of $ GL_n(k)$.
 Preliminary work of C. Soule shows that the analogue of theorem III. 1
 remains valid, the next step is to work out the analogue of the 
Lefchetz
 trace formula in this context. 
\smallskip
f) I have been told by P. Sarnak and E. Bombieri that Paul Cohen has
 considered the space of Adele classes in connection with RH, but never
 got any detail of his unpublished ideas.
\smallskip
All the results of the present paper have been announced in the 
September 98
 conference on the Riemann hypothesis held in the Schrodinger Institute 
in Vienna and 
have been published as a preprint of the Schrodinger Institute. We are 
grateful to the American
Institute of Mathematics for its sponsoring of the conference.

\vglue 1cm

\medskip

\noindent {\bf Appendix I, Proof of theorem 1}

\smallskip

In this appendix we give the proof of theorem 1. Let us first recall as 
a preliminary the results of Tate an Iwasawa
as interpreted in [W 2]

\noindent {\bf $L$ functions and homogeneous distributions on $A$}

\smallskip
In general for a non archimedean local field $K$ we use the 
notations $R$ for the maximal compact subring, $P$ for the maximal ideal 
of $R$, $\pi $ for a generator of  
the ideal $P$ ( i.e.  $P = \pi R$).

Let $k$ be a global field and $A$ the ring of Adeles of $k$. It is the 
restricted product
of the local fields $k_v$ indexed by the set of places $v$ of $k$, with 
respect to the maximal compact subrings
$R_v$. Similarly, the Bruhat-Schwartz space  $\Sc (A)$ is the restricted 
tensor product of the local Bruhat-Schwartz spaces $\Sc (k_v)$, with 
respect
to the vectors $1_{R_v}$.
\smallskip
$L$ functions on $k$ are associated to Gr\"ossencharakters, i.e. to 
characters of the Idele class group,
$$
C_k = J_k / k^* \, . \leqno (1)
$$
Let $\Xc$ be a character of the idele class group, we consider $\Xc$ as 
a 
character of $J_k$ which is 1 on $k^*$. As such 
it can be written as a product,
$$
\Xc (j) = \Pi \, \Xc_v (j_v) \qquad j = (j_v) \in J_k \, . \leqno (2)
$$
By considering the restriction of $\Xc$ to the compact subgroup
$$
G_0 = \Pi \, R_v^* \ts \,1 \sbs J_k \, , \leqno (3)
$$
it follows that for all finite $v$ but a finite number, one has
$$
\Xc_v / R_v^* = 1 \, . \leqno (4)
$$
One says  that $\Xc$ is unramified 
at $v$ when this holds.
\smallskip

Then $\Xc_v (x)$ only depends upon the module $\vert x \vert$, since
$$
k_v^* / R_v^* = \mod (k_v) \, . \leqno (5)
$$
Thus $\Xc_v$ is determined by
$$
\Xc_v (\pi_v)  \leqno (6)
$$
which does not depend upon the choice of $\pi_v$ ($\mod R_v^*$). 
\smallskip
Let $\Xc$ be a quasi-character of $C_k$, it is of the form,
$$
\Xc (x) = \Xc_0 (x) \, \vert x \vert^s \leqno (7)
$$
where $s \in \Cb$ and $\Xc_0$ is a character of $C_k$. 
The real part $\s$ of $s$ is uniquely determined by
$$
\vert \Xc (x) \vert = \vert x \vert^{\s} \, . \leqno (8)
$$
  Let $P$ be the finite set of finite places where $\Xc_0$ is ramified.
The $L$ function  $L (\Xc_0 ,s)$ is defined for $\s =\,Re(s) > 1$ as
$$
 L (\Xc_0 ,s)=\left( {\displaystyle \prod_{v \, {\rm finite} \atop v 
\notin P}} (1 - \Xc_{0,v} (\pi_v) \, q_v^{-s})^{-1} \right)= \left( 
{\displaystyle \prod_{v \, {\rm finite} \atop v 
\notin P}} (1 - \Xc_v (\pi_v))^{-1} \right)\, \leqno (9)
$$
where 
$$
\vert \pi_v \vert=\, q_v^{-1}. \leqno (10)
$$
 Let us now recall from [W 2] how $L (\Xc_0 ,s)$ appears as a 
normalization factor for homogeneous distributions on $A$. 
\smallskip

First let $K$ be a local field and $\Xc$ a quasi-character of $K^*$,
$$
\Xc (x) = \Xc_0 (x) \, \vert x \vert^s \quad , \quad \Xc_0 : K^* \ra U 
(1) \, . 
\leqno (11)
$$
A distribution $D$ on $K$ is homogeneous of weight $\Xc$ iff one has
$$
\lgl f^a , D\rgl = \Xc(a)^{-1} \, \lgl f , D\rgl \leqno (12)
$$
for all test functions $f$ and all $a$ in $K^*$, where by definition
 $$
f^a (x) = f(ax)\leqno (13)
$$
When $\s= Re(s)>0$, there exists up to normalization only one 
homogeneous distribution of weight $\Xc$ on $K$,(cf [W 2]). 
It is given by the absolutely convergent integral,
$$
\int_{K^*} f(x)  \, \Xc(x) d^* x = \D_{\Xc}(f) 
\leqno 
(14)
$$

In particular, let  $K$ be non archimedean, then, for any compactly 
supported 
locally constant function $f$ on $K$ one has,
$$
f(x) - f(\pi^{-1} x)=\, 0 \quad \fl x,\, \vert x \vert \leq \d  \leqno 
(15)
$$
thus, for any $s \in \Cb $ the integral
$$
\int_{K^*} (f(x) - f(\pi^{-1} x)) \, \vert x \vert^s \, d^* x = \D'_s 
(f) 
\leqno 
(16)
$$
with the multiplicative Haar measure $d^* x$ normalized by
$$
\lgl 1_{R^*} , d^* x \rgl = 1 \, . \leqno (17)
$$

defines a distribution on $K$  with the properties,
$$
\lgl 1_R , \D'_s \rgl = 1 \leqno (18)
$$
$$
\lgl f^a , \D'_s \rgl = \vert a \vert^{-s} \, \lgl f , \D'_s \rgl \leqno 
(19)
$$
and
$$
\D'_s = (1 - q^{-s}) \, \D_s , \leqno (20)
$$
where $\vert \pi \vert =\, q^{-1} $.
(Let us check (18)$\ldots$(20). With $f = 1_R$ one has $f (\pi^{-1} x) = 
1$ iff 
$\pi^{-1} x \in R$ i.e. $x \in \pi R = P$. Thus $\D'_s (1_R) = 
\int_{R^*} \, 
d^* 
x = 1$. Let us check (20), one has $\int f (\pi^{-1} x) \, \vert x 
\vert^s \, 
d^* x = \int f (y) \, \vert \pi \vert^s \, \vert y \vert^s \, d^* y = 
\vert \pi 
\vert^s \, \D_s (f)$. But $\vert \pi \vert < 1$, $\vert \pi \vert = {1 
\over 
q}$. Note also that for $s=1$ and $f = 1_R$ one gets $\int_{R^*} dx = 
\left( 1 
- 
{1 \over q} \right) \int_R dx$.)

\smallskip

Let then $\Xc$ be a quasi-character of $C_k$ and write as above
$$
\Xc  = \Pi \, \Xc_v \, , \quad \,\quad \Xc (x) = \Xc_0 (x) \, \vert x 
\vert^s    \leqno (21)
$$
where $s \in \Cb$ and $\Xc_0$ is a character.
Let $P$ be the finite set of finite places where it is ramified.
 For any finite place $v \notin P$, let $\D'_v (s)$ be the unique 
homogeneous distribution of weight $\Xc_v$ normalized by
$$
\lgl \D'_v (s) , 1_{R_v} \rgl = 1 \, . \leqno (22)
$$
For any $v \in P$ or any infinite place, let, for $\s=Re(s) > 0$, 
$\D'_v$ be given by 
(14) which is homogeneous of weight $\Xc_v$ but unnormalized. Then the 
infinite tensor 
product,
$$
\D'_s = \Pi \, \D'_v (s) \leqno (23)
$$
makes sense as a continuous linear form on $ \Sc (A)$
 and is homogeneous of weight $\Xc$.
\smallskip

This solution is not equal to 0 since $\D'_v \not= 0$ for any $v \in P$ 
and any 
infinite place also. It is finite by construction of the space $\Sc (A)$ 
of 
test 
functions as an infinite tensor product
$$
\Sc (A) = \ot \, (\Sc (k_v) , 1_{R_v}) \, . \leqno (24)
$$

\smallskip

\noindent {\bf Lemma 1.} (cf [W 2]) {\it For $\s=Re(s) > 1$, the 
following integral converges 
absolutely
$$
\int f(x) \, \Xc_0 (x) \, \vert x \vert^s \, d^* x = \D_s (f) \qquad \fl 
\, f 
\in \Sc (A)
$$
and $\D_s (f) = L(\Xc_0 , s) \, \D'_s (f)$.}

\medskip

\noindent {\it Proof.} To get the absolute convergence one can assume 
that $f = 1_R$ and 
$\Xc_0 = 1$. Then one has to control an infinite product of local terms, 
given locally for the Haar measure $d^* x$ on $k_v^*$ such 
that $\int_{R_v^*} d^* x = 1$, by
$$
\int_{R \cap k_v^*} \vert x \vert^s \, d^* x \qquad (s \ \hbox{real}) 
\leqno (25)
$$
which is $1 + q_v^{-s} + q_v^{-2s} + \ldots = (1 - q_v^{-s})^{-1}$. Thus 
the 
convergence for $\s > 1$ is the same as for the zeta function.
\smallskip
 To prove the 
second equality one only needs to consider the infinite tensor product 
for finite places $v 
\notin P$. Then by (20) one has $\D'_v = (1 - q_v^{-\a_v}) \, 
\D_v$ where
$$
q_v^{-\a_v} = \Xc_v (\pi) = \Xc_{0,v} (\pi) \, q_v^{-s} \leqno (26)
$$
with $\vert \pi \vert = q_v^{-1}$.

\smallskip

Thus one gets $\D_s = \left( {\displaystyle \prod_{v \, {\rm finite} 
\atop v 
\notin P}} (1 - \Xc_{0,v} (\pi) \, q_v^{-s})^{-1} \right) \, \D'_s = L 
(\Xc_0 , 
s) \, \D'_s$.~\xx

\medskip

By construction $\D'_s$ makes sense whenever $\s > 0$ and is a 
holomorphic function of $s$ 
(for fixed $f$). Let us review biefly (cf [W2]) how to extend the 
definition of $\D_s$. 

\smallskip

We let as above $k$ be a global field, we fix a non trivial additive 
character $\a$ of $A$, trivial on $k$,
$$
\a (x+y) = \a (x) \, \a (y) \in U(1) \ , \ \a (q) = 1  \quad \fl \, q 
\in k \, 
. 
\leqno (27)
$$
We then identify the dual of the locally compact additive group $A$ with 
$A$ 
itself by the pairing,
$$
\lgl x , y \rgl = \a (xy) \, . \leqno (28)
$$
One shows (cf.[W 1]) that the lattice $k \sbs A$, i.e. the discrete and 
cocompact additive 
subgroup $k$, is its own dual,
$$
\lgl x,q \rgl = 1 \qquad \fl \, q \in k \qquad  \Lra \qquad x \in k \, . 
\leqno (29)
$$
Since $A$ is the restricted product of the local fields $k_v$ one can 
write 
$\a$ 
as an infinite product,
$$
\a = \Pi \, \a_v \leqno (30)
$$
where for almost all $v$ one has $\a_v = 1$ on $R_v$.
Let us recall the definition of the space $\Sc (A)_0 $, 
$$
\Sc (A)_0 = \{ f \in \Sc (A) \ ; \ f(0) = 0 \ , \ \int f \,
dx =0 \} \leqno (31)
$$

\smallskip

\noindent {\bf Lemma 2.} {\it Let $f \in \Sc (A)_0$, then the series
$$ E(f) \, (g) = \vert g \vert^{1/2} \sum_{q\in k^*} f(qg)
\qquad \fl \, g \in C_k \, $$
converges absolutely and one has
$$
\fl 
\, n \, , \ \exists \, c ,\qquad\vert E(f) (g) \vert \leq c \, e^{-n 
\vert \log \vert g \vert \vert} \qquad \fl \, g \in C_k \,
$$
and  $ E(\wh f) \, (g) = E(f) \, (g^{-1})$.}
\medskip
\noindent {\it Proof.} Let us first recall the formal definition ([Br]) 
of the Bruhat-Schwartz space $\Sc (G)$ for an arbitrary locally compact
abelian group $G$. One considers all pairs of subgroups $G_1,G_2$ of $G$ 
such that $G_1$ is generated by a compact
neighborhood of $0$ in $G$, while $G_2$ is a compact subgroup of $G_1$ 
such that the quotient group is elementary, i.e. is of the
form $ \Rb^a \, \Tb^b \, \Zb^c \, F$  for $F$ a finite group. By 
definition the Bruhat-Schwartz space $\Sc (G)$ is the inductive limit of 
the Schwartz spaces $\Sc (G_1/G_2)$
where the latter have the usual definition in terms of rapid decay of 
all derivatives. Since $G_1$ is open in $G$, any element of $\Sc 
(G_1/G_2)$
extended by $0$ outside $G_1$ defines a continuous function on $G$. By 
construction  $\Sc (G)$ is the union of the subspaces $\Sc (G_1/G_2)$
and it is endowed with the inductive limit topology.  
\smallskip
Let $\wh G$ be the Pontrjagin dual of $G$, then the Fourier transform, 
which depends upon the normalization of the Haar measure on $G$,
gives an isomorphism of $\Sc (G)$ with $\Sc (\wh G)$.
\smallskip
          Let $\G$ be a 
lattice in the locally compact abelian group $G$. Then any function $f 
\in \Sc (G)$ is admissible for the
pair $G, \G$ in the sense of [W 1], and the Poisson summation formula 
(cf [W 1]) is the equality,
$$
\Covol \, (\G) \sum_{\g \in \G} f(\g) = \sum_{\b \in \G^{\perp}} \wh f 
(\b) 
\leqno (32)
$$
where $\G^{\perp}$ is the dual of the lattice $\G$, and 
$$
\wh f (\b) = \int f(a) \, \b (a) \, da \, . \leqno (33)
$$
Both sides of (32) depend upon the normalization of the Haar measure on 
$G$.
\smallskip 
In our case we let $A$ be as above the additive group of Adeles on $k$.
We normalize the additive Haar measure $dx$ on $A$ by
$$
\Covol (k) = 1 \, . \leqno (34)
$$
We then take $\G = xk$, for some $x \in A^{-1}$. One has
$$
\Covol \, (xk) = \vert x \vert \leqno (35)
$$
The dual $\G^{\perp}$ of the lattice $xk$, for $x$ invertible in $A$, is 
the 
lattice $\G^{\perp} = x^{-1} k$. Thus the Poisson formula (32) reads, 
for any $f \in \Sc (A)$,
$$
\vert x \vert \sum_{q \in k} f(xq) = \sum_{q \in k} \wh f (x^{-1} q) \, 
. 
\leqno 
(36)
$$
Which we can rewrite as,
$$
\vert x \vert \sum_{k^*} f(xq) = \sum_{k^*} \wh f (x^{-1} q) + \d \leqno 
(37)
$$
$\d = - \vert x \vert \, f(0) + \int f(y) \, dy$.

\smallskip 
We can then rewrite (37) as the equality, valid for all $f \in \Sc 
(A)_0$
$$
E(f) (x) = E (\wh f ) \left( {1 \over x} \right) \,\qquad  f \in \Sc 
(A)_0. \leqno (38)
$$
It remains to control the growth of $E(f) (x)$ on $C_k$, but by (38), it 
is enough to understand what happens for $ \vert x \vert $ large.
\smallskip
We only treat the case of number fields, the general case is similar. 
Let $A = A_f \ts A_{\ify}$ be the decomposition of the ring of Adeles 
corresponding to finite and
infinite places, thus ${\displaystyle 
A_{\ify}=\prod_{S_{\ify}}} k_v$ where $S_{\ify}$ is the set of infinite 
places.
\smallskip
Any element of $\Sc (A)$ is a finite linear combination of test 
functions of 
the form,
$$
f = f_0 \ot f_1 \leqno (39)
$$
where $f_0 \in \Sc (A_f) \ , \ f_1 \in \Sc (A_{\ify})$ (cf [W 5] 39), 
thus it is enough to control the growth of $E(f) (x)$ for such $f$
and $ \vert x \vert $ large. \smallskip
Let $J_{k,1} = \{ x \in J_k ; \vert x \vert = 1 \} $ be the group of 
Ideles of module one, since $J_{k,1} / k^*$ is compact (cf [W 1]),
we shall fix a compact subset $K_1$ of $J_{k,1}$ whose image in $J_{k,1} 
/ k^*$ is this compact group.
\smallskip
Let $\mu$ be the diagonal embedding:
$$
\lb \in \Rb_+^* \build \longra_{}^{\mu} (\lb , \ldots , \lb) \in 
\prod_{S_{\ify}
} k_v^* \leqno (40)
$$
which yields an isomorphism
$$
J_k = J_{k,1} \ts \Im \, \mu \, . \leqno (41)
$$
One has $f_0 \in \Sc (A_f)$, hence (cf [W 5]), $f_0 \in C_c (A_f)$ and 
we let $K_0 = \hbox{Support} \ f_0$.
Since $K_0$ is compact, one can find a finite subset $P$ of the set 
of 
finite places and $C < \ify$ such that:
$$
y \in K=(K_f)^{-1}K_0 \Ra \vert y_v \vert \leq 1 \quad \fl \, v \notin P 
\qquad , \qquad 
\vert 
y_v \vert \leq C \quad \fl \, v \, . \leqno (42)
$$
where $K_f$ is the projection of $K_1$ on $A_f$.\smallskip
 We let $\Om$ be the compact open subgroup of $ A_f$ determined by
$$
\vert a_v \vert \leq 1 \quad \fl \, v \notin P \qquad , \qquad \vert a_v 
\vert 
\leq C \quad \fl \, v \, . \leqno (43)
$$

\smallskip
By construction $E(f) (x)$ only depends upon the class of $x$ in $J_k / 
k^*$. Thus, to control the behaviour of $E(f)(x)$ for $\vert x \vert 
\ra 
\ify$, we can take $x=(x_f,x_{\ify}) \in K_1$ and consider $E(f)(\lb 
\,x)$ for $\lb \in \Rb_+^* \,$,$ \, \lb \ra 
\ify$.
Now let $q = (q_f , q_{\ify}) \in k$, then,
$$
f(q\,\lb \,x) = f_0 (q_f \, x_f) \, f_1 (q_{\ify} \,\lb\, x_{\ify}) 
\leqno (44)
$$
and this vanishes unless $q_f \, x_f \in K_0$, i.e. unless $q_f \in K$. 
But 
then by (42) one has $q_f \in \Om$. Let $\G$ be the lattice in 
${\displaystyle \prod_{S_{\ify}}} 
k_v$ determined by
$$
\G= \{q_{\ify}; \, q \in k, \, q_f  \in \Om, \} \leqno (45)
$$
The size of $E(f) (\lb \,x)$ is thus controlled ( up to the square root 
of $\vert \lb \,x \vert $ ) by
$$
C \sum_{n \in \G^*} \vert f_1 (\lb \, x_{\ify} \, n) \vert \leqno (46)
$$
where 
$x_{\ify}$ varies in the projection $K_ {\ify}$ of $K_1$ on 
${\displaystyle \prod_{S_{\ify}}} 
k_v^*$.

\smallskip
Since $f_1 \in \Sc (A_{\ify})$, this shows that  
$E(f)(x)$  decays faster than any power of $\vert x \vert$ for 
$\vert x \vert \ra \ify$.

\smallskip
We have shown that $E(f)$ has rapid decay in terms of $\vert x \vert$, 
for 
$\vert x \vert \ra \ify$. 
Using (38) and the stability of $\Sc (A)_0$ under Fourier, we see that 
it also 
has exponential decay in terms of $\vert \log \vert x \vert \vert$ when 
$\vert 
\log \vert x \vert \vert \ra \ify$.
\smallskip We then get,
\medskip
\smallskip
\noindent {\bf Lemma 3.} (cf [W 2]) {\it For $\s=Re(s) > 0$, and any 
character $\Xc_0$ of $C_k$, one has 
$$
\int E(f)(x) \, \Xc_0 (x) \, \vert x \vert^{s-1/2} \, d^* x = c L(\Xc_0 
, s) \, \D'_s (f) \qquad \fl \, f 
\in \Sc (A)_0
$$
where the  non zero constant $c$ depends upon the normalization of the 
Haar measure $d^* x$ on $C_k$.}
\medskip
\noindent {\it Proof.} For $\s=Re(s) > 1$, the equality follows from 
lemma 1, but since both sides are analytic in $s$ it holds in general. 
\smallskip
As in lemma 1, we shall continue to use the notation $\D_s (f)$ for 
$\s=Re(s) > 0$.  \smallskip

\vglue 1cm

\noindent {\bf  Approximate units in the Sobolev spaces $L_{\d}^2 
(C_k)$}

\smallskip

We first consider, for $\d > 1$, the Hilbert space $L_{\d}^2 (\Rb)$
of functions $\xi (u)$, $u 
\in \Rb$ with square norm given by
$$
\int_{\Rb} \vert \xi (u) \vert^2 \, (1 + u^2)^{\d / 2} \, du \, . \leqno 
(1)
$$
We let $\rho (u) = (1 + u^2)^{\d / 2}$. It is comparable to $(1 + \vert 
u 
\vert)^{\d}$ and in particular,
$$
{\rho (u+a) \over \rho (u)} \leq c \, \rho (a) \qquad \fl \, u \in \Rb 
\, , \ a 
\in \Rb \leqno (2)
$$
with $c=2^{\d/2}$.
\smallskip

We then let $V(v)$ be the translation operator,
$$
(V(v) \xi) (u) = \xi (u-v) \qquad \fl \, u , v \in \Rb \, . \leqno (3)
$$
One has $\int_{\Rb} \vert \xi (u-v)\vert^2 \, \rho (u) \, du = 
\int_{\Rb} \vert 
\xi (u) \vert^2 \, \rho (u+v) \, du$ so that by (2) it is less than $c 
\int_{\Rb}   \vert \xi (u)\vert^2  \, \rho (u) \,\rho (v) \, du = c \, 
\rho (v) 
\, 
\Vert \xi \Vert^2$,
$$
\Vert V(v) \Vert \leq (c \, \rho (v))^{1/2} \, . \leqno (4)
$$
This shows that $V(f) = \int f(v) \, V(v) \, dv$ makes sense as soon as
$$
\int \, \vert f(v)\vert \, \rho (v)^{1/2} \, dv < \ify \, . \leqno (5)
$$
This holds for all $f \in \Sc (\Rb)$.

\medskip

\noindent {\bf Lemma 4.} {\it There exists an approximate unit $f_n \in 
\Sc (\Rb)$, such that 
$\wh{f}_n$ has compact support, $\Vert V (f_n) \Vert \leq C \quad \fl \, 
n$, 
and
$$
V (f_n) \ra 1 \ \hbox{strongly in} \ L_{\d}^2 (\Rb) \, .
$$
}

\smallskip

\noindent {\it Proof.} Let $f$ be a function, $f \in \Sc (\Rb)$, whose 
Fourier 
transform $\wh f$ has compact support, and such that $\int f \, dx = 1$ 
(i.e. 
$\wh f (0) = 1$). Let then
$$
f_n (v) = n \, f (nv) \qquad n = 1,2,\ldots \leqno (6)
$$
One has $\int \vert f_n(v)\vert \ \, \rho (v)^{1/2} \, dv = \int \vert 
f(u)\vert \ \, \rho \left( {u 
\over n} \right)^{1/2} \, du \leq \int \vert f(u)\vert \ \, \rho 
(u)^{1/2} \, du$. Thus 
$\Vert V (f_n) \Vert$ is uniformly bounded.

\smallskip

We can assume that $\wh f$ is equal to 1 on $[-1,1]$, then $\wh{f}_n$ is 
equal 
to 1 on $[-n,n]$ and $V (f_n) \xi = \xi$ for any $\xi$ with $\Supp 
\wh{\xi} 
\sbs 
[-n,n]$ . By uniformity one gets that $V(f_n) \ra 1$ strongly.~\xx

\medskip

Let us now identify the dual $(L_{\d}^2)^*$ of the Hilbert space 
$L_{\d}^2$ 
with 
$L_{-\d}^2$ by means of the pairing,
$$
\lgl \xi , \eta \rgl_0 = \int_{\Rb} \xi (u) \, \eta (u) \, du \, . 
\leqno (7)
$$
Since $L_{\d}^2$ is a Hilbert space, it is its own dual using the 
pairing,
$$
\lgl \xi , \eta_1 \rgl = \int_{\Rb} \xi (u) \, \eta_1 (u) (1+u^2)^{\d / 
2} \, 
du 
\, . \leqno (8)
$$
If we let $\eta (u) = \eta_1 (u) (1+u^2)^{\d / 2}$, then
$$
\int \vert \eta_1 (u) \vert^2 \, (1+u^2)^{\d / 2} \, du = \int \vert 
\eta (u) 
\vert^2 \, (1+u^2)^{-\d / 2} \, du
$$
which is the natural norm square for $L_{-\d}^2$.

\smallskip

Given a quasicompact group such as $C_k$ with module,
$$
\vert \ \vert : C_k \ra \Rb_+^* \leqno (9)
$$
we let $d^* g$ be the Haar measure on $C_k$ normalized by
$$
\int_{\vert g \vert \in [1,\L]} d^* g \sim \log \L \qquad \L \ra \ify 
\leqno 
(10)
$$
and we let $L_{\d}^2 (C_k)$ be defined by the norm,
$$
\int_{C_k} \vert \xi (g) \vert^2 \, (1 + \log \vert g \vert^2)^{\d / 2} 
\, d^* 
g 
\, . \leqno (11)
$$
It is, when the module of $k$ is $\Rb_+^*$, a direct sum of spaces (1), 
labelled by 
the 
characters $\Xc_0$ of the compact group
$$
C_{k,1} = \Ker \mod \, . \leqno (12)
$$
The pairing between $L_{\d}^2 (C_k)$ and $L_{-\d}^2 (C_k)$ is given by
$$
\lgl \xi , \eta \rgl = \int \xi (g) \, \eta (g) \, d^* g \, . \leqno 
(13)
$$
The natural representation $V$ of $C_k$ by translations is given by
$$
(V(a) \xi) (g) = \xi (a^{-1} g) \qquad \fl \, g,a \in C_k \, . \leqno 
(14)
$$
It is not unitary but by (4) one has,
$$
\Vert V(g) \Vert = 0 \, \vert \log \vert g \vert \vert^{\d / 2} \, , \ 
\vert 
\log \vert g \vert \vert \ra \ify \, . \leqno (15)
$$
Finally, one has, using lemma 4 and the decomposition $ C_k = C_{k,1} 
\ts N$,
\medskip

\noindent {\bf Lemma 5.} {\it There exists an approximate unit $f_n \in 
\Sc (C_k)$, such that 
$\wh{f}_n$ has compact support, $\Vert V (f_n) \Vert \leq C \quad \fl \, 
n$, 
and
$$
V (f_n) \ra 1 \ \hbox{strongly in} \ L_{\d}^2 (C_k) \, .
$$
}

\vglue 1cm

\noindent {\bf Proof of theorem III 1}
\smallskip

We first consider the subspace of codimension 2 of $\Sc (A)$ given by
$$
f(0) = 0 \, , \ \int f \, dx = 0 \, . \leqno (1)
$$
On this subspace $\Sc (A)_0$ we put the inner product,
$$
\int_{C_k} \vert E (f) (x) \vert^2 \, (1 + \log \vert x \vert^2)^{\d / 
2} \, 
d^* 
x \, . \leqno (2)
$$
We let $U$ be the representation of $C_k$ on $\Sc (A)$ given by
$$
(U (a) \xi) (x) = \xi (a^{-1} x) \qquad \fl \, a \in C_k \, , \ x \in A 
\, . 
\leqno (3)
$$

\smallskip  We let $L_{\d}^2 (X)_0$ be the separated completion of $\Sc 
(A)_0$ for the 
inner product given by (2). The linear map $E : \Sc (A)_0 \ra L_{\d}^2 ( 
C_k )$ 
satisfies
$$
\Vert E(f) \Vert_{\d}^2 = \Vert f \Vert_{\d}^2 \leqno (4)
$$
by construction. Thus it extends to an isometry, still noted $E$,
$$
E : L_{\d}^2 (X)_0 \hra L_{\d}^2 ( C_k ) \, . \leqno (5)
$$
One has 
$$
\displaylines{
E (U (a) f) (g) = \vert g \vert^{1/2} \sum_{k^*} (U(a) f) (qg) = \vert g 
\vert^{1/2} \sum_{k^*} f(a^{-1} qg)  \cr
=\vert g \vert^{1/2} \sum_{k^*} f(q a^{-1} g) = \vert a \vert^{1/2} \, 
\vert 
a^{-1} g \vert^{1/2} \sum_{k^*} f (q a^{-1} g) = \vert a \vert^{1/2} 
(V(a) \, 
E(f)) (g) \cr
}
$$
$$
E \, U (a) = \vert a \vert^{1/2} \, V(a) \, E \, . \leqno (6)
$$

The equality (6) shows that the natural representation $U$ of $C_k$ on 
$L_{\d}^2 (X)_0$ corresponds by the isometry $E$ to the restriction of 
$\vert a 
\vert^{1/2} \, V(a)$ to the invariant subspace given by the range of 
$E$.

\smallskip

In order to understand $\Im E $ we consider its orthogonal in the dual 
space $L_{-\d}^2 ( C_k )$.\smallskip

 The compact subgroup
$$
C_{k,1} = \{ g \in C_k \, ; \, \vert g \vert = 1 \} \leqno (7)
$$
acts by the representation $V$ which is unitary when restricted to 
$C_{k,1}$. 
Thus one can decompose $L_{\d}^2 ( C_k )$ and its dual $L_{-\d}^2 ( C_k 
)$ , in the direct sum of the 
subspaces,
$$
L_{\d , \Xc_0}^2 = \{ \xi \in L_{\d}^2 ( C_k ) \, ; \, \xi (a^{-1} g) = 
\Xc_0 
(a) \, \xi (g) \qquad \fl \, g \in C_k \, , \ a \in C_{k,1}\} \leqno (8)
$$
and,
$$
L_{-\d , \Xc_0}^2 = \{ \xi \in L_{-\d}^2 ( C_k ) \, ; \, \xi (a \, g) = 
\Xc_0 
(a) \, \xi (g) \qquad \fl \, g \in C_k \, , \ a \in C_{k,1}\} \leqno (9)
$$

which corresponds to the projections $P_{\Xc_0} = \int \ov{\Xc_0} (a) \, 
V(a) \, d_1 \, a$ for $L_{\d}^2$ 
and $P_{\Xc_0}^t = \int \ov{\Xc_0} (a) \, V(a)^t \, d_1 \, a$ for the 
dual space 
$L_{-\d}^2$.

\smallskip
In (9) we used the formula
$$
(V (g)^t \, \eta) (x) = \eta (gx) \leqno (10)
$$
which follows from the definition of the transpose, $\lgl V(g) \xi , 
\eta \rgl = \lgl \xi , V(g)^t \eta \rgl$ 
using
$$
\int \xi (g^{-1} x) \, \eta (x) \, d^* x = \int \xi (y) \, \eta (gy) \, 
d^* y
$$

In these formulas one only uses the character $\Xc_0$ as a character of 
the 
compact subgroup $C_{k,1}$ of $C_k$. One now chooses, non canonically, 
an 
extension $\wt{\Xc}_0$ of ${\Xc}_0$ as a character of $C_k$ 
$$
\wt{\Xc}_0 (g) = {\Xc}_0 (g) \qquad \fl g \in C_{k,1} \, . \leqno (11)
$$
This choice is not unique but any two such extensions differ by a 
character 
which is principal, i.e. of the form: $g \ra \vert g \vert^{is_0}$, $s_0 
\in \Rb$.
\smallskip
Let us fix a factorization $ C_k= C_{k,1} \ts  \Rb_+^*   $, and fix
$\wt{\Xc}_0$ as being equal to $1$ on $ \Rb_+^*   $.
\smallskip
 We then write any element of $L_{-\d , \Xc_0}^2 (C_k)$ 
in the form
$$
g \in C_k \ra \eta (g)=\wt{\Xc}_0 (g) \, \psi (\vert g \vert) \leqno 
(12)
$$
where 
$$
\int \vert \psi (\vert g \vert)^2 \, (1 + (\log \vert g \vert)^2)^{-\d / 
2} \, 
d^* g < \ify \,\leqno (13) .
$$
This vector is in the orthogonal of $\Im E$ iff
$$
\int E(f) (x) \, \wt{\Xc}_0 (x) \, \psi (\vert x \vert) \, d^* x = 0 
\qquad 
\fl \, f \in \Sc (A)_0 \, . \leqno (14)
$$
We first proceed formally and write $\psi (\vert x \vert) = \int 
\wh{\psi} (t) 
\, \vert x \vert^{it} \, dt$ so that the left hand side of (14) becomes,
$$
\int \int E(f) (x) \, \wt{\Xc}_0 (x) \, \vert x \vert^{it} \, 
\wh{\psi} (t) \, d^* x \, dt = \int \D_{1/2 + it} (f) \, \wh{\psi} (t) 
\, dt 
\leqno (15)
$$
(using the notations of lemmas 1 and 3).

\smallskip

Let us justify this formal manipulation; since we deal with the 
orthogonal of 
an invariant subspace, we can assume that 
$$
V^t (h) \, \eta = \eta , \leqno (16)
$$ 
for some $h$ such that $\wh h $ has compact 
support. Indeed we can use lemma 5 to only consider vectors 
which belong to the range of
$$
V^t (h) = \int h(g) \, V(g)^t \, d^* g \, , \ \wh h \ \hbox{with compact 
support}. $$

\smallskip

Then, using (16), the Fourier transform 
of the tempered distribution $\psi$ on $\Rb_+^*$ has compact support in 
$\Rb$.
 Thus, since $E(f) (x)$ has rapid decay,
the equality between (14) and (15) follows from the definition of the 
Fourier transform 
of the tempered distribution $\psi$ on $\Rb_+^*$.

\smallskip Let us now describe suitable test functions $f \in \Sc (A)_0$ 
in order to test 
the distribution, 
$$
\int \D_{{1 \over 2} + it} \, \wh{\psi} (t) \, dt \leqno (17)
$$

\smallskip

We treat the case of characteristic zero, the general case is similar. 
For the finite places we take,
$$
f_0 = \build \ot_{v \notin P}^{} 1_{R_v} \ot f_{\Xc_0} \leqno (18)
$$
where $f_{\Xc_0}$ is the tensor product over ramified places of the 
functions 
equal to 0 outside $R_v^*$ and to $\ov{\Xc}_{0,v}$ on $R_v^*$. It 
follows then 
by the definition of $\D'_s $  that,
$$
\lgl \D'_s , f_0 \ot f \rgl = \int f(x) \, \Xc_{0,\ify} (x) \, \vert x 
\vert^s 
\, d^* x \leqno (19)
$$
for any $f \in \Sc (A_{\ify})$.

\smallskip

Moreover if the set $P$ of finite ramified places is not empty one has,
$$
f_0 (0) = 0 \, , \ \int_{A_f} f_0 (x) \, dx = 0 \leqno (20)
$$
so that $f_0 \ot f \in \Sc (A)_0 \qquad \fl \, f \in \Sc (A_{\ify})$.

\smallskip

Now let $\ell$ be the number of infinite places of $k$ and consider the 
map $\rho 
: (\Rb_+^*)^{\ell} \ra \Rb_+^*$ given by
$$
\rho (\lb_1 , \ldots , \lb_{\ell}) = \lb_1 \ldots \lb_{\ell} \, . \leqno 
(21)
$$
As soon as $\ell > 1$ this map is not proper. Given a smooth function 
with 
compact support, $b \in C_c^{\ify} (\Rb_+^*)$ we need to find $a \in 
C_c^{\ify} 
((\Rb_+^*)^{\ell})$ such that the direct image of the measure $a(x) \, 
d^* x$ is 
$b(y) \, d^* y$ where $d^* x = \Pi \, d^* x_i$ is the product of the 
multiplicative Haar measures.

\smallskip

Equivalently one is dealing with a finite dimensional vector space $E$ 
and a 
linear form $L : E \ra \Rb$. One is given $b \in C_c^{\ify} (\Rb)$ and 
asked to 
lift it. One can write $E = \Rb \ts E_1$ and the lift can be taken as $a 
= b 
\ot 
b_1$ where $b_1 \in C_c^{\ify} (E_1)$, $\int b_1 \, dx = 1$.
\smallskip
 Thus we can in 
(19) 
take a function $f$ of the form,
$$
f(x) = g(x) \, \ov{\Xc}_{0,\ify} \, (x) \leqno (22)
$$
where the function $g \in C_c^{\ify}(A_{\ify}) $ only depends upon 
$(\vert x \vert_v)$, $v \in S_{ \ify}$ and is 
smooth with compact support, disjoint from the closed set 
$$
\left\{ x \in {\displaystyle \prod_{v \in S_{ \ify} }} k_v 
\, 
; \, \exists \, v \, , \, x_v = 0 \right\}.
$$

\smallskip

Thus, to any function $b \in C_c^{\ify} (\Rb_+^*)$ we can assign a test 
function $f = f_b$ such that for any $s$ ($\Re \, s > 0$)
$$
\lgl \D'_s , f_0 \ot f_b \rgl = \int_{\Rb_+^*} b(x) \, \vert x \vert^s 
\, d^* x 
\, . \leqno (23)
$$
By lemma 3, we get,
$$
\displaylines{
\left\lgl \int \D_{{1 \over 2} + it} \, \wh{\psi} (t) \, dt \, , \, f_0  
\ot f_b 
\right\rgl = \left\lgl \int L \, (\Xc_0 , {1 \over 2} + it ) 
\, 
\D'_{ {1 \over 2} + it} \, \wh{\psi} (t) \, dt \, , \, f_0  \ot f_b 
\right\rgl \cr
= \int \int L \, (\Xc_0 ,{1 \over 2} + it) \, \wh{\psi} (t) \, b(x) \, 
\vert x \vert^{{1 \over 2} + it} 
\, d^* x \, dt \, . \cr
}
$$
Thus, from (14) and (15) we conclude, using arbitrary test 
functions $b$ that the Fourier transform of the distribution $L \, 
(\Xc_0 ,1/2 + it) \, \wh{\psi} (t)$ actually vanishes,
$$
L \, (\Xc_0 ,{1 \over 2} + it) \, \wh{\psi} (t)= 0 \leqno (24)
$$
 \smallskip

To justify the above equality, we need to control the growth of the $L$ 
function in the variable $t$. One has,
$$
\vert L ({1 \over 2} + it) \vert = 0 \, ( \vert t \vert^N) \qquad 
\, 
. \leqno (25)
$$
In particular, since $L \left( {1 \over 2} + it \right)$ is an analytic 
function 
of $t$ we see that it is a multiplier of the algebra $\Sc (\Rb)$ of 
Schwartz 
functions in the variable $t$. Thus the product $L \, ({1 \over 2} + it) 
\, \wh{\psi} (t) $ is still a tempered 
distribution, and so is its Fourier transform. To say that the latter 
vanishes when 
tested on arbitrary functions which are smooth with compact support 
implies 
that
it vanishes.

\smallskip
The above argument uses the hypothesis  $\Xc_0 / C_{k,1} \not= 1$.

In the case $\Xc_0 / C_{k,1}= 1$ we need to impose to the test function 
$f$ 
used in (22) the condition $\int f \, dx = 0$ which means
$$
\int b(x) \, \vert x \vert \, d^* x = 0 \, . \leqno (26)
$$
But the space of functions $b(x) \, \vert x \vert^{1/2} \in C_c^{\ify} 
(\Rb_+^*)$ such that (26) holds is still dense in the Schwartz space 
$\Sc 
(\Rb_+^*)$.
\smallskip

To understand the equation (24), let us consider an equation for 
distributions $\a (t)$ of the form
$$
\vp (t) \, \a (t) = 0 \leqno (27)
$$
where we first work with distributions $\a$ on $S^1$ and we assume that 
$\vp 
\in 
C^{\ify} (S^1)$ has finitely many zeros $x_i \in Z (\vp)$, of finite 
order $n_i$. Let $J$ be the ideal 
 of $C^{\ify} (S^1)$ generated by $\vp$. One has $\psi \in J \Lra$ order 
of 
$\psi$ at $x_i$ is $\geq n_i $. 
\smallskip

Thus the distributions $\d_{x_i}$, 
$\d'_{x_i}, \ldots , \d_{x_i}^{(n_i - 1)}$ form a basis of the space of
solutions of (27).
\smallskip

Now $\wh{\psi} (t)$ is, for $\eta $ orthogonal to $\Im(E)$ and 
satisfying (16), a distribution with compact support, and 
$L \left(\Xc_0 , {1 
\over 2} + it \right) \, \wh{\psi} (t) =0$ . Thus by the above argument 
 we get that $\wh{\psi}$ is a finite linear combination of the 
distributions,
$$
\d_t^{(k)} \ , \ L \left( \Xc_0 , {1 \over 2} + it \right) = 0 \ , \ k < 
\hbox{order of the zero}, \ k < {\d - 1 \over 2} \, . \leqno (28)
$$
The condition $k < $ order of the zero is necessary and sufficient to 
get the vanishing on the range of $E$. The condition $k < {\d - 1 \over 
2}$ 
is necessary and sufficient to ensure that $\psi$ belongs to 
$L_{-\d}^2$, i.e. 
that
$$
\int (\log \vert x \vert)^{2k} \, (1 + \vert \log \vert x \vert 
\vert^2)^{-\d / 
2} \, d^* x < \ify \leqno (29)
$$
which is $2k + \d < -1$, i.e. $k < {\d - 1 \over 2}$.

\smallskip

Conversely, let $s$ be a zero of $L(\Xc_0 , s)$ and $k > 0$ its order. 
By lemma 3 and the 
finiteness and analyticity of $\D'_s$ (for $\Re \, s > 0$) we get
$$
\left( {\part \over \part s} \right)^a \, \D_s (f) = 0 \qquad \fl \, f 
\in 
\Sc 
(A)_0 \, , \ a = 0,1,\ldots ,k-1 \, . \leqno (30)
$$
\smallskip

We can differentiate the equality of lemma 3 and get,
$$
\left( {\part \over \part s} \right)^a \, \D_s (f) = \int_{C_k} E (f) 
(x) 
\, 
\Xc_0 (x) \, \vert x \vert^{s-1/2} \, (\log \vert x \vert)^a \, d^* x \, 
. \leqno 
(31)
$$
Thus $\eta$ belongs to the orthogonal of $\Im(E)$ and satisfies (16) iff 
it is a finite linear combination of 
functions of the form,
$$
\eta_{t,a} (x)= \Xc_0(x) \,\vert x \vert^{it} \, (\log \vert x \vert)^a 
, \leqno (32)
$$
where, 
$$
 L \left( \Xc_0 , {1 \over 2} + it \right) = 0 \ , \,\,\,\ a <
\hbox{order of the zero}, \ a < {\d - 1 \over 2} \, . \leqno (33)
$$
The restriction to the subgroup $\Rb_+^*$ of $C_k$ of the transposed of 
$W$ is thus given in the above basis by:
$$
W (\lb)^t \, \eta_{t,a}  =  \sum_{b=0}^{a} \, C_a^b \,\lb^{it} \,(Log( 
\lb ))^b \, \eta_{t,a-b}  . \leqno (34)
$$

 The multiplication operator by a function with bounded derivatives is a 
bounded operator in any Sobolev space
thus one checks directly, using the density in the orthogonal of 
$\Im(E)$ of vectors satisfying (16), that if $L \left( \Xc_0 , {1 \over 
2} + is \right) \not=  0 $ then $is$ does not belong to the spectrum of
$D_{\Xc_0}^t$.
\smallskip
 This determines the spectrum of the operator $D_{\Xc_0}^t$ and hence of 
its transpose $D_{\Xc_0}$
as indicated in Theorem 1 and ends the proof of theorem 1.
\smallskip
 Let us now prove the corollary.
Let us fix $h_0 \in \Sc (C_k) $
such that $ \wh h_0 $ has compact support contained in $\{\Xc_0 \} \ts 
\Rb$ and $ \wh h_0 ( \Xc_0 , s)=1$
 for $s$ small.
\smallskip
  Let then $h_s$ be given by $ h_s(g) = \, h_0 (g)\,\vert g \vert^{is}$. 
The Fourier transform $\wh h_s$ is then the translate
of $ \wh h_0 $, and one can choose $h_0$ such that,
$$
\sum_{n \in \Zb} \, \wh h_n (\Xc_0\, ,u) =1 \qquad 
 u \in \Rb \leqno (35)
$$

 When $\vert s\vert \ra \ify $, the dimension of the range of $W^t(h_s)$
is of the order of $Log\vert s\vert$ as is the number of zeros of 
the $L$ function in the translates of a fixed interval (cf [W 3]).
\smallskip
Let $h \in \Sc (C_k)$. One has $W^t(h) = \sum_{n \in \Zb} \, W^t(h\, * 
h_n) $. 

It follows then from the polynomial growth of the norm of $W^t(g)$ that
 the operator
$$
\int h(g) \,   W(g)^t \, d^* g \leqno (36)
$$
is of trace class for any $h \in \Sc (C_k)$.
\smallskip
Moreover using the triangular form given by (34) we get its trace, and 
hence the trace of its transpose
$W(h)$ as,
$$
\hbox{Trace} \, W(h) = \sum_{{L\left( \Xc , {1\over 2} + \rho
\right) =0 \atop \rho \in i \, \Rb}} \wh h (\Xc ,\rho) \leqno (37)
$$
where the multiplicity is counted as in Theorem 1 and where
the Fourier transform $\wh h$ of $h$ is defined by,
$$
\wh h (\Xc ,\rho) = \int_{C_k} h(u) \, \wt{\Xc} (u) \, \vert u
\vert^\rho \, d^* \, u \, .\leqno (38)
 $$ 

\vglue 1cm

\vglue 1cm
\noindent {\bf Appendix II. Explicit formulas} 

\smallskip

Let us first recall the Weil explicit formulas ([W3]). One lets $k$ be 
a global field. One identifies the quotient $C_k / C_{k,1}$ with the 
range of the module,
$$
N = \{ \vert g \vert \, ; \, g \in C_k \} \sbs \Rb_+^* \, . \leqno 
(1)
$$
One endows $N$ with its normalized Haar measure $d^* x$. Given a 
function $F$ on $N$ such that, for some $b > {1 \over 2}$,
$$
\vert F (\nu) \vert = 0 (\nu^b) \quad \nu \ra 0 \ , \ \vert F (\nu) 
\vert = 0 (\nu^{-b}) \ , \ \nu \ra \ify \, , \leqno (2)
$$
one lets,
$$
\Phi (s) = \int_N F(\nu) \, \nu^{1/2 - s} \, d^* \nu \, . \leqno (3)
$$
Given a Gr\"ossencharakter $\Xc$, i.e. a character of $C_k$ and any 
$\rho$ in the strip $0 < {\rm Re} (\rho) < 1$, one lets $N(\Xc , 
\rho)$ be the order of $L (\Xc , s)$ at $s = \rho$. One lets,
$$
S (\Xc , F) = \sum_{\rho} N (\Xc , \rho) \, \Phi (\rho) \leqno (4)
$$
where the sum takes place over $\rho$'s in the above open strip. One 
then defines a distribution $\D$ on $C_k$ by,
$$
\D = \log \vert d^{-1} \vert \, \d_1 + D - \sum_v D_v \, , \leqno (5)
$$
where $\d_1$ is the Dirac mass at $1 \in C_k$, where $d$ is a 
differental idele of $k$ so that $\vert d \vert^{-1}$ is up to sign 
the discriminant of $k$ when char $(k) = 0$ and is $q^{2g-2}$ when 
$k$ is a function field over a curve of genus $g$ with coefficients 
in the finite field $\Fb_q$.

\smallskip

The distribution $D$ is given by,
$$
D(f) = \int_{C_k} f(w) \, (\vert w \vert^{1/2} + \vert w 
\vert^{-1/2}) \, d^* w \leqno (6)
$$
where the Haar measure $d^* w$ is normalized (cf. IIb).
The distributions $D_v$ are parametrized by the places $v$ of $k$ and 
are obtained as follows. For each $v$ one considers the natural 
proper homomorphism,
$$
k_v^* \ra C_k \ , \ x \ra \hbox{class of} \ (1,\ldots ,x,1 \ldots) 
\leqno (7)
$$
of the multiplicative group of the local field $k_v$ in the idele 
class group $C_k$.

\smallskip

One then has,
$$
D_v (f) = Pfw \int_{k_v^*} {f(u) \over \vert 1-u \vert} \, \vert u 
\vert^{1/2} \, d^* u \leqno (8)
$$
where the Haar measure $d^* u$ is normalized (cf. IIb), and where the 
Weil Principal value $Pfw $ of the integral is obtained as follows, for 
a 
local field $K = k_v$,
$$
Pfw \int_{k_v^*} \, 1_{R_v^*} \, {1 \over \vert 1-u \vert} \, d^* u = 
0 \, , \leqno (9)
$$
if the local field $k_v$ is non Archimedean, and otherwise:
$$
Pfw \int_{k_v^*} \vp (u) \, d^* u = PF_0 \int_{\Rb_+^*} \psi (\nu) \, 
d^* \nu \, , \leqno (10)
$$
where $\psi (\nu) = \int_{\vert u \vert = \nu} \vp (u) \, d_{\nu} u$ 
is obtained by integrating $\vp$ over the fibers, while
$$
PF_0 \int \psi (\nu) d^* \nu = 2 \log (2\pi) \, c + \lim_{t \ra \ify} 
\left( \int (1-f_0^{2t} ) \, \psi (\nu) \, d^* \nu - 2c \log t 
\right) \, , \leqno (11)
$$
where one assumes that $\psi - c \, \, f_1^{-1}$ is integrable on 
$\Rb_+^*$, and
$$
f_0 (\nu) = \inf (\nu^{1/2} , \nu^{-1/2}) 
\qquad \fl \, \nu \in \Rb_+^* , \qquad \,
f_1 = f_0^{-1} - f_0  \  .
$$
The Weil explicit formula is then,

\medskip

\noindent {\bf Theorem 1.} ([W]) {\it  With the above notations one has 
$S (\Xc , F) = \D (F (\vert w \vert)$ $\Xc (w))$.
}

\medskip

We shall now elaborate on this formula and in particular compare the 
principal values $Pfw$ with those of theorem V.3. 

\smallskip

Let us make the following change of variables,
$$
\vert g \vert^{-1/2} \, h (g^{-1}) = F (\vert g \vert) \, \Xc_0 (g) 
\, , \leqno (12)
$$
and rewrite the above equality in terms of $h$. 

\smallskip

By (3) one has,
$$
\Phi \left( {1 \over 2} + is \right) = \int_{C_k} F (\vert g \vert) 
\, \vert g \vert^{-is} \, d^* g \, , \leqno (13)
$$
thus, in terms of $h$,
$$
\int h(g) \, \Xc_1 (g) \, \vert g \vert^{1/2+is} \, d^* g = \int F 
(\vert g^{-1} \vert) \, \Xc_0 (g^{-1}) \, \Xc_1 (g) \, \vert g 
\vert^{is} 
\, d^* g \, , \leqno (14)
$$
which is equal to 0 if $\Xc_1 / C_{k,1} \not= \Xc_0 / C_{k,1}$ and for 
$\Xc_1 = \Xc_0$,
$$
\int h(g) \, \Xc_0 (g) \, \vert g \vert^{1/2+is} \, d^* g = \Phi \left( 
{1 \over 2} + is \right) \, . \leqno (15)
$$
Thus, with our notations we see that,
$$
{\rm Supp} \, \wh h \sbs \Xc_0 \ts \Rb \ , \ \wh h (\Xc_0 ,\rho) = \Phi 
(\rho) \, . \leqno (16)
$$
Thus we can write,
$$
S (\Xc_0 , F) = \sum_{L(\Xc , \rho) = 0 , \Xc \in \wh{C}_{k,1} \atop 0 
< {\rm Re} \, \rho < 1 \hfill} \wh h (\Xc , \rho) \leqno (17)
$$
using a fixed decomposition $C_k = C_{k,1} \ts N$.

\smallskip

Let us now evaluate each term in (5).

\smallskip

The first gives $(\log \vert d^{-1} \vert) \, h(1)$. One has, using (6) 
and 
(12),
$$
\displaylines{
\lgl D , F (\vert g \vert) \, \Xc_0 (g) \rgl = \int_{C_k} \vert g 
\vert^{-1/2} \, h (g^{-1}) \, (\vert g \vert^{1/2} + \vert g 
\vert^{-1/2}) \, d^* g \cr
= \int_{C_k} h(u) \, (1+ \vert u \vert) \, d^* u = \wh h (0) + \wh h 
(1) \, , \cr
}
$$
where for the trivial character of $C_{k,1}$ one uses the notation
$$
\wh h (z) = \wh h (1,z) \qquad \fl \, z \in \Cb \, . \leqno (18)
$$
Thus the first two terms of (5) give
$$
(\log \vert d^{-1} \vert) \, h(1) + \wh h (0) + \wh h (1) \, . \leqno 
(19)
$$
Let then $v$ be a place of $k$, one has by (8) and (12),
$$
\lgl D_v , F (\vert g \vert) \, \Xc_0 (g) \rgl = Pfw \int_{k_v^*} 
{h(u^{-1}) \over \vert 1-u \vert} \, d^* u \, . 
$$
We can thus write the contribution of the last term of (5) as,
$$
- \sum_v Pfw \int_{k_v^*} {h(u^{-1}) \over \vert 1-u \vert} \, d^* u \, 
. \leqno (20)
$$
Thus the equality of Weil can be rewritten as,
$$
\wh h (0) + \wh h (1) - \sum_{L(\Xc , \rho) = 0 , \Xc \in \wh{C}_{k,1} 
\atop 0 < {\rm Re} \, \rho < 1 \hfill} \wh h (\Xc , \rho) = (\log \vert 
d \vert) \, h(1) + \leqno (21)
$$
$$
\sum_v Pfw \int_{k_v^*} {h(u^{-1}) \over \vert 1-u \vert} \, d^* u \, . 
$$
Which now holds for finite linear combinations of functions $h$ of the 
form (12).

\smallskip

This is enough to conclude when $h(1) = 0$.

\smallskip

Let us now compare the Weil Principal values, with those dictated by 
theorem V.3. We first work with a local field $K$ and compare (9), (10) 
with our prescription. Let first $K$ be non Archimedean. Let $\a$ be a 
character of $K$ such that,
$$
\a / R = 1 \ , \ \a / \pi^{-1} R \not= 1 \, . \leqno (22)
$$
Then, for the Fourier transform given by,
$$
(Ff) (x) = \int f(y) \, \a (y) \, dy \, , \leqno (23)
$$
with $dy$ the selfdual Haar measure, one has
$$
F(1_R) = 1_R \, . \leqno (24)
$$

\medskip

\noindent {\bf Lemma 2.} {\it With the above choice of $\a$ one has
$$
\int' {h(u^{-1}) \over \vert 1-u \vert} \, d^* u  = Pfw \int {h(u^{-1}) 
\over \vert 1-u \vert} \, d^* u
$$
with the notations of theorem 3.
}

\medskip

\noindent {\it Proof.} By construction the two sides can only differ by 
a multiple of $h(1)$. Let us recall from theorem 3 that the left hand 
side is given by
$$
\left\lgl L ,  \, {h(u^{-1}) \over \vert u \vert} \right\rgl 
\, , \leqno (25)
$$
where $L$ is the unique extension of $ \rho^{-1}{du \over \vert 1-u 
\vert}$ whose 
Fourier transform vanishes at 1, $\wh L (1) = 0$. Thus from (9) we just 
need to check that (25) vanishes for $h = 1_{R^*}$, i.e. that
$$
\lgl L , 1_{R^*} \rgl = 0 \, . \leqno (26)
$$
Equivalently, if we let $Y = \{ y \in K \, ; \, \vert y - 1 \vert = 
1\}$ we just need to show, using Parseval, that,
$$
\lgl \log \vert u \vert \, , \, \wh{1}_Y \rgl = 0 \, . \leqno (27)
$$
One has $\wh{1}_Y (x) = \int_Y \a (xy) \, dy = \a (x) \, \wh{1}_{R^*} 
(x)$, and $1_{R^*} = 1_R - 1_P$, $\wh{1}_{R^*} = 1_R - \vert \pi \vert 
\, 1_{\pi^{-1} R}$, thus, with $q^{-1}=\vert \pi \vert $,
$$
\wh{1}_Y (x) = \a (x) \left( 1_R -{1 \over q} \, 1_{\pi^{-1} R} \right) 
(x) \, . \leqno (28)
$$
We now need to compute $\int \log \vert x \vert \, \wh{1}_Y (x) \, dx = 
A+B$,
$$
A = -{1 \over q} \int_{\pi^{-1} R^*} \a (x) \, (\log q) \, dx \ , \ B = 
\left( 1 - {1 \over q} \right) \int_R \log \vert x \vert \, dx \, . 
\leqno (29)
$$
Let us show that $A+B = 0$. One has $\int_R dx = 1$, and
$$
\displaylines{
A = - \int_{R^*} \a (\pi^{-1} y) \, (\log q) \, dy = - \log q \left( 
\int_R \a (\pi^{-1} y) \, dy - \int_P dy \right) \cr
= {1 \over q} \, \log q \ , \ \hbox{since} \ \int_R \a (\pi^{-1} y) \, 
dy = 0 \ \hbox{as} \ \a / \pi^{-1} R \not= 1 \, .
}
$$
To compute $B$, note that $\int_{\pi^n R^*} dy = q^{-n} \left( 1 - {1 
\over q} \right)$ so that
$$
B = \left( 1 - {1 \over q} \right)^2 \sum_{n=0}^{\ify} (-n \log q) \, 
q^{-n} = -q^{-1} \log q \, .
$$
and $A+B=0$.
\hfill \xx

\medskip

Let us now treat the case of Archimedean fields. We take $K = \Rb$ 
first, and we normalize the Fourier transform as,
$$
(Ff) (x) = \int f(y) \, e^{-2\pi i xy} \, dy \leqno (30)
$$
so that the Haar measure $dx$ is selfdual.

\smallskip

With the notations of (10) one has,
$$
Pfw \int_{\Rb^*} f_0^3 (\vert u \vert) \, {\vert u \vert^{1/2} \over 
\vert 1-u \vert} \, d^* u = \log \pi + \g \leqno (31)
$$
where $\g$ is Euler's constant, $\g = - \G^{'}(1)  $. Indeed integrating 
over the fibers gives $f_0^4 \ts ( 1 - 
f_0^{4} )^{-1}$, and one gets,
$$
\displaylines{
  PF_0  \int_{\Rb_+^*} f_0^4 \ts ( 1 - 
f_0^{4} )^{-1} \, d^* u =  
\Biggl(  \log (2\pi) + \lim_{t \ra \ify} \Biggl( \int_{\Rb_+^*} ( 1 - 
f_0^{2t} ) \, f_0^4  ( 1 - 
f_0^{4} )^{-1} \, d^* u \cr
- \log t \Biggl) \Biggl) = \log 2\pi + \g - \log 2 \, . \cr
}
$$
\smallskip

Now let $\vp (u) = -\log \vert u \vert$, it is a tempered distribution 
on $\Rb$ and one has,
$$
\lgl \vp , e^{-\pi u^2} \rgl = {1 \over 2} \, \log \pi + {\g \over 2} + 
\log 2 \, , \leqno (32)
$$
as one obtains from ${\part \over \part s} \int \vert u \vert^{-s} \, 
e^{-\pi u^2} \, du = {\part \over \part s} \left( \pi^{{s-1 \over 2}} 
\, \G \left( {1-s \over 2} \right) \right)$ evaluated at $s=0$, using 
$ {\G^{'} ({1 \over 2}) \over \G ({1 \over 2})} = \, - \g -2 \log2 $.
\smallskip

Thus by the Parseval formula one has,
$$
\lgl \wh{\vp} , e^{-\pi x^2} \rgl = {1 \over 2} \log \pi + {\g \over 2} 
+ \log 2 \, , \leqno (33)
$$
which gives, for any test function $f$,
$$
\lgl \wh{\vp} , f \rgl = \lim_{\ve \ra 0} \left( \int_{\vert x \vert 
\geq \ve} f(x) \, d^* x + (\log \ve) \, f(0) \right) + \lb \, f(0) 
\leqno (34)
$$
where $\lb = \log (2\pi) + \g$. In order to get (34) one uses 
the equality,
$$
\lim_{\ve \ra 0} \left( \int_{\vert x \vert \geq \ve} f(x) \, d^* x + 
(\log \ve) \, f(0) \right) = \lim_{\ve \ra 0} \left( \int f(x) \, \vert 
x \vert^{\ve} \, d^* x - {1 \over \ve} \, f(0) \right) \, , \leqno (35)
$$
which holds since both sides vanish for $f(x) = 1$ if $\vert x \vert 
\leq 1$, $f(x) = 0$ otherwise. 
\smallskip
Thus from (34) one gets,
$$
\int'_{\Rb} f(u) \, {1 \over \vert 1-u \vert} \, d^* u = \lb \, f(1) + 
\lim_{\ve \ra 0} \left( \int_{\vert 1-u \vert \geq \ve} {f(u) \over 
\vert 1-u \vert} \, d^* u + ( \log \ve )\, f(1) \right) \, . \leqno (36)
$$
Taking $f(u) = \vert u \vert^{1/2} \, f_0^3 ( \vert u \vert)$, the 
right hand side of (36) gives $\lb - \log 2 = \log \pi + \g$, thus we 
conclude using (31) that for any test function $f$,
$$
\int'_{\Rb} f(u) \, {1 \over \vert 1-u \vert} \, d^* u = Pfw \int_{\Rb} 
f(u) \, {1 \over \vert 1-u \vert} \, d^* u \, . \leqno (37)
$$
Let us finally consider the case $K = \Cb$. We choose the basic 
character $\a$ as
$$
\a (z) = \exp 2\pi i (z + \ov z ) \, , \leqno (38)
$$
the selfdual Haar measure is $dz \, d\ov z = \vert dz \wdg d \ov z 
\vert$, and the function $f(z) = \exp -2\pi \vert z \vert^2$ is 
selfdual.

\smallskip

The normalized multiplicative Haar measure is
$$
d^* z = {\vert dz \wdg d \ov z \vert \over 2\pi \vert z \vert^2} \, . 
\leqno (39)
$$
Let us compute the Fourier transform of the distribution
$$
\vp (z) = -\log \vert z \vert_{\Cb} = -2 \log \vert z \vert \, . \leqno 
(40)
$$
One has
$$
\lgl \vp , \exp - 2\pi \vert z \vert^2 \rgl = \log 2\pi + \g \, , 
\leqno (41)
$$
as is seen using ${\part \over \part \ve} \left( \int e^{-2\pi \vert z 
\vert^2} \, \vert z \vert^{-2\ve} \, \vert dz \wdg d \ov z \vert 
\right) = {\part \over \part \ve} \, ((2\pi)^{\ve} \, \G (1-\ve))$.

\smallskip

Thus $\lgl \wh{\vp} , \exp -2\pi \vert u \vert^2 \rgl = \log 2\pi + \g$ 
and one gets,
$$
\lgl \wh{\vp} , f \rgl = \lim_{\ve \ra 0} \left( \int_{\vert u 
\vert_{\Cb} \geq \ve} f(u) \, d^* u + \log \ve \, f(0) \right) + \lb' \, 
f(0) \leqno (42)
$$
where $\lb' = 2 (\log 2\pi + \g)$.

\smallskip

To see this one uses the analogue of (35) for $K= \, \Cb $, to compute 
the right hand side of (42) for 
$f(z) = \exp - 2\pi \vert z \vert^2$.

\smallskip

Thus, for any test function $f$, one has,
$$
\int'_{\Cb} f(u) \, {1 \over \vert 1-u \vert_{\Cb}} \, d^* u = \lb' \, 
f(1) + \lim_{\ve \ra 0} \left( \int_{\vert 1-u \vert_{\Cb} \geq \ve} 
{f(u) \over \vert 1-u \vert_{\Cb}} \, d^* u + (\log \ve ) f(1)\right) \, 
. 
\leqno (43)
$$
Let us compare it with $Pfw$. When one integrates over the fibers of 
$\Cb^* \build \longrightarrow_{}^{\vert \ \vert_{\Cb}} \Rb_+^*$ the 
function $\vert 1-z \vert_{\Cb}^{-1}$ one gets,
$$
{1 \over 2\pi} \int_0^{2\pi} {1 \over \vert 1 - e^{i\t} z\vert^2} \, 
d\t = {1 \over 1 - \vert z \vert^2} \ \hbox{if} \ \vert z \vert < 1 \, 
, \ \hbox{and} \ {1 \over \vert z \vert^{2} - 1 } \ \hbox{if} \ \vert z 
\vert > 1 \, . \leqno (44)
$$
Thus for any test function $f$ on $\Rb_+^*$ one has, by (10),
$$
Pfw \int f(\vert u \vert_{\Cb}) \, {1 \over \vert 1 - u \vert_{\Cb}} \, 
d^* u = PF_0 \int f(\nu) \, {1 \over \vert 1 - \nu \vert}\, d^* \nu 
\leqno (45)
$$
with the notations of (11). With $f_2(\nu) = \nu^{1 \over 2} \,f_0( \nu 
)$ we thus get, using (11),
$$
Pfw \int f_2 (\vert u \vert_{\Cb}) \, {1 \over \vert 1 - u 
\vert_{\Cb}} \, d^* u = PF_0 \int f_0 \,f_1^{-1}   d^* \nu=  2 (\log 
2\pi + \g) \, . \leqno (46)
$$
We shall now show that,
$$
\lim_{\ve \ra 0} \left( \int_{\vert 1 - u \vert_{\Cb} \geq \ve} {f_2 
(\vert u \vert_{\Cb}) \over \vert 1 - u \vert_{\Cb}} \, d^* u + \log 
\ve \right) = 0 \, , \leqno (47)
$$
it will then follow that, using (43),
$$
\int'_{\Cb} f(u) \, {1 \over \vert 1-u \vert_{\Cb}} \, d^* u = Pfw 
\int f(u) \, {1 \over \vert 1-u \vert_{\Cb}} \, d^* u \, . \leqno (48)
$$
To prove (47) it is enough to investigate the integral,
$$
\int_{\vert z \vert \leq 1 , \vert 1-z \vert \geq \ve } ((1-z) (1- \ov z 
))^{-1} \, \vert dz \wdg d \ov z \vert = j (\ve) \leqno (49)
$$
and show that $j(\ve) = \a \log \ve + o(1)$ for $\ve \ra 0$. A similar 
statement then holds for
$$
\int_{\vert z \vert \leq 1 , \vert 1-z^{-1} \vert \geq \ve } ((1-z) (1- 
\ov z 
))^{-1} \, \vert dz \wdg d \ov z \vert .
$$

\smallskip

One has $j(\ve) = \int_D \vert dZ \wdg d \ov Z \vert$, where $Z = \log 
(1-z)$ and the domain $D$ is contained in the rectangle,
$$
\{ Z = (x+iy) \, ; \, \log \ve \leq x \leq \log 2 \, , \, -{\pi \over 
2} \leq y \leq \pi / 2 \} = R_{\ve} \leqno (50)
$$
and bounded by the curve $x = \log (2 \cos y)$ which comes from the 
equation of the circle $\vert z \vert = 1$ in polar coordinates centered 
at $z=1$. One thus gets,
$$
j(\ve) = 4 \int_{\log \ve}^{\log 2} {\rm Arc} \cos (e^x / 2) \, dx \, , 
\leqno (51)
$$
when $\ve \ra 0$ one has $j(\ve) \sim 2\pi \log (1/\ve)$, which is the 
area of the following rectangle (in the measure $\vert dz \wdg d \ov z 
\vert$),
$$
\{ Z = (x+iy) \, ; \, \log \ve \leq x \leq 0 \, , \, -\pi / 2 \leq y 
\leq \pi / 2 \} \leqno (52)
$$
One has $\vert R_{\ve} \vert - 2\pi \log 2 = 2\pi \log 
(1/\ve)$. When $\ve \ra 0$ the area of $R_{\ve} \backslash D$ converges 
to
$$
4 \int_{-\ify}^{\log 2} {\rm Arc} \sin (e^x / 2) \, dx = -4 \int_0^{\pi 
/2} \log (\sin u) \, du = 2\pi \log 2 \, , \leqno (53)
$$
so that $j(\ve) = 2\pi \log (1/\ve) + o(1)$ when $\ve \ra 0$.

\smallskip

Thus we can assert that with the above choice of basic characters for 
local fields one has, for any test function $f$,
$$
\int'_K f(u) \, {1 \over \vert 1-u \vert} \, d^* u = Pfw \int f(u) \, 
{1 \over \vert 1-u \vert} \, d^* u \, . \leqno (54)
$$

\medskip

\noindent {\bf Lemma 3.} {\it Let $K$ be a local field, $\a_0$ a 
normalized character as above and $\a$, $\a (x) = \a_0 (\lb x)$ an 
arbitrary character of $K$. Let $\int'$ be defined as in theorem V.3 
relative to $\a$, then, for any test function $f$,
$$
\int'_K f(u) \, {1 \over \vert 1-u \vert} \, d^* u = \log \vert \lb 
\vert \, f(1) + Pfw \int f(u) \, {1 \over \vert 1-u \vert} \, d^* u \, 
.
$$
}

\medskip

\noindent {\it Proof.} The new selfdual Haar measure is
$$
da = \vert \lb \vert^{1/2} \, d_0 \, a \leqno (55)
$$
with $d_0 \, a$ selfdual for $\a_0$.

\smallskip

Similarly the new Fourier transform is given by
$$
\wh f (x) = \int \a (xy) \, f(y) \, dy = \int \a_0 (\lb \, xy) \, f(y) 
\, \vert \lb \vert^{1/2} \, d_0 \, y \, ,
$$
thus
$$
\wh f (x) = \vert \lb \vert^{1/2} \, \wh{f}^0 (\lb \, x) \, . \leqno 
(56)
$$
Let then $\vp (u) = -\log \vert u \vert$. Its Fourier transform as a 
distribution is given by,
$$
\lgl \wh{\vp} , f \rgl = \int (-\log \vert u \vert) \, \wh f (u) \, du 
\, . \leqno (57)
$$
One has 
$$
\int (-\log \vert u \vert) \, \wh f (u) \, du = \int (-\log \vert u 
\vert) \, \wh{f}^0 (\lb \, u) \, \vert \lb \vert \, d_0 \, u 
$$
$$
= \int (-\log \vert v \vert) \, \wh{f}^0 (v) \, d_0 \, v + \int \log 
\vert \lb \vert \, \wh{f}^0 (v) \, d_0 \, v \, .
$$
$$
= \int (-\log \vert v \vert) \, \wh{f}^0 (v) \, d_0 \, v +  \log 
\vert \lb \vert \,f(0) \,  \, .
$$

Thus the lemma follows from (54).~\xx

\medskip

Let us now pass to the global case, recall that if $\a$, $\a \not= 1$, 
is a character of $A$ such that $\a / k = 1$, there exists a 
differental idele $d = (d_v)$ such that, (cf. [W1])
$$
\a_v (x) = \a_{0,v} (d_v \, x) \leqno (58)
$$
where $\a = \Pi \, \a_v$ and each local character $\a_{0,v}$ is 
normalized as above.

\smallskip

We can thus rewrite the Weil formula (theorem 1) as,

\medskip

\noindent {\bf Theorem 6.} {\it Let $k$ be a global field, $\a$ a non 
trivial character of $A/k$ and $\a = \Pi \, \a_v$ its local factors.

\smallskip

Let $h \in \Sc (C_k)$ have compact support, then
$$
\wh h (0) + \wh h (1) - \sum_{L(\Xc , \rho) = 0 \atop 0 < {\rm Re} \rho 
< 1} \wh h (\Xc , \rho) = \sum_v \int'_{k_v^*} {h(u^{-1}) \over \vert 
1-u \vert} \, d^* u
$$
where the normalization of $\int'$ is given by $\a_v$ as in theorem V.3, 
and $\wh h (\Xc , z) = \int h(u) \, \Xc (u) \, \vert u \vert^z \, d^* 
u$.
}

\medskip

\noindent {\it Proof.} This follows from formula (21), lemma 3 and the 
equality $\log \vert d \vert = \, \build \sum_{v}^{} \log \vert d_v 
\vert$.~\xx
\medskip

\noindent {\bf Normalization of Haar measure on modulated group}   
\smallskip

We let $G$ be a locally compact abelian group with a proper morphism,
$$
g \ra \vert g \vert \, , \ G \ra \Rb_+^* \leqno (1)
$$
whose range is cocompact in $\Rb_+^*$. 

\smallskip

There exists a unique Haar measure $d^* g$ on $G$ such that
$$
\int_{\vert g \vert \in [ 1 , \L ]} d^* g \sim \log \L \qquad 
\hbox{when} \ \L 
\ra + \ify \, . \leqno (2)
$$
Let $G_0 = \Ker \mod = \{ g \in G \, ; \, \vert g \vert = 1 \}$. It is a 
compact group by hypothesis, and one can identify $G / G_0$ with the 
range $N$ 
of the module. Let us determine the measure $d^* n$ on $N \sbs \Rb_+^*$ 
such 
that (2) holds for
$$
\int f \, d^* g = \int \left( \int f (n g_0) \, d g_0 \right) d^* n 
\leqno (3)
$$
where the Haar measure $d g_0$ is normalized by
$$
\int_{G_0} d g_0 = 1 \, . \leqno (4)
$$
We let $\rho_{\L}$ be the function on $G$ given by
$$
\rho_{\L} (g) = 0 \quad \hbox{if} \ \vert g \vert \notin [1 , \L] \, , \ 
\rho_{\L} (g) = { 1 \over \log \L} \quad \hbox{if} \ g \in [1 , \L] \, . 
\leqno 
(5)
$$
The normalization (2) means that $\int \rho_{\L} \, d^* g \ra 1$ when 
$\L \ra 
\ify$.

\smallskip

Let first $N = \Rb_+^*$ then the unique measure satisfying (2) is
$$
d^* \lb = {d\lb \over \lb} \, . \leqno (6)
$$
Let then $N = \mu^{\Zb}$ for some $\mu > 1$. Let us consider the measure
$$
\int f \, d^* g = \a \sum f (\mu^n) \, . \leqno (7)
$$
We take $f = \rho_{\L}$, then the right hand side is $\a \, {N \over 
\log \L}$ 
where $N$ is the number of $\mu^n \in [1 , \L]$, i.e. $\sim {\log \L 
\over \log 
\mu}$. This shows that (2) holds iff
$$
\a = \log \mu \, . \leqno (8)
$$
Let us show more generally that if $H \sbs G$ is a compact subgroup of 
$G$ and 
if both $d^* g$ and $d^* h$ are normalized by (2) one has
$$
\int \left( \int f (hy) \, d^* h \right) \, d_0 \, y = \int f \, d^* g 
\leqno 
(9)
$$
where $d_0 \, y$ is the Haar measure of integral 1 on $G / H$,
$$
\int_{G / H} d_0 \, y = 1 \, . \leqno (10)
$$
The left hand side of (9) defines a Haar measure on $G$ and we just need 
to show that it 
satisfies (2).
\smallskip
One has
$ \Vert \rho_{\L} (\cdot \, y) - \rho_{\L} \Vert_1 \ra 0 \,\,\, 
\hbox{when} \ \L 
\ra 
\ify $ ,
 and 
$$
\int \rho_{\L} (hy) \, d^* h \ra 1 
 \qquad \hbox{when}\,\L \ra \ify \, \leqno (11)
$$
 uniformly on compact sets of $y \in G $, thus
$$
\int \left( \int \rho_{\L} (hy) \, d^* h \right) \, d_0 \, y \ra 1 
\qquad 
\hbox{when}\, \ \L \ra \ify \, . \leqno (12)
$$

\vglue 1cm
\noindent {\bf Appendix III. Distribution trace formulas }

\noindent In this appendix we recall for the convenience of the reader 
the coordinate free treatment of
distributions of [GS] and give the details of the transversality 
conditions.\smallskip
 Given a vector space $E$ over $\Rb$, $\dim
E =n$, a density is a map, $\rho \in \vert E \vert$,
$$
\rho : \wedge^n \, E \ra \Cb \leqno (1)
$$
such that $\rho (\lb v) = \vert \lb \vert \, \rho (v) \qquad
\fl \, \lb \in \Rb, \qquad\fl \, v \in\wedge^n\, E  $.
\smallskip

\noindent Given a linear map $T:E\rightarrow F$ we let $\vert T \vert 
:\vert F \vert \rightarrow  \vert E \vert $ be the corresponding linear 
map,it depends contravariantly on $T$.
\smallskip

\noindent Given a manifold $M$ and $\rho \in C_c^{\ify} (M,
\vert TM \vert)$ one has a canonical integral,
$$
\int \rho \in \Cb \, . \leqno (2)
$$
Given a vector bundle $L$ on $M$ one defines the generalized
sections on $M$ as the dual space of $C_c^{\ify} (M,L^* \ot
\vert TM \vert)$
$$
C^{-\ify} (M,L) = \ \hbox{dual of} \ C_c^{\ify} (M,L^* \ot
\vert TM \vert) \leqno (3)
$$
where $L^*$ is the dual bundle. One has a natural inclusion,
$$
C^{\ify} (M,L) \sbs C^{-\ify} (M,L) \leqno (4)
$$
given by the pairing
$$
\s \in C^{\ify} (M,L) \ , \ s \in C_c^{\ify} (M,L^* \ot \vert
TM \vert) \ra \int \lgl s,\s \rgl \leqno (5)
$$
where $\lgl s,\s \rgl$ is viewed as a density, $\lgl s , \s
\rgl \in C_c^{\ify} (M,\vert TM \vert)$.

\smallskip

\noindent One has a similar notion of
generalized section with compact support.

\smallskip

\noindent Given a smooth map $\vp : X \ra Y$, then if $\vp$
is {\it proper}, it gives a (contravariantly) associated map
$$
\vp^* : C_c^{\ify} (Y,L) \ra C_c^{\ify} (X,\vp^* (L)) \ , \
(\vp^* \, \xi)(x) = \xi (\vp (x)) \leqno (6)
$$
where $\vp^* (L)$ is the pull back of the vector bundle $L$.

\smallskip

\noindent  Thus, given a linear form on $C_c^{\ify} (X,\vp^*
(L))$ one has a (covariantly) associated linear form on
$C_c^{\ify} (Y,L)$. In particular with $L$ trivial we see
that given a generalized density $\rho \in C^{-\ify} (X,\vert
T \vert)$ one has a pushforward
$$
\vp_* (\rho) \in C^{-\ify} (Y,\vert T \vert) \leqno (7)
$$
with $\lgl \vp_* (\rho) ,\xi \rgl = \lgl \rho , \vp^* \, \xi
\rgl \qquad \fl \, \xi \in C_c^{\ify} (X)$.

\smallskip

\noindent Next, if $\vp$ is a fibration and $\rho \in
C_c^{\ify} (X,\vert T \vert)$ is a density then one can
integrate $\rho$ along the fibers, the obtained density on
$Y$, $\vp_* (\rho)$ is given as in (7) by 
$$
\lgl \vp_* (\rho) ,f \rgl = \lgl \rho, \vp^* \, f \rgl \qquad \fl
\, f \in C^{\ify} (Y) \leqno (8)
$$
but the point is that it is not only a generalized section
but a smooth section $\vp_* (\rho) \in C_c^{\ify} (Y,\vert T
\vert)$.

\smallskip

\noindent It follows  that if $f\in C^{-\ify} (Y)$ is a
generalized function, then one obtains a generalized function
$\vp^* (f)$ on $X$ by,
$$
\lgl \vp^* (f) , \rho \rgl = \lgl f,\vp_* (\rho) \rgl \qquad
\fl \, \rho \in C_c^{\ify} (X,\vert T \vert) \, . \leqno (9)
$$
In general,the pullback $\vp^* (f)$ continues to make sense provided the 
following transversality condition holds,
$$
d(\vp^*(l))\not=0 \ \quad \forall  l\in  WF(f)\, .\leqno (10)
$$
where $WF(f)$ is the wave front set of $f$ ([GS]).
The next point is the construction of the generalized section of
a vector bundle $L$ on a manifold  $X$ associated to a submanifold $Z 
\sbs X$
and a symbol,
$$
\s \in C^{\ify} (Z, L\ot \vert N_Z^* \vert) \, . \leqno (11)
$$
where $N_Z$ is the normal
bundle of $Z$.
The construction is the same as that of the current of
integration on a cycle. Given $\xi \in C_c^{\ify} (X,L^* \ot
\vert T \vert)$,  the product $\s \, \xi/Z$ is a
density on $Z$, since it is  a section of $\vert T_Z \vert = \vert
T_X \vert \ot \vert N_Z^* \vert$.One can thus integrate it over $Z$. 
\noindent When $Z=X$ one has $N_Z^* = \{ 0 \}$ and $\vert
N_Z^* \vert$ has a canonical section, so that the current
associated to $\s$ is just given by (5). When $Z =$ pt is a
single point $x \in X$ a generalized section of $L$ given by
a dirac distribution at $x$ requires not only a vector $\xi_x
\in L_x$ but also a dual density, i.e. a volume multivector
$v \in \vert T_x^* \vert$.

\smallskip

\noindent Now let $\vp : X \ra Y$ with $Z$ a submanifold of
$Y$ and $\s$ as in (11).

\smallskip

\noindent Let us assume that $\vp$ is transverse to $Z$, so
that for each $x\in X$ with $y=\vp (x) \in Z$ one has
$$
\vp_* (T_x) + T_{\vp (x)} (Z) = T_y \, Y \, . \leqno (12)
$$
Let
$$
\tau_x = \{ X \in T_x \ , \ \vp_* (X) \in T_y (Z) \} \, .
\leqno (13)
$$
Then $\vp_*$ gives a canonical isomorphism,
$$
\vp_* : T_x (X) / \tau_x \sm T_y (Y) / T_y (Z) = N_y (Z) \, .
\leqno (14)
$$
And $\vp^{-1} (Z)$ is a submanifold of $X$ of
the same codimension as $Z$  with a natural isomorphism of normal 
bundles
$$
N_{\vp^{-1} (Z)} \sm \vp^* \, N_Z \, . \leqno (15)
$$
In particular, given a (generalized) $\d$-section of a bundle
$L$ with support $Z$ and symbol $\s \in C^{\ify} (Z,L\ot
\vert N_Z^* \vert)$ one has a corresponding symbol on
$\vp^{-1} (Z)$ given by
$$
\vp^* \, \s (x) = \s (\vp (x)) \in 
 (\vp^* \, L)_x \ot \vert N_x^* \vert \leqno (16)
$$
using the isomorphism (15) i.e. $N_x^* \sm N_{\vp (x)}^*$.

\smallskip

\noindent Now for any  $\d$-section associated to $Z,\s$, the wave front 
set is contained in the conormal bundle 
of the submanifold $Z$
which shows that if $\vp$ is transverse to $Z$ the pull back $\vp^* \, 
\d_{Z,\s}$ of the
distribution on $Y$ associated to $Z,\s$ makes sense ,it is equal to  
$\d_{\vp^{-1}
(Z),\vp^* (\s)}$.

\medskip

Let us now formulate the Schwartz kernel theorem. One
considers a continuous linear map, 
$$
T : C_c^{\ify} (Y) \ra C^{-\ify} (X) \, , \leqno (17)
$$
the statement is that one can write it as
$$
(T \, \xi) \, (x) = \int k(x,y) \, \xi (y) \, dy \leqno (18)
$$
where $k(x,y) \, dy$ is a generalized section,
$$
k \in C^{-\ify} (X \ts Y \ , \ {\rm pr}_Y^* (\vert T \vert))
\, . \leqno (19)
$$
Let $f : X \ra Y$ be a smooth map, and $T=f^*$ the operator
$$
(T \, \xi) \, (x) = \xi \, (f(x)) \qquad \fl \, \xi \in
C_c^{\ify} (Y) \, . \leqno (20)
$$
Let us show that the corresponding $k$ is the $\d$-section
associated to the submanifold of $X \ts Y$ given by
$$
{\rm Graph} (f) = \{ (x,f(x)) \ ; \ x \in X \} = Z \leqno
(21)
$$
and identify its symbol, $\s \in C^{\ify} (Z, {\rm pr}_Y^*
(\vert T \vert) \ot \vert N_Z^* \vert)$.

\smallskip

\noindent Given $\xi \in T_x^* (X)$, $\eta \in T_y^* (Y)$ one
has $(\xi ,\eta) \in N_Z^*$ iff it is orthogonal to $(v,f_*
\, v)$ for any $v \in T_x (X)$, i.e. $\lgl v ,\xi \rgl + \lgl
f_* \, v , \eta \rgl =0$ so that
$$
\xi = -f_*^t \, \eta \, . \leqno (22)
$$
Thus one has a canonical isomorphism $j:T_y^* (Y) \sm N_Z^*$,
$\eta \build \ra_{}^{j} (-f_*^t \, \eta , \eta)$. The
transposed $(j^{-1})^t $ is given by $(j^{-1})^t (Y) = \,
\hbox{class of} \ (0,Y)$ in $N_Z$, $\fl \, Y \in T_y (Y)$.
Thus, there is a canonical choice of symbol $\s$,
$$
\s = \vert j^{-1} \vert \in C^{\ify} (Z, {\rm pr}_Y^*
(\vert T \vert) \ot \vert N_Z^* \vert) \, . \leqno (23)
$$
We denote the corresponding $\d$-distribution by
$$
k(x,y) \, dy = \d (y-f(x)) \, dy \, . \leqno (24)
$$
One then  checks the formula,
$$
\int \d (y-f(x)) \, \xi (y) \, dy = \xi (f(x)) \qquad \fl \,
\xi \in C_c^{\ify} (Y) \, . \leqno (25)
$$

\bigskip
\noindent Let us now consider a manifold $M$ with a flow $F_t$
$$
F_t (x) = {\rm exp} (t \, v) \, x \qquad v \in C^{\ify}
(M,T_M) \leqno (26)
$$
and the corresponding map $f$,
$$
f : M \ts \Rb \ra M \ , \ f(x,t) = F_t (x) \, . \leqno (27)
$$
We apply the above discussion with $X = M \ts \Rb$, $Y = M$.
The graph of $f$ is the submanifold $Z$ of $X \ts Y$,
$$
Z = \{ (x,t,y) \ ; \ y = F_t (x)\} \, . \leqno (28)
$$
One lets $\vp$ be the diagonal map,
$$
\vp (x,t) = (x,t,x) \ , \ \vp : M \ts \Rb \ra X \ts Y
\leqno (29)
$$
and the first issue is the transversality $\quad \vp \, {\cap \!\!\! 
^{\mid}} \, Z$.

\smallskip

\noindent We thus need to consider (12) for each $(x ,t)$
such that $\vp (x,t) \in Z$, i.e. such that $x=F_t (x)$. One
looks at the image by $\vp_*$ of the tangent space $T_x \, M
\ts \Rb$ to $M \ts \Rb$ at $(x,t)$. One lets $\part_t$ be the
natural vector field on $\Rb$. The image of $(X,\lb \,
\part_t)$ is $(X,\lb \, \part_t ,X)$ for $X\in T_x \, M , \lb
\in \Rb$.  Dividing the tangent space of $M \ts \Rb \ts M$  by the image 
of $\vp_*$ one gets an
isomorphism,
$$
(X,\lb \, \part_t , Y) \ra Y-X \leqno (30)
$$
 with $T_x \, M$.The tangent space to $Z$ is $\{ (X' ,\mu \, \part_t
, (F_t)_* \, X' + \mu \, v_{F_t (x)})$; $X' \in T_x \, M ,
\mu \in \Rb \}$. Thus the transversality
condition means that every element of $T_x \, M$ is of the
form
$$
(F_t)_* \, X - X + \mu \, v_x \qquad X \in T_x \, M \ , \ \mu
\in \Rb \, . \leqno (31)
$$
One has
$$
(F_t)_* \, \mu \, v_x = \mu \, v_x \leqno (32)
$$
so that $(F_t)_*$ defines a quotient map, the Poincar\'e
return map
$$
P : T_x / \Rb \, v_x \ra T_x / \Rb \, v_x = N_x \leqno (33)
$$
and the transversality condition (31) means exactly,
$$
1-P \qquad \hbox{is invertible}. \leqno (34)
$$
Let us make this hypothesis and compute the symbol $\s$ of
the distribution,
$$
\tau = \vp^* (\d (y-F_t (x)) \, dy) \, . \leqno (35)
$$
First, as above, let $W = \vp^{-1} (Z) = \{ (x,t) \ ; \ F_t
(x)=x\}$. The codimension of $\vp^{-1} (Z)$ in $M \ts \Rb$
is the same as the codimension of $Z$ in $M\ts \Rb \ts M$ so
it is ${\rm dim} M$ which shows that $\vp^{-1} (Z)$ is
1-dimensional. If $(x,t) \in \vp^{-1} (Z)$ then $(F_s (x),t)
\in \vp^{-1} (Z)$. Thus, if we assume that $v$ does not vanish at $x$, 
the map,
$$
(x,t) \build \ra_{}^{q} t \leqno (36)
$$
is locally constant on the connected component of $\vp^{-1}
(Z)$ containing $ (x,t)$.

\smallskip

\noindent This allows to identify the transverse space to $W = 
\vp^{-1} (Z)$ as the product,
$$
N_{x,t}^W \sm N_x \ts \Rb \leqno (37)
$$
where to $(X,\lb \, \part_t) \in T_{x,t} (M\ts \Rb)$ we
associate the pair $(\wt X ,\lb)$ given by the class of $X$ in
$N_x = T_x / \Rb \, v_x $ and $\lb \in \Rb$.

\smallskip

\noindent The symbol $\s$ of the distribution (35) is a
smooth section of $\vert N^{W*} \vert$ tensored by the pull
back $\vp^* (L)$ where $L = {\rm pr}_Y^* \, \vert T_M \vert$,
and one has
$$
\vp^* (L) \sm \vert p^* \, T_M \vert \leqno (38)
$$
where
$$
p(x,t) = x \qquad \fl \, (x,t) \in M \ts \Rb \, . \leqno (39)
$$
To compute $\s$ one needs the isomorphism,
$$
N_{(x,t)}^W \build \ra_{}^{\vp_*} T_{\vp (x,t)} (M \ts \Rb
\ts M) / T_{\vp (x,t)} (Z) = N^Z \, . \leqno (40)
$$

\smallskip

\noindent  The map $\vp_* : N_{x,t}^W \ra N^Z$ is given by
$$
\vp_* (X,\lb \, \part_t) = (1-(F_t)_*) \, X - \lb \, v \qquad
X \in N_x \ , \ \lb \in \Rb \leqno (41)
$$
and the symbol $\s$ is just
$$
\s = \vert \vp_*^{-1} \vert \in \vert p^* \, T_M \vert \ot
\vert N^{W*} \vert \, . \leqno (42)
$$
This makes sense since $\vp_*^{-1} : p^* \, T_M \ra N^W$.

\smallskip

\noindent Let us now consider the second projection,
$$
q(x,t) = t \in \Rb \leqno (43)
$$
and compute the pushforward $q_* (\tau)$ of the distribution
$\tau$.

\smallskip

\noindent By construction $\d (y-F_t (x)) \, dy$ is a
generalized section of ${\rm pr}_Y^* \, \vert T \vert$, so
that $\tau$ is a generalized section of $p^* \, \vert T \vert
= \vp^* \, {\rm pr}_Y^* \, \vert T \vert$.

\smallskip

\noindent Thus $q_* (\tau)$ is a generalized function.

\smallskip

\noindent We first look at the contribution of a periodic
orbit, the corresponding part of $\vp^{-1} (Z)$ is of the form,
$$
\vp^{-1} (Z) = V \ts \G \sbs M \ts \Rb \leqno (44)
$$
where $\G$ is a discrete cocompact subgroup of $\Rb$, while
$V \sbs M$ is a one dimensional compact submanifold of $M$.

\smallskip

\noindent To compute $q_* (\tau)$, we let $h(t) \, \vert dt
\vert$ be a 1-density on $\Rb$ and pull it back by $q$ as the
section on $M \ts \Rb$ of the bundle $q^* \, \vert T \vert$,
$$
\xi (x,t) = h(t) \, \vert dt \vert \, . \leqno (45)
$$
We now need to compute $\int_{\vp{-1} (Z)} \, \xi \, \s$. We
can look at the contribution of each component: $V \ts \{
T \}$, $T \in \G$.

\smallskip

\noindent One gets

$$
 T^{\#}  {1 \over \vert 1-P_T
\vert} \ h(T) \, . \leqno (46)
$$
Where $T^{\#}$ is the length of the primitive orbit or
equivalently the covolume of $\G$ in $\Rb$ for the Haar measure
$\vert dt \vert$. We can thus write the contributions of the
periodic orbits as
$$
\sum_{\g_p} \sum_{\G} {\rm Covol} (\G) \, {1 \over \vert 1-P_T
\vert} \ h(T) \, . \leqno (47)
$$
Where the test function $h$ vanishes at $0$.\smallskip

\noindent
The next case to consider is when the vector field
$v_x$ has an isolated $0$, $v_{x_0} =0$. In that case, the
transversality condition (31) becomes
$$
1 -(F_t)_*  \quad \hbox{invertible (at} \ x_0) \, . \leqno (48)
$$
One has $F_t (x_0) = x_0 \, $ for all $t\in \Rb$  and now the
relevant component of $\vp^{-1} (Z)$ is $\{ x_0 \} \ts \Rb$.
The transverse space $N^W$ is identified with $T_x$ and the
map $\vp_* : N^W \sm N^Z$ is given by:
$$
\vp_* = 1-(F_t)_* \, . \leqno (49)
$$
Thus the symbol $\s$ is the scalar function $\vert 1-(F_t)_*
\vert^{-1}$ . The
generalized section $q_* \, \vp^* (\d (y-F_t (x)) \, dy)$ is the
function, $t \ra \vert 1- (F_t)_* \vert^{-1}$. We can thus
write the contribution of the zeros of the flow as,
$$
\sum_{\sevenrm zeros} \int {h(t) \over \vert 1-(F_t)_*
\vert} \ dt \leqno (50)
$$
where $h$ is a test function vanishing at 0.

\smallskip

\noindent We can thus collect the contributions 47 and 50 as
$$
\sum_{\sevenrm \g} \int_{I_{\g}} {h(u) \over \vert 1-(F_u)_*
\vert} \ d^*u \leqno (51)
$$
where $h$ is as above, $I_{\g}$ is the isotropy group of the periodic 
orbit $\g$, the haar measure $d^*u$ on $I_{\g}$ is normalised so that 
the covolume of $I_{\g}$ is equal to one and we still write $(F_u)_*$ 
for its restriction to the transverse space of $\g$.
 
\smallskip

\noindent {\bf Bibliography.}

\bigskip

\item{[AB]} M.F. Atiyah and R. Bott, A Lefchetz fixed point formula for 
elliptic complexes: I, {\it Annals of Math,} {\bf 86} (1967), 374-407.

\item{[B]} M. Berry, Riemann's zeta function: a model of
quantum chaos, {\it Lecture Notes in Physics,} {\bf 263},
Springer (1986).
\item{[Bg]} A. Beurling,   A closure problem related to the  ( zeta 
function, {\it Proc. Nat. Ac. Sci.} {\bf 41} (1955),  312-314.
\item{[B-C]} J.-B. Bost and A. Connes, Hecke Algebras, Type
III factors and phase transitions with spontaneous symmetry
breaking in number theory, {\it Selecta Mathematica, New
Series } {\bf 1}, No.3 (1995), 411-457.
\item{[BG]} O. Bohigas and M. Giannoni, Chaotic motion and random matrix 
theories, {\it Lecture Notes in Physics,} {\bf 209} (1984), 1-99.
\item{[BK]} M. Berry and J. Keating, $H=qp$ and the Riemann zeros, 
'Supersymmetry and Trace Formulae: Chaos and Disorder', edited by
J.P. Keating, D.E. Khmelnitskii and I.V. Lerner (Plenum Press).
\item{[Br]} F. Bruhat,  Distributions sur un groupe localement compact 
et applications \'a l'\'etude des repr\'esentations des groupes 
$p$-adiques.{\it Bull. Soc. Math. france.} {\bf 89} (1961),  43-75.
\item{[C]} A. Connes, Noncommutative Geometry, Academic
Press (1994).
\item{[Co]} A. Connes, Formule de trace en G\'eom\'etrie non commutative 
et hypoth\`ese de Riemann, {\it C.R. Acad. Sci. Paris Ser. A-B} (1996)
\item{[D]} C. Deninger, Local $L$-factors of motives and
regularised determinants, {\it Invent. Math.,} {\bf 107}
(1992), 135-150.

\item{[G]} D. Goldfeld, A spectral interpretation of Weil's
explicit formula, {\it Lecture Notes in Math.,} {\bf 1593},
Springer Verlag (1994), 135-152.
\item{[GS]} V. Guillemin and S. Sternberg, Geometric asymptotics, {\it 
Math. Surveys}, {\bf 14}, {\it  Amer. Math. Soc., Providence, R.I.} 
(1977)
\item{[Gu]} V. Guillemin, Lectures on spectral theory of
elliptic operators, {\it Duke Math. J.,} {\bf 44}, No.3
(1977), 485-517.

\item{[H]} S. Haran, Riesz potentials and explicit sums in
arithmetic, {\it Invent. Math.,} {\bf 101} (1990), 697-703.
\item{[J]} B. Julia, Statistical theory of numbers, {\it Number Theory 
and Physics, Springer Proceedings in Physics,} {\bf 47} (1990).
\item{[K]} M. Kac, Statistical Independence in Probability, {\it 
Analysis and Number Theory, Carus Math. Monographs } {\bf 18} (1959).
\item{[KS]} N. Katz and P. Sarnak, Random matrices, Frobenius 
eigenvalues and Monodromy, (1996) , Book, to appear.
\item{[KS]} N. Katz and P. Sarnak, Zeros of zeta functions, their 
spacings and spectral nature, (1997), to appear.
\item{[LPS1]} D. Slepian and H. Pollak, Prolate spheroidal wave 
functions, Fourier analysis and uncertainty I, {\it Bell Syst. Tech. J.} 
{\bf 40} (1961).
\item{[LPS2]} H.J. Landau and H. Pollak, Prolate spheroidal wave 
functions, Fourier analysis and uncertainty II, {\it Bell Syst. Tech. 
J.} {\bf 40} (1961).
\item{[LPS3]} H.J. Landau and H. Pollak, Prolate spheroidal wave 
functions, Fourier analysis and uncertainty III, {\it Bell Syst. Tech. 
J.} {\bf 41} (1962).
\item{[M]} H. Montgomery, The pair correlation of zeros of
the zeta function, {\it Analytic Number Theory,} AMS (1973).
\item{[Me]} M.L. Mehta, Random matrices, Academic Press,(1991).
\item {[O]} A. Odlyzko, On the distribution of spacings between zeros of 
zeta functions, {\it Math. Comp.} {\bf 48} (1987), 273-308.  
\item {[P]} G. P\'olya, Collected Papers, Cambridge, M.I.T. Press 
(1974).
\item{[Pat]} S. Patterson, An introduction to the theory of the Riemann 
zeta
function, {\it Cambridge Studies in advanced mathematics}, {\bf 14}
Cambridge University Press (1988).
\item{[R]} B. Riemann, Mathematical Werke, Dover, New York (1953). 
\item{[S]} E.Seiler, Gauge Theories as a problem of constructive Quantum 
Field Theory and Statistical Mechanics, Lecture Notes in Physics {\bf 
159}
 Springer (1982).
\item{[Se]} A.Selberg, {\it Collected papers},
 Springer (1989).
\item{[W1]} A. Weil, Basic Number Theory, Springer, New
York (1974).

\item{[W2]} A. Weil, Fonctions z\^eta et distributions,
{\it S\'eminaire Bourbaki}, {\bf 312}, (1966).

\item{[W3]} A. Weil, Sur les formules explicites de la
th\'eorie des nombres, {\it Izv. Mat. Nauk.,} (Ser. Mat.)
{\bf 36}, 3-18.
\item{[W4]} A. Weil, Sur la th\'eorie du corps de classes,{\it J. Math. 
Soc. Japan}, {\bf 3}, (1951).

\item{[W5]} A. Weil, Sur certains groupes d'operateurs unitaires,{\it 
Acta Math. }, {\bf 111}, (1964).

\item{[Z]} D. Zagier, Eisenstein series and the Riemannian zeta 
function, {\it Automorphic Forms, Representation Theory and Arithmetic,} 
Tata, Bombay (1979), 275-301.

\bye